\documentclass[sort&compress,review]{elsarticle} %times

\usepackage[letterpaper,hmargin=1.75cm,tmargin=2.0cm,bmargin=2.25cm]{geometry}

\usepackage{xcolor}
\usepackage[colorlinks]{hyperref}

\hypersetup{colorlinks=true} %,linkbordercolor=red,linkcolor=green}

\usepackage{lineno}
\modulolinenumbers[5]

\usepackage{amsfonts}           % Better font

\usepackage[T1]{fontenc}        % Allows for the reader to copy "á"  %XeLaTex
\usepackage[utf8]{inputenc}     % Allow input like "á"    %XeLaTex
\usepackage[nomath]{lmodern}

\usepackage{soul}               % Strikethrough text with \st
\usepackage[super]{nth}         % Use nth(1) to get 1^st
\usepackage{amsmath}
\usepackage{amssymb}  % various useful mathematical symbols
\usepackage{amsthm}   % extended theorem environments

\usepackage{cancel}   % To be able to cancel-out terms in equations
\usepackage[normalem]{ulem}  %  Strikes out text with \sout
\usepackage{accents}  % Miscellaneous tools for mathematical accents (J. Bezos)
\usepackage{bm}       % Makes its ar­gu­ment bold

\usepackage{empheq}   % for a lot of equations instead of "cases"
\usepackage{mathtools} % defines curly brackets for sets

\usepackage{mleftright}
\usepackage{xfrac}

\usepackage{graphicx}      % complicated commands
\graphicspath{{./Figures/}}

\usepackage[caption=false]{subfig}       % Add figures side by side
\usepackage{float}
\makeatletter
\newcommand\fs@boxedtop
  {\fs@boxed
   \def\@fs@mid{\vspace\abovecaptionskip\relax}%
   \let\@fs@iftopcapt\iftrue
  }
\makeatother
\floatstyle{boxedtop}
\floatname{framedbox}{Box}
\newfloat{framedbox}{tbp}{lob} %[section]

\usepackage{caption}
\usepackage{tikz}               % To draw on latex
\usetikzlibrary{decorations.pathmorphing,decorations.markings,patterns}
\usetikzlibrary{arrows}         % For Tikz Arrows
\usetikzlibrary{positioning}    % For positioning nodes

\usepackage{multicol}           % Several columns
\usepackage{multirow}			% Several rows in table
\usepackage{booktabs}
\usepackage{tcolorbox}

\usepackage[capitalize]{cleveref}

\usepackage{python}
\usepackage{pythontex}
\usepackage{verbatim}

\journal{Elsevier (LA-UR-21-30031) }

\newtheorem{remark}{Remark}[section]

\newcommand\Tstrut{\rule{0pt}{2.6ex}}         % = `top' strut
\newcommand\Bstrut{\rule[-0.9ex]{0pt}{0pt}}   % = `bottom' strut

\newcommand*{\dt}[1]{%
  \overset{\textrm{\Large .}}{#1}}

\newcommand*{\ddt}[1]{%
  \overset{\textrm{\Large .\hspace{-0.2ex}.}}{#1}}

\newcommand{\ud}{\,\mathrm{d}}

\newcommand{\Div}{{\text{Div}}}
\newcommand{\Grad}{{\nabla}}
\newcommand{\Lapl}{{\Delta}}

\newcommand{\vect}[1]{\boldsymbol{#1}}
\newcommand{\tens}[1]{\boldsymbol{#1}}
\newcommand{\tenf}[1]{\pmb{\mathbb{#1}}}
\newcommand{\partder}[2]{\frac{\partial #1}{\partial #2}}

\newcommand{\fullder}[2]{\frac{\ud #1}{\ud #2}}

\newcommand{\avrg}[1]{ \left\{ #1 \right\} }
\newcommand{\jump}[1]{ \left[\left[ #1 \right]\right] }

\newcommand{\innerprod}[3] { \left( #1, #2 \right)_{#3}}

\newcommand{\transpose}[0]{^\text{\textbf{T}}}

\newcommand\thickbar[1]{\accentset{\rule{.4em}{.8pt}}{#1}}

\newcommand{\RefDom}{{\Omega_0}}
\newcommand{\ArbDom}{{\mathcal{U}_0}}
\newcommand{\ElRefDom}{{\Omega_{0,e}}}
\newcommand{\AppRefDom}{{\Omega^{h}_{0}}}

\newcommand{\TildeRefDom}{{\widetilde{\Omega}_0}}

\newcommand{\ExtBoun}{{\Gamma}}
\newcommand{\IntBoun}{{\widetilde{\Gamma}}}  % Includes part of exterior boundary
\newcommand{\IntOnlyBoun}{{\widehat{\Gamma}}}

\newcommand{\CurDom}{{\Omega_t}}
\newcommand{\MicroDens}{{\rho_\mu}}
\newcommand{\MicroDamp}{{\zeta_\mu}}
\newcommand{\DeformGrad}{{\tens{F}}}

\newcommand{\PF}{{d}}
\newcommand{\PFgrad}{{\nabla \PF}}
\newcommand{\PFlap}{{\Delta \PF}}

\newcommand{\virtPF}{{\delta \PF}}
\newcommand{\virtPFgrad}{{\nabla \delta \PF}}
\newcommand{\virtPFlap}{{\Delta \delta \PF}}

\newcommand{\BodyForceRef}{{\vect{b}^*}}
\newcommand{\SurfTractionRef}{{\vect{t}^*}}

\newcommand{\DispRef}{{\vect{u}}}
\newcommand{\CompInt}[3]{{\int_{#1} \!\! {#2} \ud {#3}}}
\newcommand{\SpaceTimeInt}[4]{{\int_{t_1}^{t_2} \!\! \int_{#1} \!\! \left[{#2}\right] {#3} \ud {#4} \ud \tau  }}

\newcommand{\lzr}{{\ell_0}}
\newcommand{\motion}{{\vect{\varphi}}}

\newcommand{\TotalEnerDenFun}{{\hat{W}}}
\newcommand{\DegElastEnerDen}{{\hat{W}_e}}
\newcommand{\SurfEnerDenFun}{{\hat{W}_f}}

\newcommand{\NormVecRef}{{\vect{N}}}

\newcommand{\trace}[1]{ tr \left( #1 \right) }

\newcommand{\disp}{\vect{u}}
\newcommand{\velo}{\dt{\vect{u}}}
\newcommand{\accel}{\ddt{\vect{u}}}

\let\oldequation\equation
\let\oldendequation\endequation

\renewenvironment{equation}
  {\linenomathNonumbers\oldequation}
  {\oldendequation\endlinenomath}

\begin{document}

\begin{frontmatter}

\title{A fourth-order phase-field fracture model: Formulation and numerical solution using a continuous/discontinuous Galerkin method}



\author[LANL]{Lampros Svolos\corref{cor1}}
\ead{lsvolos@lanl.gov}
\cortext[cor1]{Corresponding author}

\author[LANL]{Hashem M. Mourad} 

\author[LANLTV]{Gianmarco Manzini} 

\author[UniMichA,UniMichB,UniMichC]{Krishna Garikipati}

%\ead{hmourad@lanl.gov}

%\ead{lsvolos@lanl.gov}

%address[mysecondaryaddress]{360 Park Avenue South, New York}

\address[LANL]{Fluid Dynamics and Solid Mechanics, T-3, Theoretical Division, Los Alamos National Laboratory, Los Alamos, NM, USA}

\address[LANLTV]{Applied Mathematics and Plasma Physics, T-5, Theoretical Division,  Los Alamos National Laboratory, Los Alamos, NM, USA}

\address[UniMichA]{Department of Mechanical Engineering, University of Michigan, Ann Arbor, MI, USA}
\address[UniMichB]{Department of Mathematics, University of Michigan, Ann Arbor, MI, USA}
\address[UniMichC]{Michigan Institute for Computational Discovery \& Engineering, University of Michigan, Ann Arbor, MI, USA}

\begin{abstract}
Modeling crack initiation and propagation in brittle materials is of great importance to be able to predict sudden loss of load-carrying capacity and prevent catastrophic failure under severe dynamic loading conditions. Second-order phase-field fracture models have gained wide adoption given their ability to capture the formation of complex fracture patterns, e.g.\ via crack merging and branching, and their suitability for implementation within the context of the conventional finite element method. Higher-order phase-field models have also been proposed to increase the regularity of the exact solution and thus increase the spatial convergence rate of its numerical approximation. However, they require special numerical techniques to enforce the necessary continuity of the phase field solution. In this paper, we derive a fourth-order phase-field model of fracture in two independent ways; namely, from Hamilton's principle and from a higher-order micromechanics-based approach. The latter approach is novel, and provides a physical interpretation of the higher-order terms in the model. In addition, we propose a continuous/discontinuous Galerkin (C/DG) method for use in computing the approximate phase-field solution. This method employs Lagrange polynomial shape functions to guarantee $C^0$-continuity of the solution at inter-element boundaries, and enforces the required $C^1$ regularity with the aid of additional variational and interior penalty terms in the weak form. The phase-field equation is coupled with the momentum balance equation to model dynamic fracture problems in hyper-elastic materials. Two benchmark problems are presented to compare the numerical behavior of the C/DG method with mixed finite element methods.
\end{abstract}

%\MSC[2010] 00-01\sep  99-00

\end{frontmatter}

%\linenumbers

\section{Introduction}
\label{sec:Intro}

Applications in a variety of important sectors, such as aerospace and energy, require materials and structures to withstand extreme service conditions. In many such applications, it may be required to design against catastrophic failure under accident scenarios, and the material behavior in the post-failure regime may also be of interest. These needs motivate efforts aimed at the development of modeling and simulation methodologies that are capable of predicting material failure in an accurate and computationally efficient manner, and despite much work conducted over the past few decades, and a voluminous literature covering this area of research, predictive modeling of crack initiation and propagation in materials and structures remains one of the most significant challenges in solid mechanics \cite{wu_chapter_2020}.

% Physics (and modeling)

Fracture is an important mode of material failure which is frequently encountered in practice. It is often regarded as a process consisting of a nucleation/initiation phase followed by a crack propagation phase \cite{besson2010continuum,rabczuk2013computational}, where pre-existing micro-scale defects or inclusions typically serve as nucleation sites. From a modeling standpoint, this necessitates the development and use of a fracture criterion combined with a crack propagation model that must both incorporate accurate mechanistic information in order to achieve the required predictive ability. Moreover, while the physical mechanisms that culminate in fracture may differ from one material to another, complicating the modeling task, fracture invariably involves the introduction of (discrete) discontinuities into a medium that is initially continuous, thus presenting a separate challenge from a computational perspective.

% Discrete Models of material failure

Methodologies for the computational treatment of fracture include cohesive zone modeling \cite{xu1994numerical,CAMACHO19962899,park2011cohesive} where discontinuities are typically introduced into the mesh between neighboring elements, potentially giving rise to mesh dependence and other non-physical artifacts \cite{rimoli2012mesh, versino2015thermodynamically}. In other computational treatments, such as extended finite element methods (XFEM) \cite{daux2000arbitrary,sukumar2000extended,fries2010extended} and their embedded finite element (EFEM) counterparts \cite{armero1996analysis,linder2007finite,mourad2017modeling,jin2019three}, discontinuities are introduced within individual elements in the mesh. This requires explicit tracking of these discontinuities as they propagate through the computational domain, which can be a difficult task in the case of complex crack topologies, especially in three-dimensional problems; see \cite{jin2019three,jin2019comparative} and references therein. Comparative studies covering a number of methods for the representation of fracture in a computational setting can be found in the literature; e.g., see \cite{jirasek2000comparative,song2008comparative}.

% Phase-field models (review)
Phase field (PF) modeling of fracture overcomes the aforementioned difficulties in an elegant manner. The PF approach to fracture was originally proposed by \citet{bourdin_numerical_2000} as a generalization of the variational formulation of \citet{francfort_revisiting_1998}, which in turn regards fracture as a competition between bulk elastic energy and crack surface energy \cite{griffith1921vi}. In this framework, the solution to the fracture problem is obtained simply by minimizing the total potential energy of the solid body under consideration. In addition, the approach makes use of an auxiliary field, denoted herein by $\PF$, to obtain a regularized\footnote{In fact, PF models have been shown to be related to gradient regularized continuum damage models \cite{de2016gradient}.} representation of crack surfaces, obviating the need for injecting discontinuities into the kinematic solution fields.

The PF approach to modeling fracture received significant attention over the past 10 years or so, due to its  ability to capture complicated fracture patterns, such as crack merging and branching, in an elegant manner. In the seminal work of \citet{miehe2010thermodynamically}, a framework based on continuum mechanics and thermodynamic arguments was introduced to model brittle fracture, and was further developed to treat dynamic brittle failure \cite{borden2012phase,schluter2014phase,ambati2015review,mandal2020evaluation}. Extensions of the PF approach to enable the treatment of cohesive fracture were presented in \cite{bourdin2008variational,verhoosel2013phase,vignollet2014phase,geelen2019phase}, and thermodynamically-consistent frameworks were developed in \cite{miehe_phase_2015,borden_phase-field_2016,hu2021variational} for the PF modeling of fracture in ductile solids, taking into account the coupling between plastic deformation and damage. The latter developments are notable since, in contrast to the brittle case, no variational theory exists in the case of ductile fracture \cite{ambati2015phase}. Shear localization, a precursor to ductile fracture, can also be captured using PF models of this (thermodynamically-consistent) type coupled with elastic-viscoplastic constitutive relations, as shown in \cite{mcauliffe_unified_2015,zhang2021phase}. Recent works also study the coupling between damage and heat transfer in brittle \cite{badnava2018h,yan2021new} and ductile materials \cite{dittmann2020phase,svolos2020thermal,svolos2021anisotropic}.

% Higher order PF models - goals
One of the oft-cited drawbacks of PF fracture models is their high computational cost, due in large part to the need for a highly refined mesh to accurately resolve the gradients that develop in the PF solution, $\PF$, in the vicinity of a crack~\cite{wu_chapter_2020}. Adoption of the popular second-order PF fracture model exacerbates this issue, as it leads to the development of a cusp in the solution field at the crack surface, which has a negative impact on the rate of convergence of the solution with mesh refinement (i.e.\ the spatial convergence rate). This lack of regularity of the solution obtained using the second-order theory has motivated \citet{borden_higher-order_2014} to develop their fourth-order PF fracture model. Adopting this higher-order PF theory involves a trade-off, however, as it necessitates the use of specialized computational techniques that ensure $C^1$-continuity of $\PF$ in order to attain convergence. The conventional finite element method (FEM) with a classical ($C^0$-continuous) Lagrange polynomial basis is not an adequate choice in this context, for instance. Its use to treat higher-order PF problems (including those governed by the Cahn--Hilliard equation \cite{wells_discontinuous_2006}) results in discontinuities in the gradient of the phase field, and thus constitutes a variational crime \cite{strang_variational_1972} since the weak form in such problems includes at least second-order spatial derivatives (Laplacians and Hessians) thereof. 

Different categories of computational methods have been developed over the years specifically for the treatment of such problems (with $C^1$-continuous or higher-regularity solutions). These include the FEM with Hermite polynomial basis functions \cite{stogner2008approximation} or B-spline basis functions (i.e.\ isogeometric analysis) \cite{cottrell2006isogeometric,bartezzaghi2015isogeometric}, virtual element methods \cite{antonietti2016c,antonietti2020conforming} and discontinuous Galerkin (DG) methods \cite{georgoulis2009discontinuous}, in addition to mixed FEM and continuous/discontinuous Galerkin (C/DG) methods, which are discussed in detail in the present work.

% Mixed
Mixed finite element formulations for fourth-order phase field problems are founded on introducing an additional scalar field. It is cast in a multi-field variational formulation, whose Euler-Lagrange equations weakly relate this added field to the chemical potential associated with the derivative of the free energy with respect to the auxiliary field, $d$. As is well-understood, however, the stability of mixed methods is conditional, depending on the discretization of the different fields. The sufficient criteria for stability are the celebrated Ladyzhenskaya-Babu\v{s}ka-Brezzi (LBB) stability conditions \cite{ladyzhenskaya1969mathematical,babuvska1973finite,brezzi1974existence,auricchio2017mixed}. This can prove a drawback that detracts from their ease of implementation for fourth-order phase field problems. While advances have steadily been made over the years \cite{malkus1978mixed,franca1988two,arnold1990mixed,boffi2013mixed} in the formulation of mixed methods, this fundamental obstacle remains.

% DG methods
Discontinuous Galerkin methods were initially introduced in \cite{reed1973triangular} for hyperbolic problems (specifically, the steady-state neutron transport equation). Since then, the development of DG methods has been an active research field, resulting in different formulations for partial differential equations of different types/orders. A historical overview of this topic can be found in \cite{arnold2002unified,hesthaven2007nodal}. In interior penalty discontinuous Galerkin (IPDG) methods, the basis functions do not satisfy the continuity requirements by design. Instead, inter-element ($C^0$ and/or $C^1$) continuity requirements are enforced weakly with the aid of extra penalty terms added to the variational formulation. The origin of these methods for high-order equations can be traced back to \cite{babuvska1973nonconforming,baker1977finite} where interior penalties are used to weakly impose $C^1$ inter-element continuity for fourth-order problems. C/DG methods are related to IPDG methods (with the added advantage of satisfying variational consistency \cite{wells2007c0}), and present a natural approach to high-order partial differential equations \cite{engel2002continuous,wells2004discontinuous,molari2006discontinuous,wells_discontinuous_2006,wells2007c0,hansbo2011posteriori,larsson2017continuous}. Moreover, while the most frequent criticism of DG methods relates to the proliferation of degrees of freedom especially in 3D \cite{HUGHES20062761}, C/DG methods elegantly counter this drawback by employing a $C^0$-continuous basis akin to the conventional FEM, and accounting for gradient discontinuities using the same degrees of freedom of the $C^0$ mesh in a manner that preserves variational consistency. Like their fully discontinuous counterparts, C/DG formulations also do need stabilization. However, unlike the fundamental restriction on choice of spaces that is delineated by the LBB condition, stability of C/DG methods is ensured by interior penalty terms \cite{engel2002continuous,wells2004discontinuous,molari2006discontinuous,wells2006discontinuous,hartmann2008optimal} or lifting operators \cite{georgoulis2009discontinuous,carstensen2009unifying}.

% Goal of this work
In the present paper, we develop a modeling and simulation framework for dynamic brittle fracture problems. To this end, we derive a fourth-order PF model of fracture in two independent ways; namely, from Hamilton's principle and from a higher-order micromechanics-based approach. The latter approach is novel and provides a physical interpretation of the higher-order terms in the PF model. In addition, we propose a C/DG method for use in computing the approximate PF solution. This method employs Lagrange polynomial shape functions to guarantee $C^0$-continuity of the solution at inter-element boundaries, and enforces the required $C^1$ regularity with the aid of additional variational and interior penalty terms in the weak form. The PF equation is coupled with the momentum balance equation to model dynamic fracture problems in hyper-elastic solids. Making use of the conventional FEM to solve the momentum balance equation, we compare different methods for the treatment of the fourth-order PF equation. Specifically, we present a comparison between the proposed C/DG method and two mixed finite element methods which, also, have not been applied previously in conjunction with a fourth-order PF fracture model. For the coupled field problem, we adopt a partitioned/staggered solution approach where each of the two governing equations is solved individually. 

% Structure of paper
The remainder of this paper is organized as follows. Phase-field modeling of brittle fracture is discussed in Sect.~\ref{sec:Formulation}. In this section, we present the fourth-order PF model equations along with the hyperelastic constitutive equations used to describe the undamaged material behavior. The proposed computational formulations are described in detail in Sect.~\ref{sec:Computational}.  This includes three spatial discretization schemes, based on mixed finite element and continuous/discontinuous Galerkin methods, in addition to the staggered scheme used to step the solution of the coupled equations forward in time. Numerical simulations involving two standard benchmark problems are carried out and presented in Sect.~\ref{sec:Results}, to assess the relative merits and compare the performance characteristics of the three different computational formulations in the context of dynamic brittle fracture problems. We also present sensitivity analyses examining the influence of the penalty parameter used on the C/DG solution. Finally, a summary and concluding remarks are presented in Sect.~\ref{sec:Conclusions}.

%\clearpage
\section{Formulation: a fourth-order phase-field fracture model in finite elasticity }
\label{sec:Formulation}

\subsection{Notation and primary unknown fields}
In this section, we introduce notation and define primary fields which are necessary for the description of a fourth-order phase-field fracture model in finite hyper-elasticity.

\subsubsection{Phase-field regularization of discrete crack surfaces}

\begin{figure}[t] 
\centering
\subfloat[\label{fig:DiscreteCrackSmall}]{%
\includegraphics[width=0.335\textwidth]{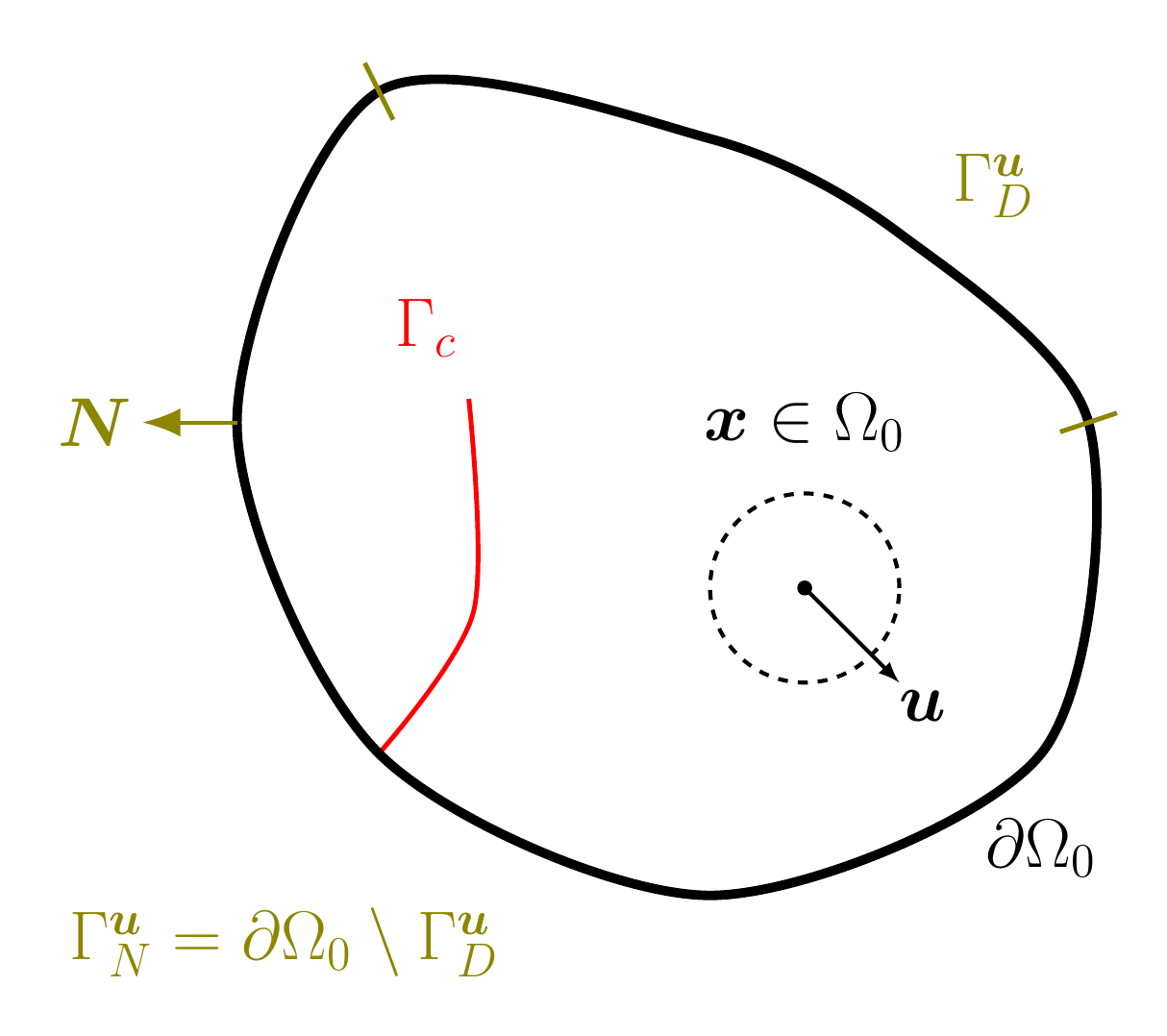}  %0.29
}
\hspace{1.0cm}
\subfloat[\label{fig:DiffCrackSmall_BC}]{%
\includegraphics[width=0.41\textwidth]{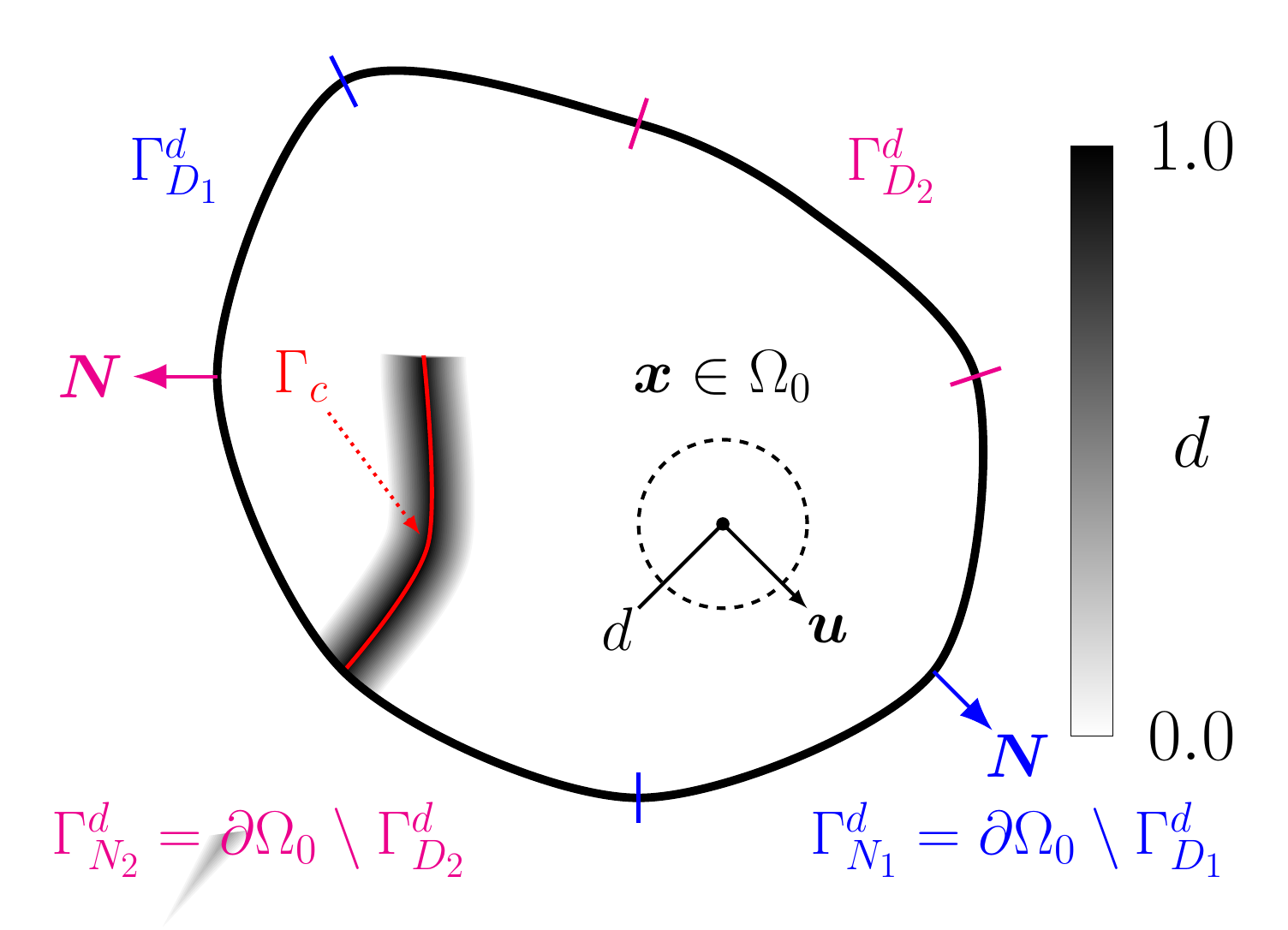}
}
\captionsetup{width=.90\textwidth}
\caption{ Unknown fields and associated boundary conditions: (a) deformable body with a discrete crack of length $\Gamma_c$ (b) approximation of the crack by introducing a phase field $\PF$ which distinguishes between intact ($d=0$) and fully damaged ($d=1$) material states}
\label{fig:FreeBodyWithCrackAndBC}
\end{figure}

We consider a deformable body, represented by $\RefDom$ in the reference configuration, with a boundary denoted by $\partial \RefDom$ as shown in \cref{fig:DiscreteCrackSmall}. Under given loading conditions, over a time interval of interest $\mathcal{T}$, the body deforms and may develop cracks. This state is described by the unknown displacement field $\vect{u}$ and \emph{discrete} crack, denoted by $\Gamma_c$. This setting forms a free-discontinuity problem since the unknown (displacement) field is discontinuous across unknown sets (crack surfaces).

In fracture mechanics, the Griffith criterion \cite{griffith1921vi} for crack evolution can be reformulated in a variational setting as presented in \cite{bourdin2008variational}. In this context, the \emph{discrete} crack can be approximated by a phase field $\PF$, defined as follows
\begin{equation}
\left. \begin{aligned}
\PF \text{ : } & \RefDom \times \mathcal{T} \rightarrow [0,1] \\
		  & (\vect{X},t) \mapsto d(\vect{X},t)  \\
\end{aligned} \right\}, 
\end{equation}
where $\PF=0$ characterizes intact material and $\PF=1$ denotes the fully damaged state at a material point $\vect{X}$ in the reference configuration $\RefDom$,  as depicted in \cref{fig:DiffCrackSmall_BC}. 

Using the phase field $\PF$, the crack length can be approximated by introducing the following regularization
\begin{equation} \label{eq:GammaApproximationDomainIntegral}
\Gamma_c \approx \Gamma_{\lzr}= \int_\RefDom \gamma_{\lzr}(\PF, \PFgrad, \PFlap) \ud V,
\end{equation}
where $\gamma_{\lzr}$ is the so-called crack surface density function per unit  reference volume of the solid and $\lzr$ is a length-scale parameter with which the formulation is regularized \cite{miehe_phase_2015}. 

%\lsc{Here you can write the 1-D minimization problem.}

In the literature, different forms of $\gamma_{\lzr}$ have been introduced depending on the phase-field solutions of the following one-dimensional minimization problem:
\begin{equation}
\left.
\begin{aligned}
\min_{d} \quad & \Gamma_{\lzr} \\
\textrm{s.t.} \quad  & d(0)=1 \\
\textrm{and} \quad &  d(X) \rightarrow 0 \textrm{ as } |X| \rightarrow \infty \\
\end{aligned}
\right\} .
\end{equation}
These phase-field solutions are devised to represent `diffusive' cracks with steep gradients. Assuming specific crack density functions, it is possible to derive the aforementioned solutions in 1D problems and generalize the model to higher dimensions.

\begin{figure}[H] 
\centering
\includegraphics[width=0.45\textwidth]{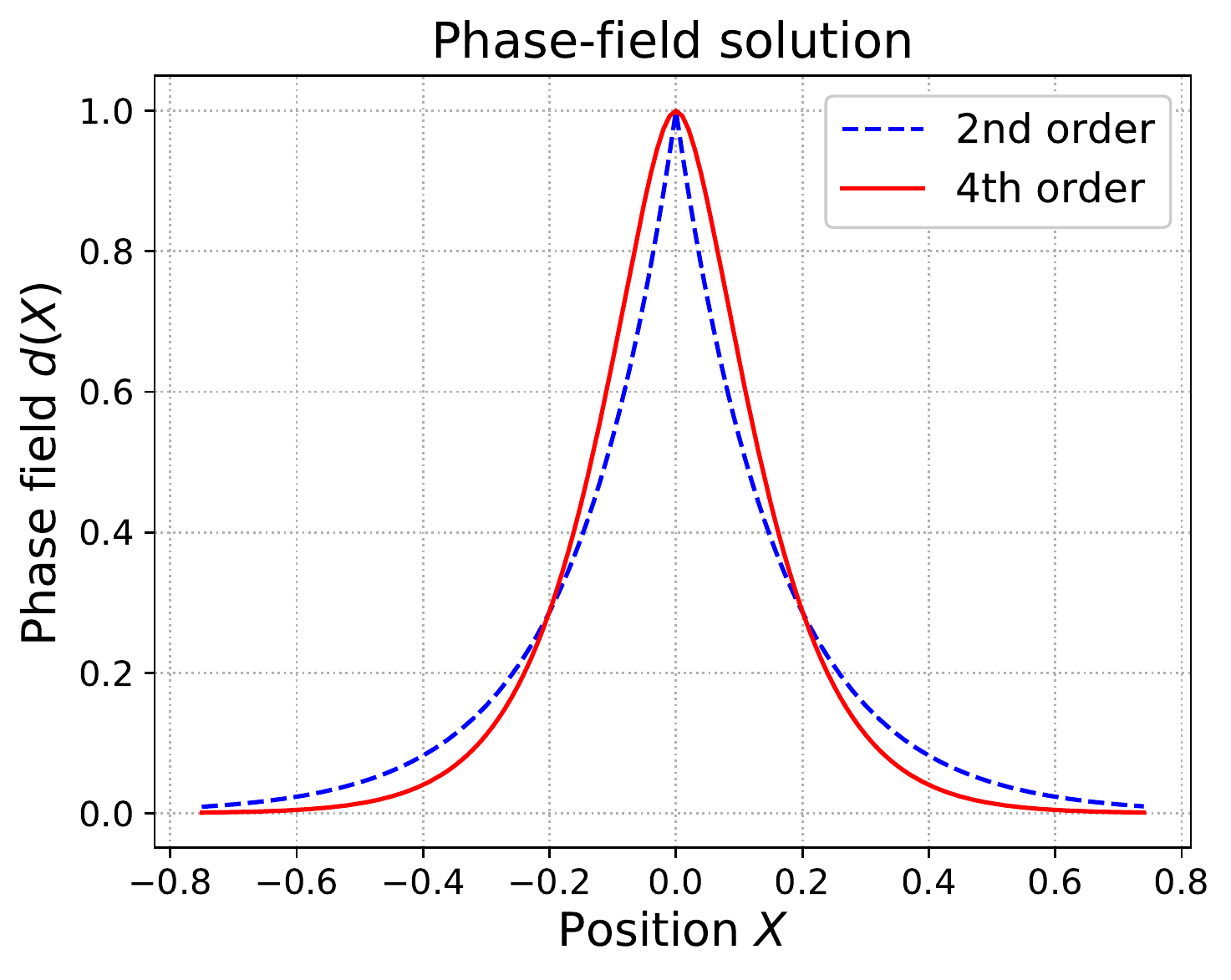}
\captionsetup{width=.90\textwidth}
\caption{ Comparison between the second-order and fourth-order solutions of the one-dimensional phase-field equation.}
\label{fig:PhaseFieldSolutionProfiles1D}
\end{figure}

For instance, in one dimension, the phase-field solution is given by the $C^0$ exponential function, $\PF(X)=\exp\left(-|X|/2 \lzr\right)$, when the crack density function reads
\begin{equation}
\gamma_{\lzr}(\PF, \PFgrad) = \frac{1}{4\lzr} \left[ d^2 +  4 \lzr^2 |\PFgrad|^2 \right] .
\end{equation}
The above assumption yields the commonly used second-order PF model \cite{miehe2010thermodynamically}. On the other hand, the $C^1$ function 
\begin{equation}
\PF(X)=\left(1+\frac{|X|}{\lzr} \right) \, \exp\left(-|X|/\lzr\right) ,
\end{equation}
is obtained by introducing the crack density function 
\begin{equation} \label{eq:CrackDensityFunction}
\gamma_{\lzr}(\PF, \PFgrad, \PFlap)  = \frac{1}{4\lzr} \left[ d^2+2\lzr^2 \vert \PFgrad \vert^2 + \lzr^4 (\PFlap)^2 \right].
\end{equation}
\citet{borden_higher-order_2014} have derived their fourth-order PF theory following this approach. Specialized numerical schemes can be employed with this higher-order model to attain higher spatial convergence rates, since the solution is smoother in comparison (and differentiable at $X=0$ in contrast) to the second-order PF theory. For comparison, the second- and fourth-order model solutions are computed using the same length-scale $\lzr = 0.08$, and plotted in \cref{fig:PhaseFieldSolutionProfiles1D}. It is worth noting that although the higher-order form \eqref{eq:CrackDensityFunction} is adopted throughout the present work, the theoretical development in the following sections is presented generally without specifying the form of $\gamma_{\lzr}$.

%\kgc{Need to explain that the 1D form in (4) is sought to be reproduced by defining the crack density function in \eqref{eq:CrackDensityFunction}?}  

%The 1D forms are sought to be reproduced by defining the crack density functions. 

% for the sake of generality.

\subsubsection{Kinematics}
\label{sssec:DeformationNKinematics}
In finite deformation theory, the motion of the deformable body is mathematically described by a mapping $\motion$ between initial and current particle positions. This function is defined as follows
\begin{equation}
\left. \begin{aligned}
\motion \text{ : } & \RefDom \times \mathcal{T} \rightarrow \CurDom \\
		  & (\vect{X},t) \mapsto \vect{x}=\motion(\vect{X},t)  \\
\end{aligned} \right\},
\end{equation} 
and maps a material point $\vect{X}\in\RefDom$ to its deformed position $\vect{x}$ in the current configuration $\CurDom := \motion(\RefDom,t)$ at time $t \in \mathcal{T}$. In this context, the displacement field $\DispRef$ can be expressed as  
\begin{equation}
\DispRef(\vect{X},t)=\motion(\vect{X},t)-\vect{X}.
\end{equation}

The deformation gradient, defined by $\DeformGrad = \nabla \motion$, can be multiplicatively decomposed into the spherical component $\tens{\Theta}_\circ$ and the deviatoric component $\thickbar{\DeformGrad}$ as follows
\begin{equation}
\DeformGrad = \tens{I}+ \nabla \DispRef =  \tens{\Theta}_\circ \thickbar{\DeformGrad} = J^\frac{1}{3} \, \thickbar{\DeformGrad},
\end{equation}
where $\tens{I}$ is the second-order identity tensor and $J= det[\DeformGrad] >0$ and the components are expressed as
\begin{equation}
\left. \begin{aligned}
\tens{\Theta}_\circ  & = J^\frac{1}{3} \tens{I}  \\
\thickbar{\DeformGrad} & =J^{-\frac{1}{3}} \DeformGrad	   \\
\end{aligned} \right\}.
\end{equation}
The determinant of the deformation gradient, $J$, gives the ratio of the elemental deformed volume $\ud \upsilon$ (spatial configuration) to the corresponding material volume $\ud V$ (material configuration) at a material point, and consequently the constraint $J>0$ always applies. Moreover, the conditions $J>1$ and $J<1$ are respectively associated with the states of volumetric dilatation and contraction.

The right Cauchy-Green tensor $\tens{C}$ is defined by
\begin{equation}
\tens{C} =\DeformGrad\transpose \DeformGrad,
\end{equation}
and the part associated with the deviatoric component is given by
\begin{equation}
\thickbar{\tens{C}} = \thickbar{\DeformGrad}\transpose \thickbar{\DeformGrad}.
\end{equation}
Finally, the Green-Lagrange strain tensor $\tens{E}$ is defined as follows:
\begin{equation}
\tens{E} = \frac{1}{2} \left(\tens{C} - \tens{I} \right).
\end{equation}

%A detailed survey on generalized strain and stress measures is provided in \cite{curnier1991generalized}.

\subsubsection{Boundary conditions}
\label{ssec:NotationForBoundaryConditions}

For each unknown field $Y\in\{ \DispRef, \PF \}$, the exterior boundary can be split into two mutually exclusive parts $\ExtBoun_D^Y$ and $\ExtBoun_N^Y$ to impose essential and natural boundary conditions. This is expressed as follows
\begin{equation} \label{eq:BoudaryConditionDefinitionCondition}
\partial \RefDom = \ExtBoun_D^Y \cup \ExtBoun_N^Y \text{ such that } \ExtBoun_D^Y \cap \ExtBoun_N^Y = \emptyset.
\end{equation}
Depending on the order of a governing equation, multiple essential and natural boundary conditions (and consequently more sets $\ExtBoun_{D_i}^Y$ and $\ExtBoun_{N_i}^Y$, $i \in \mathbb{N}$) may be required to define problems with unique solutions (i.e., well-posed problems). 

Here, the evolution of the displacement field is governed by a second-order hyperbolic equation. Hence, along with initial conditions, boundary conditions on $\ExtBoun_D^\DispRef$ and $\ExtBoun_N^\DispRef$ are sufficient for the problem's well-posedness. These sets are used to prescribe displacements and tractions on the exterior boundary as shown in \cref{fig:DiscreteCrackSmall}. However, the fourth-order PF equation requires two pairs of sets, namely ($\ExtBoun_{D_1}^\PF$, $\ExtBoun_{D_2}^\PF$), and ($\ExtBoun_{N_1}^\PF$, $\ExtBoun_{N_2}^\PF$) as shown in \cref{fig:DiffCrackSmall_BC}. The first tuple corresponds to essential boundary conditions for $\PF$ and $\PFgrad \cdot \vect{N}$, while the second tuple is related to natural boundary conditions which involve second and third spatial derivatives of $\PF$ \cite{stogner2008approximation}. Applying \cref{eq:BoudaryConditionDefinitionCondition} to the aforementioned sets, we arrive at the relations
\begin{equation}
\left. \begin{aligned}
\ExtBoun_{N_1}^\PF  & = \partial \RefDom \setminus  \ExtBoun_{D_1}^\PF  \\
\ExtBoun_{N_2}^\PF  & = \partial \RefDom \setminus  \ExtBoun_{D_2}^\PF  \\
\end{aligned} \right\}.
\end{equation}

The notation introduced above relates to the case of mixed boundary conditions, where Dirichlet or Neumann boundary conditions must hold on disjoint parts of the boundary. Robin boundary conditions, which are linear combinations of Dirichlet and Neumann boundary conditions, can be applied instead on $\ExtBoun_{R_i}^Y = \partial \RefDom$ where $Y\in\{ \DispRef, \PF \}$ and $i \in \mathbb{N}$. In this case, the boundary need not be split into disjoint parts.

\subsection{Governing equations}
\label{ssec:GoverningEquations}

In this section, we first derive the governing equations of dynamic fracture problems in finite elasticity by employing Hamilton's principle. This first-principles derivation is of importance for studying dynamic coupled problems because it is consistent and can be generalized by introducing additional terms in the potential energy functional. Second, we present a novel derivation of the fourth-order phase-field model using a micromechanics-based approach. The latter approach provides a micromechanical interpretation of the higher-order terms in the model. Finally, we provide the hyper-elastic constitutive law which includes a softening mechanism caused by damage.

\subsubsection{Derivation from Hamilton's principle}
\label{ssec:HamiltonDerivationGovEqFourthOrderTheory}

Derived from D'Alembert's principle \cite{wu_chapter_2020}, the extended Hamilton principle for deformable bodies states that
\begin{equation} \label{eq:HamiltonPrinciple}
\int_{t_1}^{t_2} \left( \delta \mathcal{L} + \delta \mathcal{W}_{ext} \right) \ud \tau = 0 \qquad \forall\, t_1, t_2 \in \mathbb{R}^+,
\end{equation}
where the Lagrangian $\mathcal{L}$ is the difference between the kinetic energy functional $\mathcal{K}$ and the potential energy functional $\mathcal{E}$, i.e.\ $\mathcal{L} = \mathcal{K} - \mathcal{E}$, $\delta$ denotes the first variation of a functional (as defined in the calculus of variations), regarded as a virtual field of the system, and $\mathcal{W}_{ext}$ is the external work of the body forces $\BodyForceRef$ and tractions $\SurfTractionRef$ acting on the body and its boundary, respectively, and defined as folows:
\begin{equation} \label{eq:ExternalEnergy}
\mathcal{W}_{ext} = \int_\RefDom \BodyForceRef \cdot \DispRef \ud V + \int_{\ExtBoun_N^\DispRef} \SurfTractionRef \cdot \DispRef \ud A.
\end{equation}

%\cite{Baddour2007Hamilton,wu_chapter_2020}

To guarantee mass conservation, the mass density $\rho_0$ is independent of the phase-field variable, and consequently, the kinetic energy functional is given by
\begin{equation} \label{eq:KineticEnery}
\mathcal{K} = \int_\RefDom\frac{1}{2} \rho_0 \dt{\DispRef} \cdot \dt{\DispRef} \ud V + \int_\RefDom\frac{1}{2} \MicroDens (\dt{\PF})^2 \ud V,
\end{equation}
where the micro-density $\MicroDens$ accounts for the local inertia of material at the crack tip (micro-inertia effects) \cite{stumpf_micromechanical_2003}. The potential energy functional has the following general form,
\begin{equation} \label{eq:PotentialEnergyFunctionalGeneralForm}
\mathcal{E} = \int_\RefDom \TotalEnerDenFun \ud V = \int_\RefDom \rho_0 \hat{\psi} \ud V,
\end{equation}
where $\TotalEnerDenFun$ denotes the total energy density function per unit reference volume. It relates the energy density, stored at a material point due to different physical processes (e.g.\ elastic deformation and crack formation), to the state variables. Moreover, $\hat{\psi}$ is the corresponding specific Helmholtz free energy function, defined by the relation $\TotalEnerDenFun=\rho_0 \hat{\psi}$. In this work, it is assumed that the total energy density function $\TotalEnerDenFun$ can be additively decomposed, as shown below, into the elastic strain energy density function $\DegElastEnerDen$, which is degraded due to damage/fracture, and the surface energy density function $\SurfEnerDenFun $, related to the creation of new crack surfaces:
\begin{equation} \label{eq:TotalEnergyDensityAdditiveDecomposition}
\TotalEnerDenFun (\DeformGrad, \PF, \PFgrad, \PFlap) = \DegElastEnerDen(\DeformGrad,\PF) + \SurfEnerDenFun(\PF, \PFgrad, \PFlap).
\end{equation}
To define $\SurfEnerDenFun$, we assume that the exact surface energy can be approximated by a domain integral over $\RefDom$:
%\begin{equation} \label{eq:WfDefinitionRelations}
 %\int_{\Gamma_c} G_c \ud A = G_c \Gamma_c \approx \int_\RefDom G_c \gamma_{\lzr}(\PF, \PFgrad, \PFlap) \ud V =: \int_\RefDom \SurfEnerDenFun \ud V   \Rightarrow \SurfEnerDenFun = G_c \gamma_{\lzr}(\PF, \PFgrad, \PFlap) \, .
%\end{equation}
\begin{equation} \label{eq:WfDefinitionRelations}
 \int_{\Gamma_c} G_c \ud A = G_c \Gamma_c \approx \int_\RefDom \SurfEnerDenFun \ud V,
\end{equation}
and making use of the approximate relation in \cref{eq:GammaApproximationDomainIntegral}, we obtain
\begin{equation} \label{eq:DefineSurfEnerDenFun}
\SurfEnerDenFun = G_c \gamma_{\lzr}(\PF, \PFgrad, \PFlap).
\end{equation}
% is the summation of the degraded elastic free energy functional (which can be stored in the body) and the surface energy functional as follows
Combining \cref{eq:PotentialEnergyFunctionalGeneralForm,eq:TotalEnergyDensityAdditiveDecomposition,eq:WfDefinitionRelations}, we find the potential energy functional that describes the physics of the fracture process:
\begin{equation} \label{eq:PotentialEnergy}
\mathcal{E} = \underbrace{\int_\RefDom \DegElastEnerDen (\DeformGrad,\PF) \ud V}_\text{Elastic energy functional} + \underbrace{\int_\RefDom G_c \gamma_{\lzr}(\PF, \PFgrad, \PFlap) \ud V}_\text{Surface energy functional}.
\end{equation}
%\begin{equation} 
%\mathcal{E} = \int_{\RefDom \setminus \Gamma_c} \!\!\!\!\!\!\!\! W_e(\DeformGrad) \ud V+\int_{\Gamma_c} G_c \ud A
%\end{equation}
The degraded elastic strain energy density function $\DegElastEnerDen$ depends on the choice of the constitutive law, which is provided in \cref{ssec:ElasticConstitutiveLaw}. Using \cref{eq:KineticEnery,eq:PotentialEnergy,eq:ExternalEnergy}, the integrand in \cref{eq:HamiltonPrinciple} is expressed as
\begin{equation} \label{eq:IntegrandHamiltonPrinciple}
\begin{split}
\delta \mathcal{L} + \delta \mathcal{W}_{ext} & = \delta \mathcal{K} - \delta \mathcal{E} + \delta \mathcal{W}_{ext} = \\
& = \CompInt{\RefDom}{\rho_0\dt{\DispRef} \cdot \delta \dt{\DispRef}}{V} + \CompInt{\RefDom}{\MicroDens  \dt{\PF} \, \delta \dt{\PF}}{V}
- \CompInt{\RefDom}{\partder{\DegElastEnerDen}{\DeformGrad} : \nabla \delta \DispRef}{V} - \CompInt{\RefDom}{\partder{\DegElastEnerDen}{\PF} \delta \PF}{V}   \\
& - \CompInt{\RefDom}{ G_c \partder{\gamma_{\lzr}}{\PF} \delta \PF }{V} - \CompInt{\RefDom}{ G_c \partder{\gamma_{\lzr}}{\PFgrad}  \virtPFgrad }{V} - \CompInt{\RefDom}{ G_c \partder{\gamma_{\lzr}}{\PFlap} \virtPFlap }{V} \\
& + \CompInt{\RefDom}{\BodyForceRef \cdot \delta \DispRef}{V} + \CompInt{\ExtBoun_N^\DispRef}{ \SurfTractionRef \cdot \delta \DispRef }{A}. \\
\end{split}
\end{equation}

After some algebraic manipulations and the substitution of \cref{eq:IntegrandHamiltonPrinciple} into \cref{eq:HamiltonPrinciple} (see \ref{apdx:HamiltonDerivationDetails} for details), we derive the following governing equations (in local form)
\begin{equation} \label{eq:GoverningEqsHamilton}
\left. \begin{aligned}
& \rho_0\ddt{\DispRef} = \Div \left( \partder{\DegElastEnerDen}{\DeformGrad} \right) + \BodyForceRef  \\
&  \MicroDens  \ddt{\PF} = - \partder{\DegElastEnerDen}{\PF} - G_c \partder{\gamma_{\lzr}}{\PF} + G_c \Div \left(\partder{\gamma_{\lzr}}{\PFgrad} \right) - G_c \Lapl \left( \partder{\gamma_{\lzr}}{\PFlap} \right) \\
\end{aligned} \right\}   \text{ in } \RefDom
\end{equation}
in conjunction with the boundary conditions
\begin{equation} \label{eq:BoundaryConditionsHamilton}
\left. \begin{aligned}
& \SurfTractionRef =  \partder{\DegElastEnerDen}{\DeformGrad} \cdot  \NormVecRef \text{ on } \ExtBoun_N^\DispRef   \\
& \left[  \Grad \left( \partder{\gamma_{\lzr}}{\PFlap} \right) -  \partder{\gamma_{\lzr}}{\PFgrad} \right] \cdot \NormVecRef   = 0  \text{ on } \partial \RefDom \\
& \partder{\gamma_{\lzr}}{\PFlap} = 0  \text{ on } \partial \RefDom \\
\end{aligned} \right\}.
\end{equation}
where $\Div$ denotes the divergence operator with respect to the reference configuration, and $\NormVecRef$ is the outward normal vector on the boundary $\partial \RefDom$. It is noteworthy that both governing equations and boundary conditions are derived from Hamilton's principle. This approach leads to the same set of equations as the variational derivation of Ref.~\cite{borden_higher-order_2014}.

\subsubsection{Micromechanics-based derivation}
\label{ssec:MicroMechanicsBasedInterpretationFourthOrderT}

% in a micro-mechanical framework
%A micro-mechanics-based interpretation of the fourth-order theory

By postulating a balance of micro-forces and using it in conjunction with a dissipation inequality (second law of thermodynamics), it is possible to derive \emph{second-order} phase-field model equations, e.g.\ as shown in \cite{mcauliffe_unified_2015,borden_phase-field_2016}. The critical assumption is that cracks propagate under the influence of micro-forces \cite{stumpf_micromechanical_2003}. In this section, we present this alternative approach to the derivation of the governing equations \eqref{eq:GoverningEqsHamilton} and boundary conditions \eqref{eq:BoundaryConditionsHamilton} of the \emph{fourth-order} PF model, with the goal of elucidating the physical meaning of the higher-order terms appearing in these equations.  

%\kgc{Gurtin, Cermelli, Fried and co. may have derived phase field models from microforce balances. It may be that this hasn't been done for phase field fracture, however.}

First, the standard balances of macro-forces and moments are expressed in a Lagrangian description \cite{tadmor_miller_elliott_2011} as follows
\begin{equation} \label{eq:MacrobalanceLocalLagr}
\left. \begin{aligned}
 \rho_0 \ddt{\DispRef} & = \Div \tens{P} + \BodyForceRef   \\
 \tens{S} & = \tens{S}\transpose = \DeformGrad^{-1} \tens{P}
\end{aligned} \right\} \, \text{ in } \RefDom,
\end{equation}
where $\tens{P}$ and $\tens{S}$ are the first and second Piola-Kirchhoff stresses. 

Inspired by the Toupin-Mindlin strain gradient theory \cite{toupin1962elastic,mindlin_micro-structure_1964}, a balance of microforces with higher-order micro-stresses is postulated in the following form \cite{engel_continuousdiscontinuous_2002}:
\begin{equation}  \label{eq:MicrobalanceLocalLagrHighOrd}
\MicroDens  \ddt{\PF} = -\Delta \Phi + \Div \thickbar{\vect{\Sigma}} + M^*,
\end{equation}
where $\MicroDens$ accounts for micro-inertia effects, $ \thickbar{\vect{\Sigma}}$ is the micro-stress vector, $\Phi$ denotes the higher-order micro-stress, and $M^* \in \mathbb{R}$ describes internal scalar-valued micro-forces in material form.

% $\vect{\Sigma}$ denotes the micro-stress vector,

%To the best of our knowledge, this is the first work in which fourth-order phase-field fracture models are derived using a micromechanics-based approach.
% The main deviation from the literature

%Next, the fourth-order PF fracture model equations are derived by introducing a novel micro-stress decomposition, which sheds light on higher-order damage terms. Specifically, the original idea is to decompose the micro-stress vector into a conservative component  $\widehat{\vect{\Sigma}}$, determined by a scalar potential function $\Phi$, and the non-conservative remaining part $ \thickbar{\vect{\Sigma}}$ as follows
%\begin{equation}
%\vect{\Sigma} =  \widehat{\vect{\Sigma}} + \thickbar{\vect{\Sigma}} = \left( - \nabla \Phi + \thickbar{\vect{\Sigma}} \right).
%\end{equation} 
%Substituting this expression into \cref{eq:MicrobalanceLocalLagr}, a novel high-order micro-balance equation is derived 

%+ \int_\ArbDom \Phi \Delta \dt{\PF} \ud V

%The internal energy of the system is denoted by
%\begin{equation}
%\mathcal{I} = \int_\RefDom \rho_0 \hat{e}  \ud V
%\end{equation}
Defining the internal energy per unit mass, $\hat{e}$, in terms of the specific Helmholtz free energy and entropy per unit mass, $\hat{\eta}$:
\begin{equation} \label{eq:LegendreTransformation}
\rho_0 \hat{e} = \rho_0 \hat{\psi} + \rho_0 T \hat{\eta} ,
\end{equation}
where $T$ denotes the absolute temperature, and following the derivation in \ref{apdx:MicromechanicalDerivationDetails}, the local form of the energy balance equation reads
\begin{equation} \label{eq:LocalEnergyBalance}
\rho_0 \dt{\hat{e}} =  \tens{P} : \dt{\DeformGrad} + \thickbar{\vect{\Sigma}} \cdot \nabla \dt{\PF} - M^* \dt{\PF} + \Phi \Delta \dt{\PF}.
\end{equation}

%A mechanical form of the second law of thermodynamics states that the summation of kinetic and potential energies over the reference volume $\RefDom$ increases at a rate not greater than the power of all external forces. This law is expressed via the dissipation inequality 
%\begin{equation} \label{eq:DissipationInequalityIntergralForm}
%\dt{\mathcal{K}} + \dt{\mathcal{E}}  \leq  \int_{\RefDom} \BodyForceRef \cdot \dt{\DispRef} \ud V + \int_{\partial \RefDom}  ( \tens{P} \cdot \NormVecRef ) \cdot  \dt{\DispRef} \ud A  + \int_{\partial \RefDom} \dt{\PF} \,  \thickbar{\tens{\Sigma}} \cdot \NormVecRef   \ud A + \int_{\partial \RefDom} \Phi \,  \nabla \dt{\PF} \cdot  \NormVecRef  \ud A.
%\end{equation}
%It is noteworthy that (i) the first two terms of the power in the right-hand-side of \cref{eq:DissipationInequalityIntergralForm} is equal to $\dt{\mathcal{W}}_{ext}$ (defined in \cref{eq:ExternalEnergy}); (ii) the last two terms are associated with the two components of micro-stress vector; (iii) internal micro-forces do not contribute to the external power.

%Multiplying \cref{eq:MacrobalanceLocalLagr} by $\dt{\DispRef}$, \cref{eq:MicrobalanceLocalLagrHighOrd} by $\dt{\PF}$, and integrating over the domain $\RefDom$, the dissipation inequality (\cref{eq:DissipationInequalityIntergralForm}) can be simplified in its local form

For isothermal processes, the second law of thermodynamics can be simplified as follows \cite{tadmor_miller_elliott_2011}:
\begin{equation} \label{eq:DissipationInequalityIsothermal}
\rho_0 T \dt{\hat{\eta}} \ge 0 \quad  \Rightarrow \quad \dt{\TotalEnerDenFun} = \rho_0 \dt{\hat{\psi}} \leq \rho_0 \dt{\hat{e}}  .
\end{equation} 
By substituting \cref{eq:LocalEnergyBalance} into \eqref{eq:DissipationInequalityIsothermal}, we arrive at the dissipation inequality for our system
\begin{equation} \label{eq:LocalDissipationInequality}
\dt{\TotalEnerDenFun} - \tens{P} : \dt{\DeformGrad} - \thickbar{\vect{\Sigma}} \cdot \nabla \dt{\PF} + M^* \dt{\PF} - \Phi \Delta \dt{\PF} \leq 0.
\end{equation}

We assume that constitutive equations are functions of $(\DeformGrad,\PF, \PFgrad, \PFlap)$ in agreement with the premise in \cref{eq:TotalEnergyDensityAdditiveDecomposition}, that is
\begin{equation} \label{eq:ConstitutiveRelationsDependenceAssumpt}
\left. \begin{aligned}
\TotalEnerDenFun &=\TotalEnerDenFun(\DeformGrad,\PF, \PFgrad, \PFlap) \\
\tens{P} &= \tens{P}(\DeformGrad,\PF, \PFgrad, \PFlap) \\
\thickbar{\vect{\Sigma}} &= \thickbar{\vect{\Sigma}}(\DeformGrad,\PF, \PFgrad, \PFlap) \\
M^* &= M^*(\DeformGrad,\PF, \PFgrad, \PFlap) \\
\Phi &= \Phi(\DeformGrad,\PF, \PFgrad, \PFlap) 
\end{aligned} \right\}.
\end{equation}
Substituting the rate of total energy density function $\TotalEnerDenFun$ which is expressed by
\begin{equation}
\dt{\TotalEnerDenFun} = \partder{\TotalEnerDenFun}{\DeformGrad} : \dt{\DeformGrad} + \partder{\TotalEnerDenFun}{\PF}  \dt{\PF} + \partder{\TotalEnerDenFun}{\PFgrad} \cdot \nabla \dt{\PF} + \partder{\TotalEnerDenFun}{\PFlap}  \Delta\dt{\PF}
\end{equation}
into \cref{eq:LocalDissipationInequality} and collecting the terms into groups, the local dissipation inequality reads
\begin{equation}
\left( \partder{\TotalEnerDenFun}{\DeformGrad} - \tens{P}  \right) : \dt{\DeformGrad} +  \left( \partder{\TotalEnerDenFun}{\PF} + M^* \right)  \dt{\PF} +\left( \partder{\TotalEnerDenFun}{\PFgrad}  - \thickbar{\vect{\Sigma}} \right) \cdot \nabla\dt{\PF} + \left( \partder{\TotalEnerDenFun}{\PFlap} - \Phi \right)  \Delta\dt{\PF}    \leq 0.
\end{equation}
Following the Coleman-Noll procedure \cite{coleman_thermodynamics_1963} and assuming that the only dissipative process is associated with $\dt{\PF}$, we arrive at the relations
\begin{equation} \label{eq:RestrictionsBasedonThermodyna}
\left. \begin{aligned}
\tens{P} & = \partder{\TotalEnerDenFun}{\DeformGrad} =\partder{\DegElastEnerDen}{\DeformGrad} \\
M^* & = - \partder{\TotalEnerDenFun}{\PF} - \MicroDamp  \dt{\PF}   \\
\thickbar{\vect{\Sigma}} &= \partder{\TotalEnerDenFun}{\PFgrad} = \partder{\SurfEnerDenFun}{\PFgrad}  \\
\Phi &=  \partder{\TotalEnerDenFun}{\PFlap} =   \partder{\SurfEnerDenFun}{\PFlap}
\end{aligned} \right\} ,
\end{equation}
where the micro-viscosity parameter $\MicroDamp \ge 0$ accounts for the local damping of material at the crack tip (micro-damping effects). These effects produce entropy in the system through the dissipative process which is associated with $\dt{\PF}$.

% = - \partder{\DegElastEnerDen}{\PF} - \partder{\SurfEnerDenFun}{\PF} - \MicroDamp  \dt{\PF}

The high-order micro-force balance equation (\cref{eq:MicrobalanceLocalLagrHighOrd}) is reformulated using \cref{eq:RestrictionsBasedonThermodyna} and the additive decomposition of $\TotalEnerDenFun$ (\cref{eq:TotalEnergyDensityAdditiveDecomposition}) as
\begin{equation} \label{eq:HigherOrderEquationGeneral}
\MicroDens  \ddt{\PF} + \MicroDamp  \dt{\PF}  = - \Delta \left(  \partder{\SurfEnerDenFun}{\PFlap} \right)  + \Div \left( \partder{\SurfEnerDenFun}{\PFgrad} \right)   - \partder{\DegElastEnerDen}{\PF} - \partder{\SurfEnerDenFun}{\PF},
\end{equation}
which is identical with \cref{eq:GoverningEqsHamilton} after the substitution of \cref{eq:DefineSurfEnerDenFun} and $\MicroDamp = 0$. For comparison purposes, micro-damping effects are ignored in the latter equation since Hamilton's principle is only applied to conservative systems (i.e. without dissipative mechanisms).

The boundary conditions of phase-field equation (given in \cref{eq:BoundaryConditionsHamilton}) can also be re-derived by imposing the following conditions (the reader is referred to \ref{apdx:BoundaryconditionsPFequation} for the derivation)
\begin{equation} \label{eq:BoundaryConditionsThermodyDerivCond}
\left. \begin{aligned}
\left( - \nabla \Phi + \thickbar{\vect{\Sigma}} \right) \cdot \NormVecRef  =0 &     \\
\Phi = 0   &   
\end{aligned} \right\}  \text{ on } \partial \RefDom.
\end{equation}
The first one expresses a microtraction-free condition whereas the second zeroes out the higher-order micro-stress $\Phi$ on the boundary. This observation provides a micromechanics-based interpretation of fourth-order phase-field models.

%This observation provides a physical meaning to the current derivation based on micro-mechanics. 
%This interpretation allows us to extend the model and 

This interpretation delineates the physical meaning of the boundary conditions in higher-order damage theories and allows us to impose additional physics-based boundary conditions. For instance, the micro-stress vector $\thickbar{\vect{\Sigma}}$ can satisfy a traction-free condition as follows
\begin{equation} \label{eq:AdditionalBCHOPFThermoDynamics}
\thickbar{\vect{\Sigma}} \cdot \NormVecRef = 0 \text{ on } \partial \RefDom.
\end{equation}
This equation can simplify considerably the boundary conditions of our system and imposes the condition that the damage does not diffuse across external boundaries.

%Finally, we can show that the work done by the conservative micro-stress vector is zero on each curve in the reference domain which connects two points of boundary $\partial \RefDom$. Let $C_\gamma$ be any curve in $\RefDom$ which starts at $\bar{\vect{p}} \in \partial \RefDom$ and ends at $\bar{\vect{q}} \in \partial \RefDom$, the line integral
%\begin{equation}
%\int_{C_\gamma} \widehat{\vect{\Sigma}} \cdot \ud \vect{r} = - \int_{C_\gamma} \nabla \Phi \cdot \ud \vect{r} = \Phi(\bar{\vect{p}}) - \Phi(\bar{\vect{q}}) = 0 ,
%\end{equation}
%where $\vect{r}$ is a bijective parametrization of the curve $C_\gamma$, and $\Phi(\bar{\vect{p}}) = \Phi(\bar{\vect{q}}) = 0$, since $\bar{\vect{p}},\bar{\vect{q}} \in \partial \RefDom$ (\cref{eq:BoundaryConditionsThermodyDerivCond}). A special case of the aforementioned result holds for $C_\gamma = \partial \RefDom$.

\begin{remark} \label{rmk:GeneralFormulationComm}
The derivation of high-order phase-field models using the micromechanics-based approach comprises a powerful alternative to the derivation from Hamilton's principle since (i) dissipative terms can be naturally included (e.g., micro-damping effects associated with  $\dt{\PF}$) and (ii) additional physics-based boundary conditions can be imposed.
\end{remark}

\subsection{Hyper-elastic constitutive law}
\label{ssec:ElasticConstitutiveLaw}

To introduce the elastic constitutive law in this model, we specify the elastic strain energy density function $\DegElastEnerDen$. Specifically, the latter function determines the stress-strain relation through \cref{eq:RestrictionsBasedonThermodyna}
\begin{equation} \label{eq:FirstPKstressDefGradientRelationF}
\tens{P} =\partder{\DegElastEnerDen(\DeformGrad,\PF)}{\DeformGrad}.
\end{equation}

% the response should be invariant under changes in observer

Although constitutive equations were assumed functions of $(\DeformGrad,\PF, \PFgrad, \PFlap)$ in \cref{eq:ConstitutiveRelationsDependenceAssumpt}, the response should be independent of the frame of reference (i.e., constitutive model should be frame-indifferent/objective). To this end, the dependence of constitutive equations on state variables is modified as follows
\begin{equation}
\left. \begin{aligned}
\TotalEnerDenFun &=\TotalEnerDenFun(\tens{C},\PF, \PFgrad, \PFlap) \\
\tens{S} &= \tens{S}(\tens{C},\PF, \PFgrad, \PFlap) \\
\thickbar{\vect{\Sigma}} &= \thickbar{\vect{\Sigma}}(\tens{C},\PF, \PFgrad, \PFlap) \\
M^* &= M^*(\tens{C},\PF, \PFgrad, \PFlap) \\
\Phi &= \Phi(\tens{C},\PF, \PFgrad, \PFlap) 
\end{aligned} \right\},
\end{equation}
and taking into account the relation $\tens{P} : \dt{\DeformGrad} = \frac{1}{2} \tens{S} : \dt{\tens{C}}$, \cref{eq:FirstPKstressDefGradientRelationF} is reformulated into
\begin{equation} \label{eq:SecondPKStressRelationC}
\tens{S} = 2 \partder{\DegElastEnerDen (\tens{C},\PF) }{\tens{C}}.
\end{equation}
%where  and $\tens{C}$ denotes the right Cauchy-Green tensor.

Adopting a finite strain elasticity model with uncoupled, volumetric/deviatoric response as presented in \cite{simo2006computational}, the \emph{undamaged} elastic strain energy density function is decomposed as follows
\begin{equation} \label{eq:GeneralizedVolDevContributions}
W_e(\tens{C}) = U(J) + \thickbar{W}(\thickbar{\tens{C}}),
\end{equation}
where $U$ and $\thickbar{W}$ denote the volumetric and deviatoric contributions to the energy, and their arguments were defined in \cref{sssec:DeformationNKinematics}. In this work, we select their forms as follows
\begin{equation} \label{eq:VolDevContributionsSpec}
\left. \begin{aligned}
U(J) & = \frac{\kappa}{2}  \left[ \frac{1}{2} \left( J^2 -1 \right) - \ln{J} \right] \\
\thickbar{W}(\thickbar{\tens{C}}) & = \frac{\mu}{2} \left[ tr(\thickbar{\tens{C}}) - 3 \right] \\
\end{aligned} \right\},
\end{equation}
where $\kappa=\frac{E}{3(1-2\nu)}$ is the bulk modulus and $\mu=\frac{E}{2(1+\nu)}$ is the shear modulus, which are given as a function of Young's modulus $E$ and Poisson's ratio $\nu$. This type of energy split in isochoric (volume preserving) and volumetric (shape preserving) parts is used to model elastic materials such as foam or rubber but it is also appropriate choice to couple it with plasticity models \cite{ambati2016phase}.

%Alternative decomposition for incompressible rubberlike solids are introduced in \cite{ogden1972large}.

\begin{figure}[H] 
\centering
\includegraphics[width=0.45\textwidth]{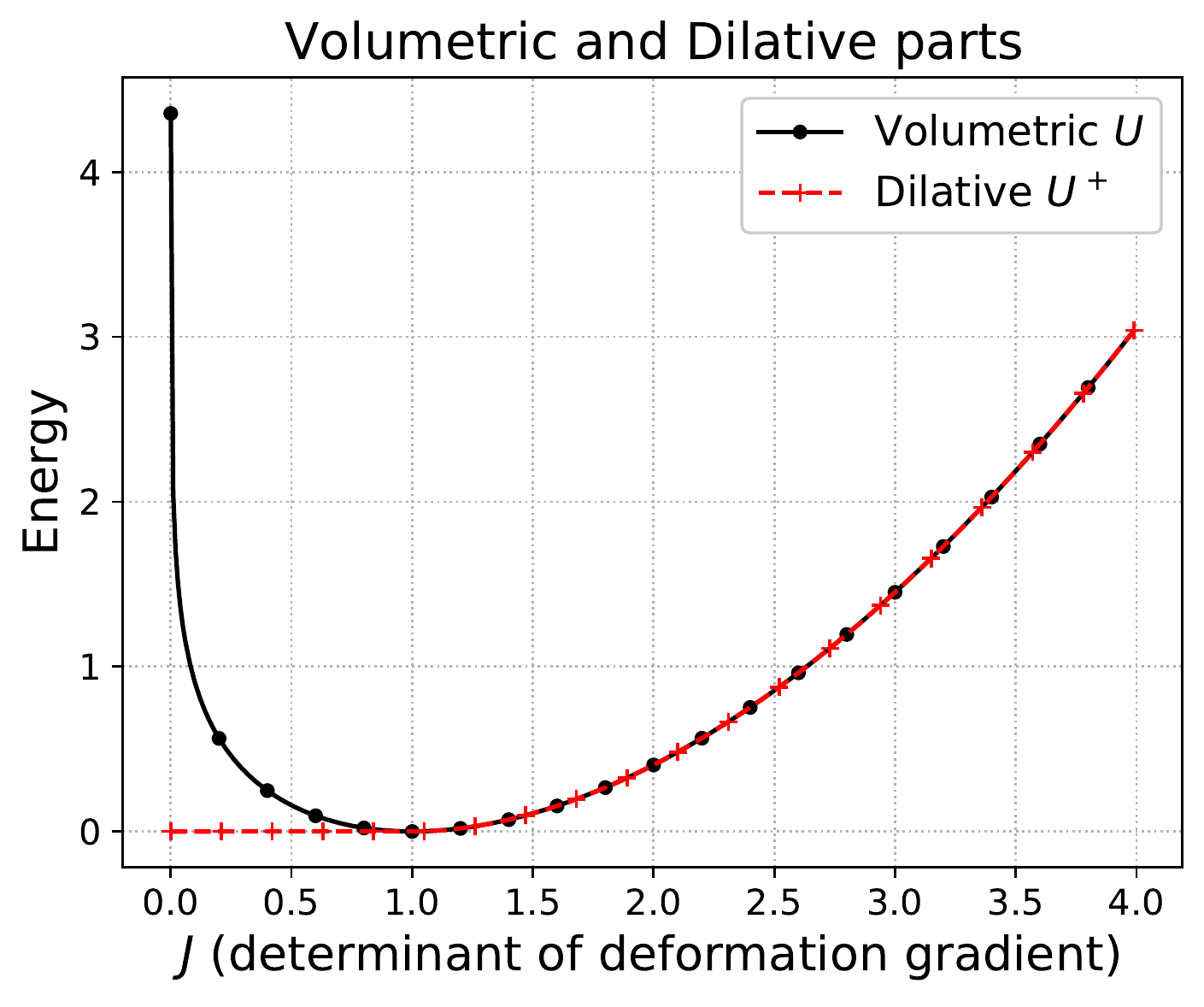}
\captionsetup{width=.90\textwidth}
\caption{Volumetric and dilative parts of the elastic strain energy function. Notice that the dilative term $U^+(J)$ is continuously differentiable function. }
\label{fig:VolDilativePartsPlot1DasFunctionofDetF}
\end{figure}

Assuming that the volumetric contribution $U$ at a material point is not released under volume contraction (i.e., $J < 1$), we define the dilative term $U^+(J)$  using the Heaviside function $H$ as
%we split $U$ into dilative and contractive components ($U^+$ and $U^-$ respectively)
%\begin{equation}
%U(J)  = U^+(J)  + U^-(J) 
%\end{equation}
%The dilative component can be expressed using the Heaviside function $H$ as
\begin{equation} \label{eq:DilativeTermDef0}
U^+(J) = H(J-1) \, U(J) = 
\left\{
	\begin{array}{ll}
		0  & \mbox{if }  0< J < 1 \\
		U(J) & \mbox{if } J  \geq 1
	\end{array}
\right.
\end{equation}
%and in more compact form using Macaulay brackets $\langle \cdot \rangle$
%\begin{equation}
%U^+(J) = \int_{1}^{J} \left\langle U^\prime(\Theta) \right\rangle  \ud \Theta
%\end{equation}
%where the operator is defined by 
%\begin{equation}
%\langle x \rangle =
%\left\{
%	\begin{array}{ll}
%		0  & \mbox{if } x < 0 \\
%		x & \mbox{if } x \geq 0
%	\end{array}
%\right\} \, .
%\end{equation}
The volumetric and dilative parts of the elastic strain energy function are depicted in \cref{fig:VolDilativePartsPlot1DasFunctionofDetF}. It is noteworthy that the dilative term $U^+(J)$ is continuously differentiable function although the Heaviside function is used in the definition. With this regularity stable numerical results are expected in the transition from dilation to expansion.

With the above volumetric decomposition, we can define the \emph{degraded} elastic strain energy density function as follows
\begin{equation}
\DegElastEnerDen(\tens{C},\PF) =
\left\{
	\begin{array}{ll}
		U(J) +  g(\PF)  \thickbar{W}(\thickbar{\tens{C}})   & \mbox{if }  0< J < 1 \\
		g(\PF) \left[ U(J)  + \thickbar{W}(\thickbar{\tens{C}}) \right]  & \mbox{if } J  \geq 1
	\end{array}
\right.
\end{equation}
and expressed in more compact form after some algebraic manipulations as
%\begin{equation}
%\DegElastEnerDen(\tens{C},\PF) = U^-(J) + g(\PF) \left[ U^+(J) + \thickbar{W}(\thickbar{\tens{C}}) \right] = W_e(\tens{C}) + \left[ g(\PF) - 1 \right] \left[ U^+(J) + \thickbar{W}(\thickbar{\tens{C}}) \right]
%\end{equation}
\begin{equation} \label{eq:CompactFormOfdegradedElasticStrainEnergy}
\DegElastEnerDen(\tens{C},\PF) = W_e(\tens{C}) + \left[ g(\PF) - 1 \right] \left[ U^+(J) + \thickbar{W}(\thickbar{\tens{C}}) \right],
\end{equation}
where $g(\PF)$ is the so-called degradation function which introduces the transition rule from unbroken to fully broken response. It is noteworthy that the volumetric component of elastic energy is not degraded by the phase field in case of volume contraction, and consequently, the interpenetration of crack surfaces is prevented. In this work, we use the standard quadratic degradation function
\begin{equation} \label{eq:QuadraticDegrFunction}
g(\PF) = (1-\eta_0) (1 - \PF)^2 + \eta_0,
\end{equation}
where $\eta_0$ provides a residual stiffness and prevents ill-posedness of the problem when the crack is fully developed $\PF=1$. Additional options of degradation functions can be found for brittle fracture in \cite{kuhn_degradation_2015}.
%and for cohesive fracture in \cite{lorentz2011convergence,lorentz2017nonlocal}.

To find the material response, the second Piola-Kirchhoff in \cref{eq:SecondPKStressRelationC} is calculated
%\begin{equation} \label{eq:SecondPKStressFunctionofCdUplus}
%\tens{S} = \tens{S}_e(\tens{C}) + \left[ g(\PF) - 1 \right] \left[ \tens{S}^+(\tens{C}) + \thickbar{\tens{S}}(\tens{C}) \right] = 2 \partder{W_e(\tens{C})}{\tens{C}}  + 2 \left[ g(\PF) - 1 \right] \left[ \partder{U^+(J)}{\tens{C}} + \partder{\thickbar{W}(\thickbar{\tens{C}})}{\tens{C}}  \right]
%\end{equation}
\begin{equation} \label{eq:SecondPKStressFunctionofCdUplus}
\tens{S} = \tens{S}_e(\tens{C}) + \left[ g(\PF) - 1 \right] \left[ \tens{S}^+(\tens{C}) + \thickbar{\tens{S}}(\tens{C}) \right],
\end{equation}
where the stresses $\tens{S}_e(\tens{C})$, $\tens{S}^+(\tens{C})$, and $\thickbar{\tens{S}}(\tens{C})$ are defined by the derivatives as follows
%\begin{equation}
%\tens{S}_e(\tens{C}) = 2 \partder{W_e(\tens{C})}{\tens{C}} = J U^\prime(J) \tens{C}^{-1} + \thickbar{\tens{S}}(\tens{C}) = \frac{\kappa}{2} \left( J^2 -1 \right) \, \tens{C}^{-1} + \thickbar{\tens{S}}(\tens{C})
%\end{equation}
%with
%\begin{equation}
%\thickbar{\tens{S}}(\tens{C}) =  2\partder{\thickbar{W}(\thickbar{\tens{C}})}{\tens{C}} = J^{-\frac{2}{3}} \left[ \mu \tens{I} - \frac{\mu}{3} tr(\tens{C}) \, \tens{C}^{-1} \right]
%\end{equation}
%and
%\begin{equation}
% \tens{S}^+(\tens{C}) = 2 \partder{U^+(J)}{\tens{C}}  = J \left\langle U^\prime(J) \right\rangle \tens{C}^{-1} = \frac{\kappa}{2} H(J-1) \left( J^2 -1 \right) \tens{C}^{-1}
%\end{equation}

%(undamaged elastic stress)
%\begin{equation}
%\tens{S}_e(\tens{C}) = 2 \partder{W_e(\tens{C})}{\tens{C}}  = \tens{S}_\circ(\tens{C})   + \thickbar{\tens{S}}(\tens{C})
%\end{equation}
%with
%(volumetric stress)
%\begin{equation}
% \tens{S}_\circ(\tens{C})  =  2\partder{U(J)}{\tens{C}} = J U^\prime(J) \tens{C}^{-1} = \frac{\kappa}{2} \left( J^2 -1 \right) \, \tens{C}^{-1}
%\end{equation}
%(isochoric stress)
%\begin{equation}
%\thickbar{\tens{S}}(\tens{C}) =  2\partder{\thickbar{W}(\thickbar{\tens{C}})}{\tens{C}} = \mu J^{-\frac{2}{3}} \left[ \tens{I} - \frac{1}{3} tr(\tens{C}) \, \tens{C}^{-1} \right]
%\end{equation}
%and
%(dilative stress)
%\begin{equation}
% \tens{S}^+(\tens{C}) = 2 \partder{U^+(J)}{\tens{C}}  = H(J-1) \tens{S}_\circ(\tens{C})
%\end{equation}

\begin{equation} \label{eq:StressesInCompactFormandExplan}
\left.
	\begin{array}{ll}
		\mbox{(undamaged elastic stress)} \quad   & \tens{S}_e(\tens{C}) = 2 \partder{W_e(\tens{C})}{\tens{C}}  = \tens{S}_\circ(\tens{C})   + \thickbar{\tens{S}}(\tens{C}) \\
		\mbox{(volumetric stress)}  &  \tens{S}_\circ(\tens{C})  =  2\partder{U(J)}{\tens{C}} = J U^\prime(J) \tens{C}^{-1} = \frac{\kappa}{2} \left( J^2 -1 \right) \, \tens{C}^{-1} \\
		\mbox{(isochoric stress)}  &  \thickbar{\tens{S}}(\tens{C}) =  2\partder{\thickbar{W}(\thickbar{\tens{C}})}{\tens{C}} = \mu J^{-\frac{2}{3}} \left[ \tens{I} - \frac{1}{3} tr(\tens{C}) \, \tens{C}^{-1} \right] \\
		\mbox{(dilative stress)}  &   \tens{S}^+(\tens{C}) = 2 \partder{U^+(J)}{\tens{C}}  = H(J-1) \tens{S}_\circ(\tens{C})
	\end{array} 
\right\}.
\end{equation}

The reader is referred to \cite{wriggers2008nonlinear} for more details on calculations of the above derivatives. Moreover, other volumetric/deviatoric energy splits are presented in \cite{amor2009regularized,ulloa2019phase}.

\section{Numerical methods: finite element approximations}
\label{sec:Computational}

In this section, we develop numerical methods for the solution of hyper-elastic fracture problems modeled by the fourth-order phase-field equation. In \cref{ss:SummaryStrongForms}, we present the final form of the governing equations and constitutive laws. Specifically, (i)~the momentum balance equation and constitutive law are summarized, (ii)~the phase-field equation is reformulated, making use of a history variable to enforce irreversibility (i.e.\ to prevent crack healing), and (iii)~the phase-field equation is additionally expressed as a coupled system of two second-order partial differential equations (with a splitting approach). In \cref{ss:DomainDiscretizationNinterelementBound}, the domain is discretized using finite elements and notation is introduced to describe discontinuities across inter-element boundaries. In \cref{ss:WeakForms}, semidiscrete weak forms of the governing equations are derived employing (i)~a standard Galerkin formulation for the momentum balance equation, (ii)~a continuous/discontinuous Galerkin method for the fourth-order phase-field equation, and (iii)~a mixed finite element formulation for the phase-field equation based on the splitting approach. In \cref{ss:NumericalImplementStagTimeInt}, we provide details on the numerical integration of the governing equations using a staggered scheme (i.e., a weak coupling between the unknown fields).

\subsection{Summary: strong forms}
\label{ss:SummaryStrongForms}

% In this section, we summarize the governing equations and constitutive laws of the problem. First, we present the momentum balance equation and the dependence of the second Piola–Kirchhoff stress on the phase field. Second, the fourth-order phase-field equation is recapitulated, defining a history variable to enforce crack irreversibility. Finally, the phase-field equation is rewritten as a coupled system of two second-order partial differential equations, based on a standard splitting approach. 

%$\ExtBoun_D^\DispRef$ and $\ExtBoun_N^\DispRef$

%and the dependence of the second Piola–Kirchhoff stress on the phase field. 

\subsubsection{Momentum balance equation \& constitutive law}

Combining \cref{eq:MacrobalanceLocalLagr} with \cref{eq:SecondPKStressFunctionofCdUplus}, the momentum equation reads
%\begin{equation} \label{eq:MomBalEqStrongForm}
%\left. \begin{aligned}
%\rho_0 \ddt{\DispRef} & = \Div \tens{P} + \BodyForceRef   \\
%\tens{S} & = \tens{S}\transpose = \DeformGrad^{-1} \tens{P} 
%\end{aligned} \right\} \, \text{ in } \RefDom
%\end{equation}
\begin{equation} \label{eq:MomBalEqStrongForm}
\rho_0 \ddt{\DispRef}  = \Div \tens{P} + \BodyForceRef  \, \text{ in } \RefDom
\end{equation}
with the constitutive law expressed based on the symmetric second Piola–Kirchhoff stress as follows
\begin{equation}
\tens{S} = \DeformGrad^{-1} \tens{P}  =  \tens{S}_e(\tens{C}) + \left[ g(\PF) - 1 \right] \left[ \tens{S}^+(\tens{C}) + \thickbar{\tens{S}}(\tens{C}) \right],
\end{equation}
where its components were given in \cref{eq:StressesInCompactFormandExplan}.
%\begin{equation}  \label{eq:ConstLawStrongForm}
%\left. \begin{aligned}  
%\tens{S}(\tens{C},\PF) & = \tens{S}_e(\tens{C}) + \left[ g(\PF) - 1 \right] \left[ \tens{S}^+(\tens{C}) + \thickbar{\tens{S}}(\tens{C}) \right] \\
%\tens{S}_e(\tens{C}) &=  \frac{\kappa}{2} \left( J^2 -1 \right) \, \tens{C}^{-1} + \thickbar{\tens{S}}(\tens{C}) \\
%\thickbar{\tens{S}}(\tens{C}) & = J^{-\frac{2}{3}} \left[ \mu \tens{I} - \frac{\mu}{3} tr(\tens{C}) \, \tens{C}^{-1} \right] \\
%\tens{S}^+(\tens{C}) &= \frac{\kappa}{2} H(J-1) \left( J^2 -1 \right) \tens{C}^{-1}  \\ 
%\end{aligned} \right\} \, \text{ in } \RefDom
%\end{equation}
Finally, the boundary conditions are written 
\begin{equation}
\left. \begin{aligned}
\SurfTractionRef = \tens{P} \cdot \NormVecRef =  \tens{F} \tens{S} \cdot  \NormVecRef =  \tens{F} \cdot \vect{s}^* \text{ on } \ExtBoun_N^\DispRef \\
\DispRef = \thickbar{\DispRef}  \text{ on }  \ExtBoun_D^\DispRef
\end{aligned} \right\},
\end{equation}
where $\vect{s}^* $ is the prescribed traction associated with $\tens{S}$, and  $ \thickbar{\DispRef}$ is the prescribed displacement on $\ExtBoun_D^\DispRef$.

\subsubsection{Phase-field equation}

% \left( \PFlap \right)

% $\DegElastEnerDen$ (which is degraded due to fracture) and the surface energy density function $\SurfEnerDenFun $

% the equation derived in \cref{eq:GoverningEqsHamilton}
The fourth-order phase-field equation is derived using the novel approach, presented in \cref{ssec:MicroMechanicsBasedInterpretationFourthOrderT}. First, we find the surface energy density function $\SurfEnerDenFun$ from \cref{eq:DefineSurfEnerDenFun} by selecting the fourth-order crack density function (\cref{eq:CrackDensityFunction}). Second, the \emph{degraded} elastic strain density function $\DegElastEnerDen$ is given in \cref{eq:CompactFormOfdegradedElasticStrainEnergy}. Substituting the aforementioned density functions into \cref{eq:HigherOrderEquationGeneral}, we arrive at the strong form of the phase-field equation:
\begin{equation} \label{eq:PFEquationQuadraticDegrFuncGeneral}
\MicroDens  \ddt{\PF} + \MicroDamp  \dt{\PF} + \alpha_2  \Delta^2 \PF  +  \alpha_1  \PFlap  +  \alpha_0  \PF  + g^\prime(\PF) \left[ U^+(J) + \thickbar{W}(\thickbar{\tens{C}}) \right]  = 0 ,
\end{equation}
where the operator $\Delta^2$ is defined as $\Delta^2 \PF = \Delta  \left( \PFlap \right)$ and the coefficients are given as follows $\alpha_2 = \frac{G_c}{2} \lzr^3$, $\alpha_1 = -  G_c \lzr $, and $\alpha_0 =  \frac{G_c}{2 \lzr}$.
%as follows
%\begin{equation}
%\left. \begin{aligned}
%\alpha_2 &= \frac{G_c}{2} \lzr^3 \\
%\alpha_1 &= -  G_c \lzr  \\
%\alpha_0 &=  \frac{G_c}{2 \lzr} 
%\end{aligned} \right\},
%\end{equation}

The boundary conditions of the phase-field equation can be derived from \cref{eq:BoundaryConditionsThermodyDerivCond} as shown in \ref{apdx:BoundaryconditionsPFequation} by exploiting \cref{eq:CrackDensityFunction} and \cref{eq:RestrictionsBasedonThermodyna}. That is
\begin{equation} \label{eq:BoundaryCondPhaseFieldHighOrderThermoA}
\left. \begin{aligned}
\nabla \left(  \alpha_2 \PFlap + \alpha_1  \PF  \right) \cdot \vect{N}  =0  \text{ on } \partial \RefDom = \ExtBoun_{R_2}^\PF &  \\
  \PFlap  = 0  \text{ on } \partial \RefDom = \ExtBoun_{N_1}^\PF &
\end{aligned} \right\}   .
\end{equation}
% \cref{eq:BoundaryConditionsHamilton} or 
%\lzr^2  \nabla  \PFlap \cdot  \vect{N}   = 2  \nabla \PF   \cdot \vect{N} 

\begin{remark}
Although the proposed formulation is thermodynamically consistent, unloading can result in crack healing (i.e., a decrease of the phase-field values). To enforce crack irreversibility, we follow a similar approach to \cite{miehe2010phase} and replace the summation of strain energy density functions with a history variable $\mathcal{H}$, defined as
\begin{equation}
\mathcal{H} = \max_t \left[ U^+(J) + \thickbar{W}(\thickbar{\tens{C}}) \right].
\end{equation}
\end{remark}

%\kgc{Does the max arise consistently?}.

It is noteworthy that the strong form is the same as the fourth-order phase-field equation presented in \cite{borden_higher-order_2014}. In our work, we derive the strong form based on two different approaches (Hamilton's principle and the micromechanics-based approach) and we show that these two approaches are equivalent. 

%Moreover, the unprecedented derivations, given in \cref{ssec:MicroMechanicsBasedInterpretationFourthOrderT}, delineates the physical interpretation of fourth-order damage theories, based on micro-mechanics principles.

The boundary conditions in \cref{eq:BoundaryCondPhaseFieldHighOrderThermoA} consist of one Robin boundary condition and one Neumann boundary condition. However, they do not explicitly prevent `crack diffusion' across the external boundaries. Given that the micromechanics-based approach is more general (see \cref{rmk:GeneralFormulationComm}), we can impose this physical constraint using \cref{eq:AdditionalBCHOPFThermoDynamics}. After some algebraic manipulations, the aforementioned condition is expressed in the following form 
\begin{equation}
\Grad d \cdot \vect{N} = 0  \text{ on }  \ExtBoun_{D_2}^\PF,
\end{equation}
and the boundary $\ExtBoun_{D_2}^\PF$ can be treated as `crack insulated.' In this way, the new boundary conditions read
\begin{equation} \label{eq:NewBoundaryConditionPhysicsIntr}
\left. \begin{aligned}
 \Grad d \cdot \vect{N} = 0  \text{ on }  \ExtBoun_{D_2}^\PF &  \\
 \nabla \left(  \alpha_2 \PFlap + \alpha_1  \PF  \right) \cdot \vect{N}  =0  \text{ on }  \ExtBoun_{R_2}^\PF &  \\
  \PFlap  = 0  \text{ on } \partial \RefDom  &
\end{aligned} \right\}   .
\end{equation}
%\partial \RefDom =
% the Robin boundary condition can be replaced and

%\begin{equation} 
% \MicroDens  \ddt{\PF} + \alpha_3  \Delta \left( \PFlap \right) +  \alpha_2  \PFlap  +  \alpha_1  \PF   = 2 \mathcal{H} \text{ in } \RefDom
%\end{equation}

%= \frac{G_c}{4\lzr} \left[- 2 \lzr^4 \Delta \left( \PFlap \right)  + 4 \lzr^2 \PFlap  -2 \PF \right]

%\begin{equation} 
%\left. \begin{aligned}
%&  \left( 4 \lzr^2 \PFgrad- 2 \lzr^4 \nabla \PFlap \right) \cdot \vect{N} = - \nabla \left( 2 \lzr^4 \PFlap- 4 \lzr^2 \PF  \right) \cdot \vect{N} =0   \\
% &  - 2 \lzr^4 \PFlap= 0 \Rightarrow \PFlap = 0  
%\end{aligned} \right\}  \text{ on } \partial \RefDom
%\end{equation}

%\subsection{Weak forms}

\subsubsection{Recasting the phase-field equation based on a splitting approach} \label{sss:SplittingApproach}

%We express the phase-field equation in two different ways, because distinct finite element approximations are applied to each form in \cref{ss:WeakForms}.

The solution of the fourth-order phase-field equation, given in \cref{eq:PFEquationQuadraticDegrFuncGeneral}, requires high-order approximation spaces when using conforming finite element formulations. However, a standard alternative technique is to introduce a split form of the phase field equation in which two second-order equations replace the fourth-order equation. Specifically, \cref{eq:PFEquationQuadraticDegrFuncGeneral} can be rewritten as a coupled system of two partial differential equations as follows 
\begin{equation} \label{eq:SecondOrderSplittingSystems}
\left. \begin{aligned}
\PFlap  - \lambda_0 \psi  & = 0\\
\MicroDens  \ddt{\PF}  + \MicroDamp  \dt{\PF} + \lambda_0 \alpha_2 \Delta \psi + \lambda_0 \alpha_1 \psi  +   \alpha_0  \PF  + g^\prime(\PF) \mathcal{H} & = 0 
\end{aligned} \right\},
\end{equation}
and the corresponding boundary conditions as
\begin{equation} \label{eq:BoundaryCondPhaseFieldHighOrderMixed}
\left. \begin{aligned}
\Grad d \cdot \vect{N}  = 0 \text{ on } \ExtBoun_{D_2}^\PF & \\
 \nabla \left(  \alpha_2 \lambda_0  \psi + \alpha_1  \PF  \right) \cdot \vect{N}  =0  \text{ on }  \ExtBoun_{R_2}^\PF & \\
 \psi  = 0  \text{ on } \partial \RefDom  &
\end{aligned} \right\}  .
\end{equation}
% \nabla \psi  \cdot \vect{N} & =0   \\
The parameter $\lambda_0$ is introduced in \cref{eq:SecondOrderSplittingSystems} to control the coefficients of the second-order terms. Depending on their order of magnitude, the numerical solution can be made more efficient by adjusting the value of $\lambda_0$. In the present work, we choose $\lambda_0 = 1$ because we do not observe any convergence difficulties associated with the order of magnitude of these terms. In any case, the phase-field solution $\PF$ is independent of $\lambda_0$. 
% (\cref{eq:NewBoundaryConditionPhysicsIntr}) 

\subsection{Domain discretization and notation of inter-element boundaries}
\label{ss:DomainDiscretizationNinterelementBound}

In this section, we discretize the reference domain into finite elements, and we provide definitions and define elements of notation that are essential to the development of weak forms in the next section. Specifically, the notation of inter-element boundaries is introduced and operators on these boundaries are defined.

The geometry of the deformable body $\RefDom$ in the initial configuration is partitioned as follows
\begin{equation}
\RefDom =  \overline{\bigcup_{e=1}^{n_{el}} \ElRefDom},
\end{equation} 
where its domain is subdivided into $n_{el}$ mutually exclusive \emph{open} sets $\ElRefDom$ (for each finite elements $e = 1, \dots, n_{el}$  with boundary $\ExtBoun_e = \partial \ElRefDom$). The overline denotes the closure of this set in order to include the element boundaries in the domain. Excluding the element boundaries, the union of the element interiors is denoted as
\begin{equation}
\TildeRefDom = \bigcup_{e=1}^{n_{el}} \ElRefDom .
\end{equation} 

%The aforementioned boundary is split into two parts ($\ExtBoun_D$ and $\ExtBoun_N$) depending on the imposed type of boundary conditions 
%\begin{equation}
%\ExtBoun = \ExtBoun_D \cup \ExtBoun_N \text{ such that } \ExtBoun_D \cap \ExtBoun_N = \emptyset,
%\end{equation}
%where Dirichlet conditions are imposed on $\ExtBoun_D$ and Neumann conditions on $\ExtBoun_N$ . 

%In addition, we denote the boundary of each element by . 

The inter-element boundaries play an important role in discontinuous Galerkin methods. Hence, the union of inter-element boundaries is expressed in terms of the intersection of the boundaries $\ExtBoun_e$ of individual elements:
\begin{equation}
\IntOnlyBoun = \bigcup_{\substack{e_1,e_2=1 \\ e_2>e_1}}^{n_{el}} \left( \ExtBoun_{e_1} \cap \ExtBoun_{e_2} \right).
\end{equation}

The exterior boundary of the whole domain is denoted by $\ExtBoun=\partial \AppRefDom$ and is split into parts $\ExtBoun_{D_i}^Y$ and $\ExtBoun_{N_i}^Y$ with $i \in \mathbb{N}$ and $Y\in\{ \DispRef, \PF \}$ as described in \cref{ssec:NotationForBoundaryConditions}. The union of inter-element boundaries and the boundary $\ExtBoun_{D_2}^\PF$ is commonly used in the development that follows, and therefore, we introduce the notation
\begin{equation}
\IntBoun = \IntOnlyBoun \cup  \ExtBoun_{D_2}^\PF.
\end{equation}

% and the outward unit normal vector by $\vect{n}_e^+$ whereas the inward unit normal vector by  $\vect{n}_e^- (= - \vect{n}_e^+)$

Two other constructs that appear frequently in the weak form of the problem are the average of a given quantity evaluated at an inter-element boundary, and the jump in the value of such a quantity across an inter-element boundary. For this reason, we define the average and jump operators, denoted by $\avrg{\bullet}$ and $\jump{\bullet}$, respectively, as follows: 
\begin{equation} \label{eq:AvrgJumpOperatorsDefinitions}
\left. \begin{aligned}
\avrg{\vect{v}} &:= \frac{1}{2} \left( \vect{v}^+ + \vect{v}^- \right) &\text{ and } &\avrg{\lambda} := \frac{1}{2} \left( \lambda^+ +  \lambda^- \right)  \\
\jump{\vect{v}}& := \vect{v}^+ \cdot \vect{N}^+ + \vect{v}^- \cdot \vect{N}^-  &\text{ and }  &\jump{\lambda} := \lambda^+ \vect{N}^+ + \lambda^- \vect{N}^-
\end{aligned} \right\},
\end{equation}
where $\lambda$ is a scalar, $\vect{v}$ is a vector, the superscripts $\pm$ refer to the element faces on the two sides of each inter-element boundary ($\vect{N}^\pm$ denote the outward unit normal vectors associated with these sides with $\vect{N}^- = - \vect{N}^+$). Based on these definitions, it can also be shown that the following identity holds:
\begin{equation}
\jump{ \lambda \, \vect{v} }=  \avrg{\lambda}  \jump{\vect{v}}  +  \jump{\lambda} \cdot \avrg{\vect{v}}.
\end{equation}
We also introduce the standard compact notation for $L_2$ inner products on a domain or boundary $S$ as follows
\begin{equation}
\innerprod{u}{w}{S} = \int_S u w \ud \tilde{\mu}(S)
\end{equation}
where $\tilde{\mu}(S)$ denotes the standard measure of the corresponding set $S$. Consequently, the integration by parts of functions $\lambda : \RefDom \rightarrow \mathbb{R}$ and $\vect{v}: \RefDom \rightarrow \mathbb{R}^n$ with discontinuities over a inter-element set $\IntOnlyBoun$ reads
\begin{equation}
\begin{split}
\innerprod{\lambda}{\Div \vect{v}}{\RefDom}  + \innerprod{\Grad \lambda }{\vect{v}}{\RefDom} & =   \innerprod{\lambda}{ \vect{v} \cdot \vect{N}}{\ExtBoun}   + \int_{\IntOnlyBoun} \jump{ \lambda \, \vect{v}}  \ud A \\ 
& =   \innerprod{\lambda}{ \vect{v} \cdot \vect{N}}{\ExtBoun}  +\innerprod{\avrg{\lambda}}{\jump{\vect{v}}}{\IntOnlyBoun}  + \innerprod{\jump{\lambda}}{\avrg{\vect{v}}}{\IntOnlyBoun}
\end{split}
\end{equation}
where $\ExtBoun$ denotes the exterior boundary on which the functions are smooth (without including any inter-element boundaries).

\subsection{Semidiscrete weak forms of the governing equations} 
\label{ss:WeakForms}

In this section, we derive weak forms of the governing equations, which were given in \cref{ss:SummaryStrongForms}. First, the standard Galerkin formulation is used to solve the momentum equation since it is a second-order partial differential equation with no higher-order regularity requirements. Second, a C/DG formulation is applied to the fourth-order phase field equation. Finally, we present a more conventional scheme for the solution of the phase-field equation, in which we take advantage of the second-order splitting approach to derive a mixed finite element formulation.

%However, this is not the case for the fourth-order phase-field equation.

\subsubsection{Standard Galerkin formulation of momentum balance equation}

Suppose that $S_\DispRef^h$ is the space in which the approximate displacement solution of \cref{eq:MomBalEqStrongForm} at a given time is sought, where
\begin{equation}
S_\DispRef^h = \left\{ \DispRef_h \, | \, \DispRef_h \in H^1(\RefDom),  \DispRef_h = \thickbar{\vect{u}} \text{ on } \ExtBoun_D^\DispRef \right\},
\end{equation}
with $H^1$ being the standard Sobolev space (that consists of square integrable functions with square integrable derivatives in a weak sense). Let $\vect{w}_h \in V_\DispRef^h$ be the corresponding test function, belonging to the space
\begin{equation}
V_\DispRef^h=\left\{ \vect{w}_h  \, | \, \vect{w}_h \in H^1(\RefDom), \vect{w}_h=\vect{0} \text{ on }  \ExtBoun_D^\DispRef \right\}.
\end{equation}

% $\ExtBoun_{D_i}^Y$ and $\ExtBoun_{N_i}^Y$ with $i \in \mathbb{N}$ and $Y=\{ \DispRef, \PF \}$

Multiplying the strong form by a vector-valued test function $\vect{w}_h$ (virtual displacement) and integrating over the domain $\RefDom$, we arrive at the relation
\begin{equation} \label{eq:WeakFormStep1}
\int_\RefDom \rho_0 \ddt{\DispRef}_h \cdot \vect{w}_h \ud V -\int_\RefDom \Div \tens{P} \cdot \vect{w}_h \ud V - \int_\RefDom \BodyForceRef  \cdot  \vect{w}_h \ud V = 0.
\end{equation}
By partial integration of the second term in \cref{eq:WeakFormStep1}, application of the divergence theorem, and introduction of the traction boundary condition, the weak form of the problem is expressed by
\begin{equation}
R_\DispRef(\vect{w}_h,\motion) = \innerprod{\vect{w}_h}{\rho_0 \ddt{\DispRef}_h}{\RefDom} + \innerprod{\Grad \vect{w}_h}{\tens{P}}{\RefDom} - \innerprod{\vect{w}_h}{\BodyForceRef}{\RefDom} - \innerprod{\tens{w}_h}{\tens{F} \cdot \vect{s}^*}{\ExtBoun_N^\DispRef} = 0.
\end{equation}
The second term can be rewritten using the symmetric property of tensor $\tens{S}$ as follows
\begin{equation}
\tens{P} : \Grad \tens{w}_h = \tens{S} : \left[ \frac{1}{2} \left( \tens{F}\transpose \Grad \vect{w}_h + \Grad\transpose \vect{w}_h \tens{F} \right) \right] = \tens{S} : \delta \tens{E}_h = \frac{1}{2} \tens{S} : \delta \tens{C}_h,
\end{equation}
where $\delta$ denotes the variation of a tensor which is obtained using a directional derivative. For instance, the variation of strain tensor $\tens{E}$ is calculated as follows
\begin{equation}
\delta \tens{E}_h = D \, \tens{E} \cdot \vect{w}_h = \left.  \fullder{ \,\,\,}{\varepsilon} \left[ \tens{E}(\motion + \varepsilon \tens{w}_h) \right]\right\vert_{\varepsilon=0} = \frac{1}{2} \left( \tens{F}\transpose \Grad \vect{w}_h + \Grad\transpose \vect{w}_h \tens{F} \right),
\end{equation}
where $D \, \tens{E} \cdot \vect{w}_h$ denotes the directional derivative of the tensor $\tens{E}$ in the direction $\vect{w}_h$. In this way, the weak form is derived.

%\begin{framedbox}[H]  %htbp
%\caption{Weak form of linear momentum balance equation} 
%\begin{minipage}{0.95\textwidth}
%%\begin{tcolorbox}
Find $\DispRef_h \in S_\DispRef^h \times [0,T]$ such that
\begin{equation}  \label{eq:WeakFormLinearMomBalance}
\left. \begin{aligned}
& R_\DispRef(\vect{w}_h,\motion;\PF) = \innerprod{\vect{w}_h}{\rho_0 \ddt{\DispRef}_h}{\RefDom} + \innerprod{\delta \tens{E}_h}{\tens{S}}{\RefDom} - \innerprod{\vect{w}_h}{\BodyForceRef}{\RefDom} - \innerprod{\tens{w}_h}{\tens{F} \vect{s}^*}{\ExtBoun_N^\DispRef} = 0, \quad \forall \vect{w}_h \in V_\DispRef^h \\
& \innerprod{\vect{w}_h}{\DispRef_h(\vect{X},0)}{\RefDom} = \innerprod{\vect{w}_h}{\DispRef_0(\vect{X})}{\RefDom}, \quad \forall \vect{w}_h \in V_\DispRef^h
\end{aligned} \right\},
\end{equation}
where $\DispRef_0$ denotes the initial displacement, and the stress tensor $\tens{S}=\tens{S}(\tens{C},\PF)$ is calculated from \cref{eq:SecondPKStressFunctionofCdUplus}. The linearization of the variational formulation is presented in \ref{apdx:Linearization}.

\subsubsection{Continuous/Discontinuous Galerkin formulation of the phase-field equation}

The C/DG formulation we present in this section is similar to the one applied strain gradient elasticity problems and to the analysis of thin beams and plates in \cite{engel_continuousdiscontinuous_2002}, as well as the one used for the solution of the Cahn-Hilliard equation in \cite{wells_discontinuous_2006}.

Suppose that the approximate phase-field solution of \cref{eq:PFEquationQuadraticDegrFuncGeneral} is in the space $S_{\PF,k}^h$, defined by
\begin{equation} \label{eq:TrialFunctionsSd1}
S_{\PF,k}^h = \left\{ \PF_h \, | \, \PF_h \in H^1(\RefDom),  \PF_h  \in \mathbb{P}^k(\ElRefDom) \quad \forall e \right\}.
\end{equation}
Here, $\mathbb{P}^k(\ElRefDom)$ denotes the space of polynomials of degree at most $k$ with real coefficients on each element $\ElRefDom$. Let $c_h \in V_{\PF,k}^h$ be the corresponding test function, belonging to the space
\begin{equation}  \label{eq:TestFunctionsSd1}
V_{\PF,k}^h=\left\{ c_h  \, | \, c_h \in H^1(\RefDom), c_h  \in \mathbb{P}^k(\ElRefDom) \quad \forall e \right\},
\end{equation}
Although the fourth-order phase-field equation mandates the use of a space of approximate solutions in $H^2$ (with the standard continuous Galerkin method), the order of regularity in $S_{\PF,k}^h $ is reduced. In this continuous/discontinuous Galerkin method, the functions in the approximate solution space are continuous, but their first- and higher-order derivatives are discontinuous. The continuity requirements on the derivatives are imposed weakly, i.e.\ with the aid of terms added to the weak form to enforce these requirements on inter-element boundaries. The weak form of this formulation is summarized as follows.

%\begin{framedbox}[H]  %htbp
%\caption{Weak form of the fourth-order phase-field equation in C/DG formulation} 
%\begin{minipage}{0.95\textwidth}
%%\begin{tcolorbox}  
Find $\PF_h \in S_{\PF,k}^h \times [0,T]$ such that
\begin{equation}  \label{eq:CDGWeakFormCDGMethod}
\begin{split}
 R^{(1)}_\PF(c_h,\PF_h;\mathcal{H}) &= \innerprod{c_h}{\MicroDens  \ddt{\PF}_h}{\RefDom} +  \innerprod{c_h}{\MicroDamp  \dt{\PF}_h}{\RefDom}  + \innerprod{\Delta c_h }{\alpha_2 \PFlap_h}{\TildeRefDom}  - \innerprod{\nabla c_h }{ \alpha_1 \PFgrad_h }{\RefDom} \\ 
& + \innerprod{c_h}{\alpha_0  \PF_h}{\RefDom} + \innerprod{c_h}{ g^\prime(\PF_h) \mathcal{H}}{\RefDom}  - \innerprod{\jump{\nabla c_h}}{\alpha_2 \avrg{\PFlap_h}}{\IntBoun} \\ 
& - \innerprod{\avrg{\Delta c_h}}{\alpha_2 \jump{\PFgrad_h}}{\IntBoun} + \innerprod{\frac{\beta_{s2}}{\avrg{h_e}} \jump{\nabla  c_h}}{ \jump{\PFgrad_h}}{\IntBoun}  = 0 , \quad \forall c_h \in V_{\PF,k}^h,
\end{split}
\end{equation}
where $\beta_{s2}$ is a penalty parameter and $h_e$ is an element length-scale. The initial condition of the phase-field equation, $\PF_0$, is imposed by
\begin{equation}
\innerprod{c_h}{\PF_h(\vect{X},0)}{\RefDom} = \innerprod{c_h}{\PF_0(\vect{X})}{\RefDom}, \quad \forall c_h \in V_{\PF,k}^h  \, .
\end{equation}
%%\end{tcolorbox}
%\end{minipage}
%\end{framedbox}

%\MicroDens  \ddt{\PF} + \alpha_3  \Delta \left( \PFlap \right) +  \alpha_2  \PFlap  +  \alpha_1  \PF   = 2 \mathcal{H} \text{ in } \RefDom

%\TotBoun = \ExtBoun \cup \IntBoun \, .

The interior penalty term in \cref{eq:CDGWeakFormCDGMethod} plays an important role in the weak imposition of solution regularity. The penalty parameter $\beta_{s2}$  is a function of $\alpha_2$ (and consequently of the material parameters $G_c$ and $\lzr$). We also note that \cref{eq:CDGWeakFormCDGMethod} is a symmetric formulation. The reader is referred to \ref{apdx:CDG_Consistency} for proofs of consistency of the C/DG method and dependence on $\alpha_2$.

% depends on an element length-scale $h_e$ and penalty parameter $\beta_{s2}$.

\begin{remark}
The C/DG formulation involves not only the standard FEM integration over element domains but also integration over inter-element boundaries.
\end{remark}

\subsubsection{Mixed finite element formulation of the phase-field equation}

%A common strategy to solve a fourth-order equation is to replace it by two second-order equations and 
In this section, we begin with the phase-field formulation based on the splitting approach, and described by the coupled second-order system in \cref{eq:SecondOrderSplittingSystems}. We use the standard finite element method for each equation separately. To this end, we introduce the following function spaces in addition to the spaces ($S_{\PF,k}^h$ and $V_{\PF,k}^h$) which were defined in \cref{eq:TrialFunctionsSd1,eq:TestFunctionsSd1}:
\begin{equation} \label{eq:TrialTestFunctionsSpsi1}
S_{\psi,k}^h = V_{\psi,k}^h = \left\{ \psi_h \, | \, \psi_h \in H^1(\RefDom),  \psi_h  \in \mathbb{P}^k(\ElRefDom) \quad \forall e, \,  \psi^h = 0 \text{ on } \partial \RefDom \right\}
\end{equation}
% \label{eq:TestFunctionsSpsi1}
The approximate phase-field solution $\PF_h$ to \cref{eq:SecondOrderSplittingSystems} is sought in space $S_{\PF,k}^h$ and its approximate scaled Laplacian $\psi_h$ in $S_{\psi,k}^h$. The corresponding test functions $c_h$ and $\chi_h$ respectively belong to spaces $V_{\PF,k}^h$ and $V_{\psi,k}^h$. The weak form of this formulation is summarized as follows.

%\begin{equation} \label{eq:TrialTestFunctionsSpsi1}
%\left. \begin{aligned}
%S_{\psi,k}^h & = \left\{ \psi_h \, | \, \psi_h \in H^1(\RefDom),  \psi_h  \in \mathbb{P}^k(\ElRefDom) \quad \forall e, \,  \psi^h = 0 \text{ on } \partial \RefDom \right\} \\
%V_{\psi,k}^h & = \left\{ \chi_h  \, | \, \chi_h \in H^1(\RefDom), \chi_h  \in \mathbb{P}^k(\ElRefDom) \quad \forall e, \, \chi^h = 0 \text{ on } \partial \RefDom \right\}
%\end{aligned} \right\}.
%\end{equation}

Find $\PF_h \in S_{\PF,k}^h \times [0,T]$ and $\psi_h \in S_{\psi,k}^h \times [0,T]$ such that
\begin{equation} \label{eq:MixedFEMPFWeakForm}
\left. \begin{aligned}
  R^{(2)}_\PF(c_h,\PF_h) & =   \innerprod{\nabla c_h }{\PFgrad_h }{\RefDom}  + \lambda_0 \innerprod{c_h}{\psi_h}{\RefDom}  + \frac{\alpha_2}{\alpha_1} \lambda_0 \innerprod{c_h}{\nabla \psi_h \cdot \vect{N}}{\ExtBoun_{R_2}^\PF } = 0 , \quad \forall c_h \in V_{\PF,k}^h  \\
  R^{(2)}_\psi(\chi_h,\psi_h;\mathcal{H}) & = \innerprod{\chi_h}{\MicroDens  \ddt{\PF}_h}{\RefDom} +  \innerprod{\chi_h}{\MicroDamp  \dt{\PF}_h}{\RefDom}  - \lambda_0 \alpha_2  \innerprod{\nabla \chi_h }{\nabla \psi_h }{\RefDom}  \\
  & + \lambda_0  \alpha_1 \innerprod{\chi_h}{\psi_h}{\RefDom} + \alpha_0 \innerprod{\chi_h}{\PF_h}{\RefDom}   + \innerprod{\chi_h}{ g^\prime(\PF_h) \mathcal{H}}{\RefDom}  = 0 , \quad \forall \chi_h \in V_{\psi,k}^h \\
\end{aligned} \right\},
\end{equation}
and the corresponding initial condition $\PF=\PF_0$, is imposed by
\begin{equation}
\innerprod{c_h}{\PF_h(\vect{X},0)}{\RefDom}  = \innerprod{c_h}{\PF_0(\vect{X})}{\RefDom}, \quad \forall c_h \in V_{\PF,k}^h.
\end{equation}

\begin{remark}
Although the function space $S_{\PF,k}^h$, defined in \cref{eq:TrialFunctionsSd1}, is the same in both mixed and C/DG formulations, the C/DG formulation involves second-order derivatives while the mixed FEM formulation involves only first-order derivatives. Hence, there are different low bounds of polynomial degree $k$ that can be used. Namely, the bounds are given for C/DG formulation: $k \ge 2$ and for mixed formulation: $k \ge 1$. 
\end{remark}

\begin{remark}
The stability of mixed methods is conditional, depending on the discretization of the different fields. Although we do not prove the stability of the proposed mixed schemes, no spurious oscillations or other instabilities have been observed in the numerical results. Furthermore, similar mixed finite element methods have been used successfully elsewhere in the literature, e.g.\ to treat the Cahn-Hilliard equation \cite{elliott1989second,diegel2016stability} which is very similar in form to the fourth-order phase-field equation considered here. 
\end{remark}

%\kgc{Remarks on LBB?}

\subsection{Numerical implementation: staggered scheme \& time integration}
\label{ss:NumericalImplementStagTimeInt}
Considering the previous discrete weak forms, we derive two techniques to solve the governing equations. In both techniques, the momentum balance is imposed with the same weak form as expressed in \cref{eq:WeakFormLinearMomBalance}. However, the distinction between the proposed techniques stems from the way we treat the fourth-order phase-field equation. The novel technique utilizes the continuous/discontinuous Galekrin weak form in \cref{eq:CDGWeakFormCDGMethod} and is compared with the more conventional mixed finite element method based on the weak form \eqref{eq:MixedFEMPFWeakForm}. The techniques are briefly outlined as follows.

% of coupled hyper-elastic fracture mechanics problems

\begin{table}[H]
\centering
\caption{Different techniques to solve the governing equations of coupled hyper-elastic fracture mechanics problems}
\label{tbl:StrategiesHyperElasticFractureMechanicsProbs}
%\resizebox{0.8\textwidth}{!}{%
\begin{tabular}{|l|c||l|c|}
\hline
      \multicolumn{2}{|c||}{Momentum equation} &
      \multicolumn{2}{c|}{Phase-field equation}  \\
      \hline
      Method & Weak forms & Method & Weak forms  \Tstrut\Bstrut\\
\hline
Standard FEM & \cref{eq:WeakFormLinearMomBalance} & C/DG Method & \cref{eq:CDGWeakFormCDGMethod}   \Tstrut\\
Standard FEM & \cref{eq:WeakFormLinearMomBalance} & Mixed FEM & \cref{eq:MixedFEMPFWeakForm}   \Tstrut\Bstrut\\
\hline
\end{tabular}
%}
\end{table}

\subsubsection{Staggered scheme}
For each technique, a staggered time integration scheme (which introduces a weak coupling between the governing equations) is used. Following a similar approach to the one presented in \cite{larsen2010existence,bourdin2011time}, governing equations are solved independently at each time step. Specifically, we first find the solution to the momentum equation. Once nodal displacements are known, the history variable $\mathcal{H}$ is updated based on the strain measures computed at each Gauss point. Finally, the phase-field equation is solved using the calculated history variable. We repeat the aforementioned process for each time step.

The staggered scheme provides the flexibility of adopting different time integration schemes for each governing equation. Phase-field methods typically require highly refined meshes to capture the steep gradients in the neighborhood of cracks. Hence, explicit time integration schemes may introduce significant limitations on the size of time steps due to the CFL conditions. 

For second-order phase-field models, an interesting approach is presented in \cite{kamensky2018hyperbolic,MOUTSANIDIS2018114}, where damage evolution is governed by a hyperbolic partial differential equation with less restrictive CFL conditions compared with the parabolic equations. For fourth-order phase-field models, these conditions may be more prohibitive to apply explicit schemes. A thorough investigation of alternative strategies to improve convergence conditions is left for future work.

%\kgc{The phase field equation here is second-order in time, which makes it hyperbolic, doesn't it? Also, Refs \cite{kamensky2018hyperbolic,MOUTSANIDIS2018114} claim to use microforce balance.}

\begin{remark}
A major advantage of the quadratic degradation function given in \cref{eq:QuadraticDegrFunction} is that it causes the phase-field equation to be linear, obviating the need for an iterative procedure and simplifying the solution in a staggered scheme. 
\end{remark}

\subsubsection{Time integration}

In this paper, we use only implicit time integration schemes because CFL conditions inhibit an efficient parallel scaling performance when running on large numbers of processors \cite{keyes2006implicit}. Specifically, the implicit Newmark method with parameters $\beta \in (0,0.5]$ and $\gamma \in (0,1] $ is employed to advance the momentum equation in time. In this scheme, accelerations and velocities at time $t^{(n+1)}$ are approximated as follows
\begin{equation} \label{eq:NewmarkScheme}
\left. \begin{aligned}
& \accel^{(n+1)}_h = \frac{1}{\beta (\Delta t)^2} \left( \disp^{(n+1)}_h - \disp^{(n)}_h  \right) - \frac{1}{\beta \Delta t} \velo^{(n)}_h - \frac{1-2\beta}{2 \beta} \accel^{(n)}_h \\
& \velo^{(n+1)}_h = \frac{\gamma}{\beta \Delta t}  \left( \disp^{(n+1)}_h - \disp^{(n)}_h  \right) + \left( 1 - \frac{\gamma}{\beta} \right) \velo^{(n)}_h +  \left( 1 - \frac{\gamma}{2 \beta} \right) \accel^{(n)}_h  \\
\end{aligned} \right\},
\end{equation}
where $\Delta t$ denotes the current time step. A fully discrete nonlinear system for the displacements $\disp^{(n+1)}_h$ can be derived after the substitution of \cref{eq:NewmarkScheme} into the weak form of momentum equation in \cref{eq:WeakFormLinearMomBalance}. Finally, we assume that micro-inertia and micro-damping effects are negligible (i.e., $\MicroDens = \MicroDamp = 0$). The aforementioned assumption leads to an elliptic phase-field equation with no need of time discretization.

%In this part, we discretize the momentum balance in time. The total time interval $[0,T]$ is divided into $N_T$ subintervals with time step $\Delta t = \frac{T}{N_T}$ and time instants $t^{(n)}=n \Delta t$.  In this implicit scheme, 

%\subsection{Spatial discretization method}

%\subsection{Time integration method}

%\subsection{Linearizations}

\section{Numerical results}
\label{sec:Results}

In this section, the numerical performance of the proposed C/DG method is compared to the conventional mixed FEM to model fracture. Specifically, three staggered schemes are examined. For all schemes, the momentum balance equation is solved for displacements employing the standard FEM with four-node quadrilateral elements [Q4] (i.e., the standard bilinear shape functions). However, three different methods are used to solve the fourth-order phase-field equation: (i) the proposed continuous/discontinuous Galerkin method with nine-node quadrilateral elements (Q9), (ii) second-order mixed FEM with nine-node quadrilateral elements (Q9Q9), and (iii) standard mixed FEM with four-node quadrilateral elements (Q4Q4). Here, the element type is written twice in parentheses for mixed FEMs since the two unknown fields (namely $\PF$ and $\psi$) are approximated using the same shape functions, while the C/DG method uses the shape functions to approximate only the primitive variable, $\PF$. The numerical schemes are abbreviated using the corresponding element types as shown in \cref{tbl:MethodNotationsSchemes}. Note that the mixed FEM [Q4](Q4Q4) is the most conventional technique used in the literature to solve a fourth-order equation, by first splitting it into a system of two second-order equations. The mixed FEM [Q4](Q9Q9) is included in this study for direct comparison with the C/DG method [Q4](Q9), since both of these schemes employ bi-quadratic interpolation of the phase field at the element level. 

\begin{table}[b]
\centering
\caption{Numerical schemes for solving the governing equations with a staggered approach}
\label{tbl:MethodNotationsSchemes}
%\resizebox{0.8\textwidth}{!}{%
\begin{tabular}{|l|c||l|c||l|}
\hline
      \multicolumn{2}{|c||}{Momentum Equation} &
      \multicolumn{2}{c||}{Phase-Field Equation} &
      Combined Notation  \\
      \hline
       Method & Element Type & Method & Element Type & Staggered Approach  \Tstrut\Bstrut\\
\hline
Standard FEM & [Q4] & C/DG Method & (Q9) & Scheme [Q4](Q9)   \Tstrut\\
Standard FEM & [Q4] & Second-order Mixed FEM & (Q9Q9) &  Scheme [Q4](Q9Q9) \Tstrut\\
Standard FEM & [Q4] & Mixed FEM & (Q4Q4) & Scheme [Q4](Q4Q4)    \Tstrut\Bstrut\\
\hline
\end{tabular}
%}
\end{table}

The numerical performance of the aforementioned schemes is studied with the aid of two dynamic fracture problems: (i) a shear experiment and (ii) a double cantilever beam experiment. For each example problem, two congruent meshes are constructed in order to compute the solution using the three proposed numerical schemes: one mesh consisting of Q4 elements and the second comprised of Q9 elements. Problem-specific details are provided below regarding the geometry, boundary conditions, and discretization schemes used in each problem. Material parameters are chosen based on numerical examples previously presented in the literature. First, we study the influence of the penalty term on the phase-field solution by changing the penalty parameter $\beta_{s2}$, and observe that there is a range of $\beta_{s2}$ in which numerical results are insensitive to the value of this parameter. After the selection of a value in the aforementioned range, we compare the numerical solutions obtained using the schemes outlined in \cref{tbl:MethodNotationsSchemes}. 
%In addition, the influence of the penalty term on the proposed formulation is studied through changing the penalty parameter $\beta_{s2}$. Finally, a scalability analysis is run to present the efficiency of the C/DG method compared with the standard mixed formulation.

For the purpose of implementing the three computational schemes described above, a MOOSE-based application was developed using the Tensor Mechanics module \cite{permann2020moose}. The staggered time-integration scheme is achieved using the MultiApp and Transfer systems \cite{gaston2015physics}. All the simulations presented are carried out on a computer cluster with Intel Xeon CPUs (E5-2695) running at 2.10GHz.

% \cite{permann2020moose,rovinelli2020identify}

% the automatic differentiation capabilities \cite{lindsay2021automatic}

% First Example

\subsection{Dynamic shear problem}

\subsubsection{Problem description and discretization}

\begin{figure}[t] 
\centering
\subfloat[\label{fig:Geometry_Ex1}]{%
\includegraphics[width=0.35\textwidth]{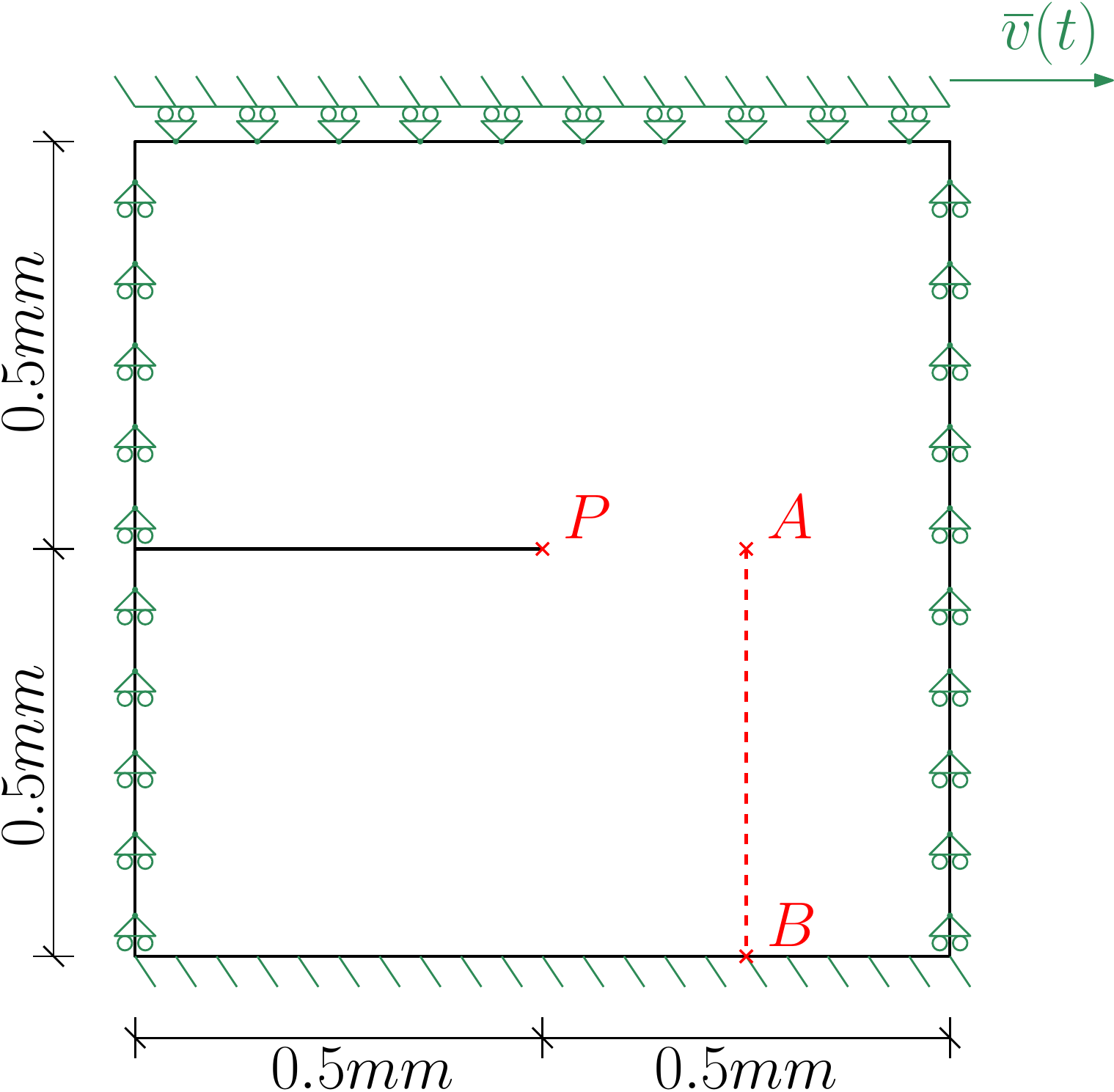}  %0.29
}
\hspace{1.0cm}
\subfloat[\label{fig:Load_Ex1}]{%
\includegraphics[width=0.35\textwidth]{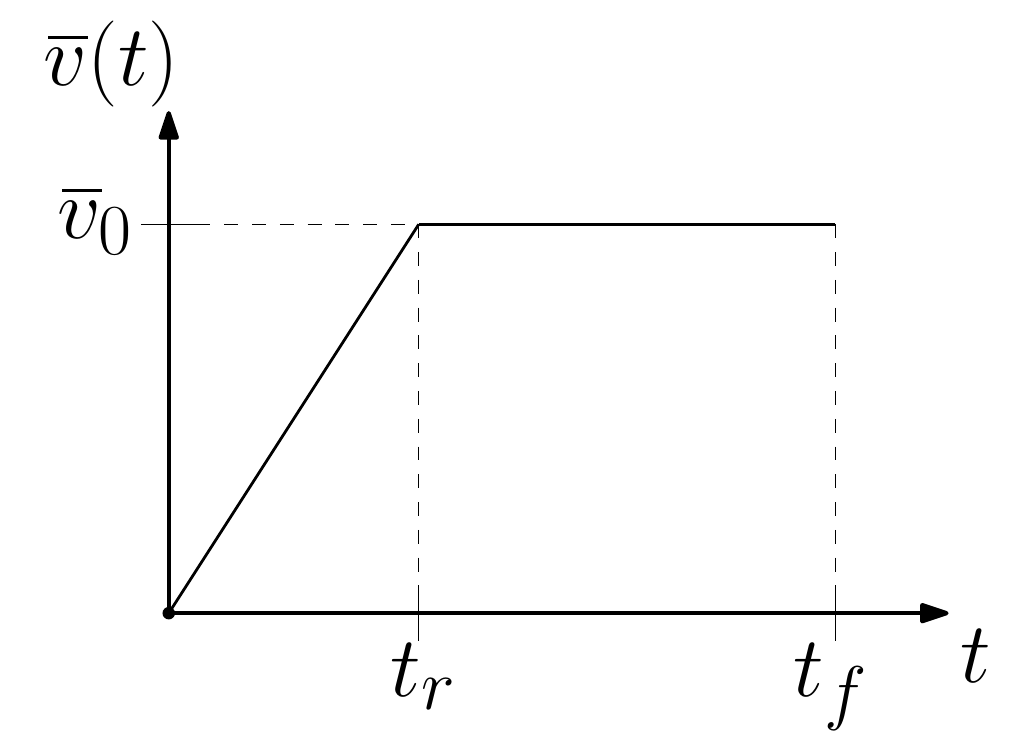}
}
\captionsetup{width=.90\textwidth}
\caption{ (a) Geometry and boundary conditions of a square plate under shear deformation, (b) The time-dependent velocity field, prescribed on the upper side of the plate.}
\label{fig:GeometryLoad_Ex1}
\end{figure}

In this example, a displacement-controlled shear experiment of a square plate is conducted. The geometry and boundary conditions of the problem are depicted in \cref{fig:Geometry_Ex1}. The plate is modeled as a two-dimensional square with $1mm$ side length under plane strain conditions. The plate has an initial crack, represented by the notch (discontinuity) extending from the left side of the plate to point $P$. The bottom side of the plate is completely fixed (i.e., displacements set to zero in the x- and y-directions), while permitting displacements in the x-direction only on the three other sides. A shear deformation is imposed by prescribing, everywhere on the upper side of the plate, a velocity $\thickbar{v}(t)$ which depends on time as shown in \cref{fig:Load_Ex1}, with ramp time $t_r=25.0\times10^{-6}$~s, final time $t_f=65.0\times10^{-6}$~s and maximum velocity $\thickbar{v}_0=260$~mm/s.

%while they are fixed only vertically  (i.e., modeled as rollers). 
%and restricting it in the vertical direction. The time-dependent horizontal displacement field has a ramp form as shown in \cref{fig:Load_Ex1}.  

\begin{figure}[t] 
\centering
\subfloat[\label{fig:MeshQ4_Ex1} (Q4) mesh]{%
\includegraphics[width=0.35\textwidth]{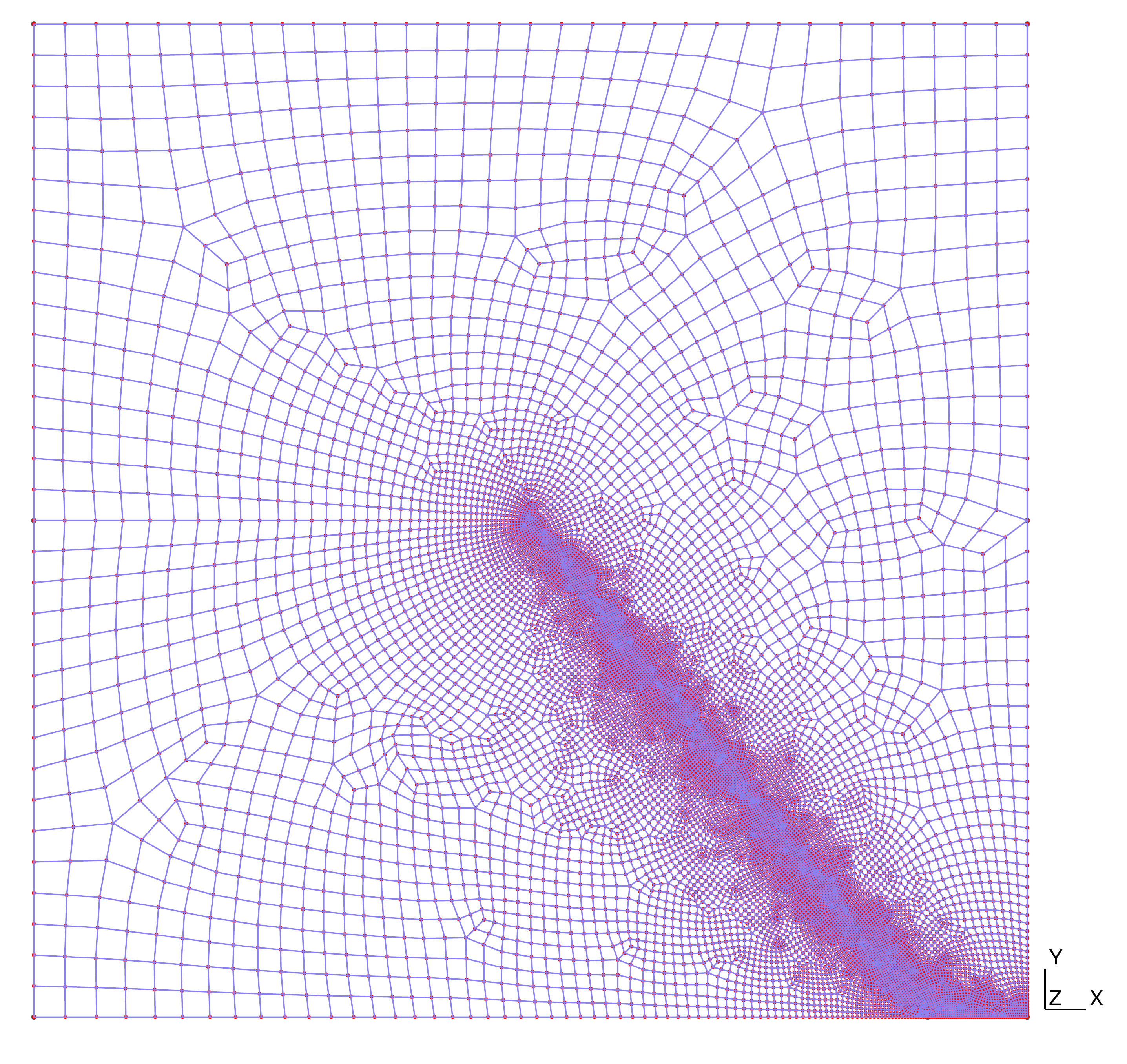}  %0.29
}
\hspace{1.0cm}
\subfloat[\label{fig:MeshQ9_Ex1} (Q9) mesh]{%
\includegraphics[width=0.35\textwidth]{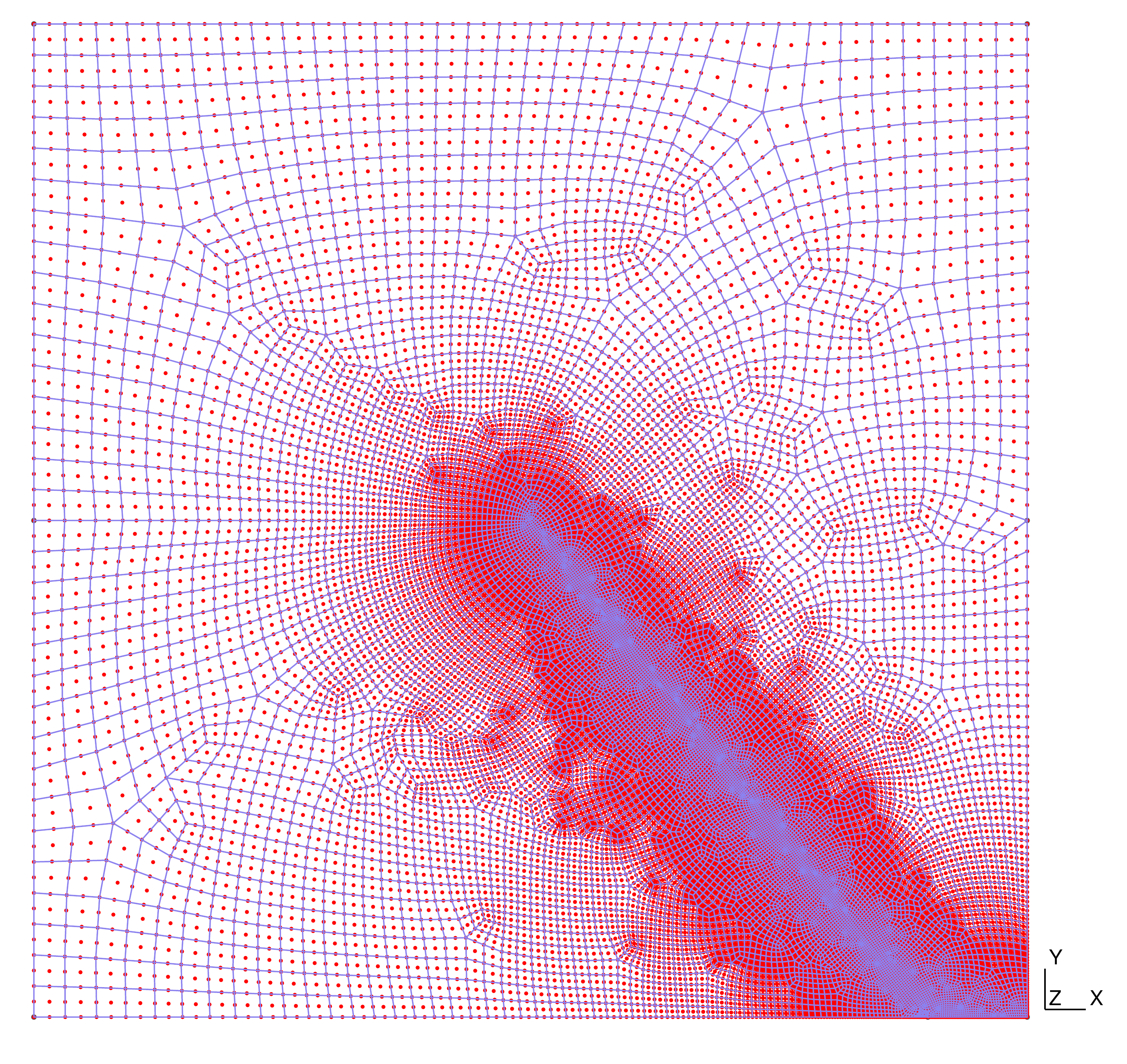}
}
\captionsetup{width=.90\textwidth}
\caption{Unstructured meshes in this problem, refined along the anticipated crack path: (a) four-node quadrilateral elements, and (b) nine-node quadrilateral elements.  }
\label{fig:MeshesQ4andQ9_Ex1}
\end{figure}

In this problem, two unstructured meshes are constructed to discretize the domain of the plate: (i) the mesh shown in \cref{fig:MeshQ4_Ex1}, consisting of 11,948 four-node quadrilateral (Q4) elements and 12,412 nodes, and (ii) the mesh in \cref{fig:MeshQ9_Ex1}, composed of 11,948 nine-node quadrilateral elements (Q9) and 48,718 nodes. Both meshes are highly refined around the area where the crack is expected to propagate. Solution of the linear momentum balance equation involves two unknown displacements per node, while the phase-field equation involves either one unknown per node in the case of the C/DG method or two unknowns per node in the case of mixed FEM. The total number of degrees of freedom (DoFs) for each scheme are 72,371~DoFs for C/DG Scheme [Q4](Q9), 120,504~DoFs for mixed FEM Scheme [Q4](Q9Q9), and 48,476~DOFs for mixed FEM Scheme [Q4](Q4Q4). Although Scheme [Q4](Q4Q4) involves the smallest number of DoFs, it is also the least accurate since it employs interpolation functions of the lowest order (bi-linear). On the other hand, Scheme [Q4](Q9Q9) uses bi-quadratic interpolation of the phase field, and as a result, the number of DoFs increases significantly. Despite the higher accuracy of this scheme, its elevated computational cost can be a hurdle in large problems. It is noted that despite using bi-quadratic interpolation functions for the phase field, The C/DG Scheme [Q4](Q9) involves a much smaller number of DoFs compared to Scheme [Q4](Q9Q9). This is because the C/DG scheme enforces the necessary continuity conditions by adding variational/penalty terms to the weak form, instead of introducing additional auxiliary fields.

%reported in \cref{tbl:DoFsPerSchemes}. 

%\begin{table}[H]
%\centering
%\caption{The total number of Degrees of Freedom for each scheme}
%\label{tbl:DoFsPerSchemes}
%%\resizebox{0.3\textwidth}{!}{%
%\begin{tabular}{|l|r|}
%\hline
%Numerical Schemes  &  DoFs  \Tstrut\Bstrut\\
%\hline
%Scheme  [Q4](Q9) & 72,371  \Tstrut\\
%Scheme [Q4](Q9Q9) & 120,504   \Tstrut\\
%Scheme [Q4](Q4Q4)  & 48,476   \Tstrut\Bstrut\\  
%\hline
%\end{tabular}
%%}
%\end{table}

%For the C/DG formulation a higher order mesh is required to capture the higher order derivatives. 
%Finally, the standard mixed formulation is solved on (i) a (Q4Q4) mesh with 23,896 elements and 24,824 nodes and (ii) a (Q9Q9) mesh with 23,896 elements and 97,436 nodes.
%The latter mesh is selected because the solution has the same convergence rate with the C/DG method.

%Units: mm, s, N, MPa, Mg, K, mJ
\begin{table}[b]
\centering
\caption{Material properties and simulation parameters}
\label{tbl:simulation_parameters_Ex1}
\resizebox{0.55\textwidth}{!}{%
\begin{tabular}{|lccc|}
\hline
Property & Notation  &  Value  &  Unit  \Tstrut\Bstrut\\
\hline
Young's Modulus & $E$ & 210.0E+3 & $MPa$    \Tstrut\\
Poisson's Ratio & $\nu$ & 0.3 &  -  \Tstrut\\
Mass Density & $\rho$ & 8.0E-9 & $Mg$ ${(mm)}^{-3}$  \Tstrut\\
Critical Energy Release Rate & $G_c$ & 2.7 & $mJ$ ${{(mm)}^{-2}}$ \Tstrut\\
Length Scale & $\lzr$ & 3.75E-3 & $mm$ \Tstrut\\
Residual Stiffness Parameter & $\eta_0$ & 1.0E-6 & - \Tstrut\\
Newmark beta & $\beta$ & 0.3025 & - \Tstrut\\
Newmark gamma & $\gamma$ & 0.6 & - \Tstrut\Bstrut\\
\hline
\end{tabular}
}
\end{table}

The material and simulation parameters used in this problem are provided in \cref{tbl:simulation_parameters_Ex1}. Young's modulus, Poisson's ratio and the mass density are assigned values representative of stainless steel. The values of the phase-field fracture parameters are on the same order of magnitude as those found in the literature. The Newmark parameters are chosen so as to damp higher frequency modes without negatively affecting the stability of the numerical integration scheme. 

\subsubsection{Sensitivity analysis: Influence of the penalty parameter}
\label{sss:SensitivityAnalysisEx1}

\begin{figure}[t] %[H] % "[t!]" placement specifier just for this example
\centering
\hspace{3.0cm}
\includegraphics[width=0.35\textwidth]{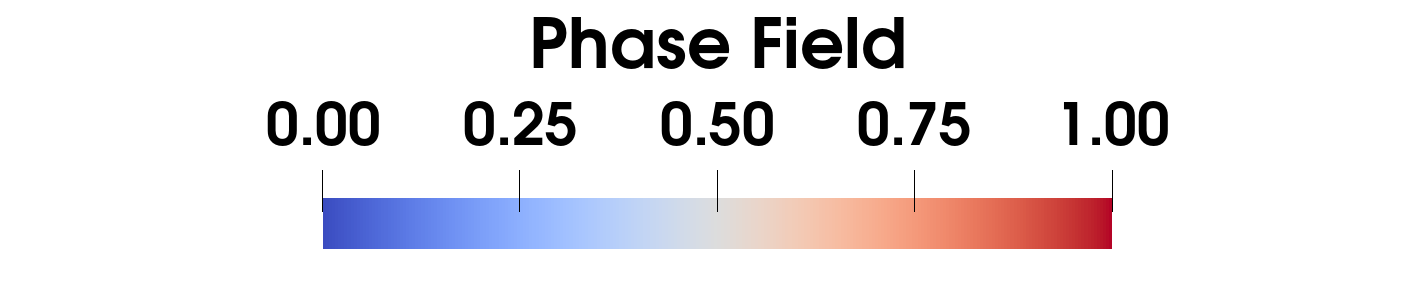}
\hspace{3.0cm}

\vspace{-0.35cm}

\centering

\subfloat[\label{fig:CDG_beta1_Ex1}  $\beta_{s2} = 5.0\times10^{-5}$  ]{%
\includegraphics[width=0.30\textwidth]{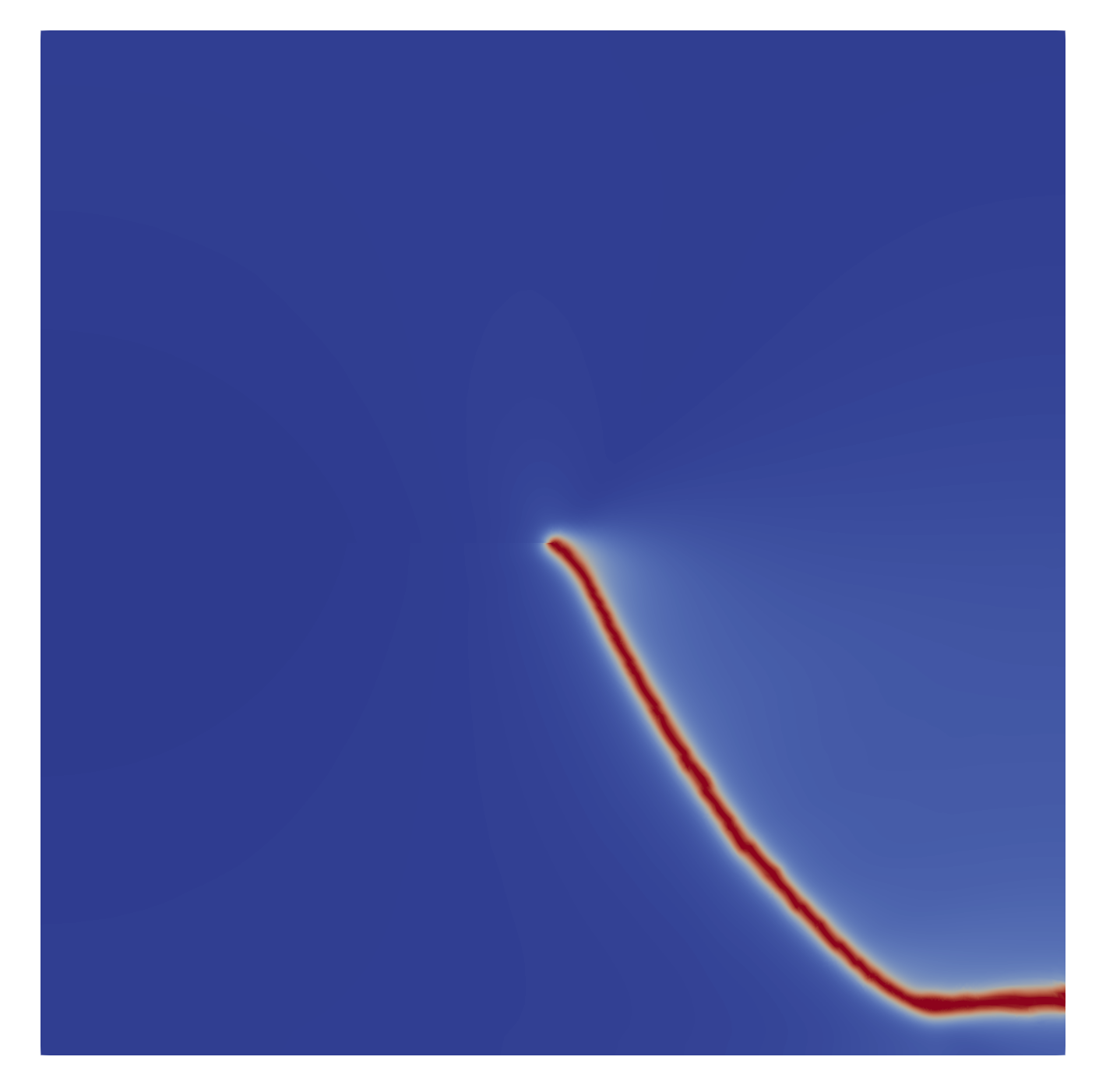}
}
\hspace{1.0cm}
\subfloat[\label{fig:CDG_beta2_Ex1}  $\beta_{s2} = 20.0\times10^{-5}$]{%
\includegraphics[width=0.30\textwidth]{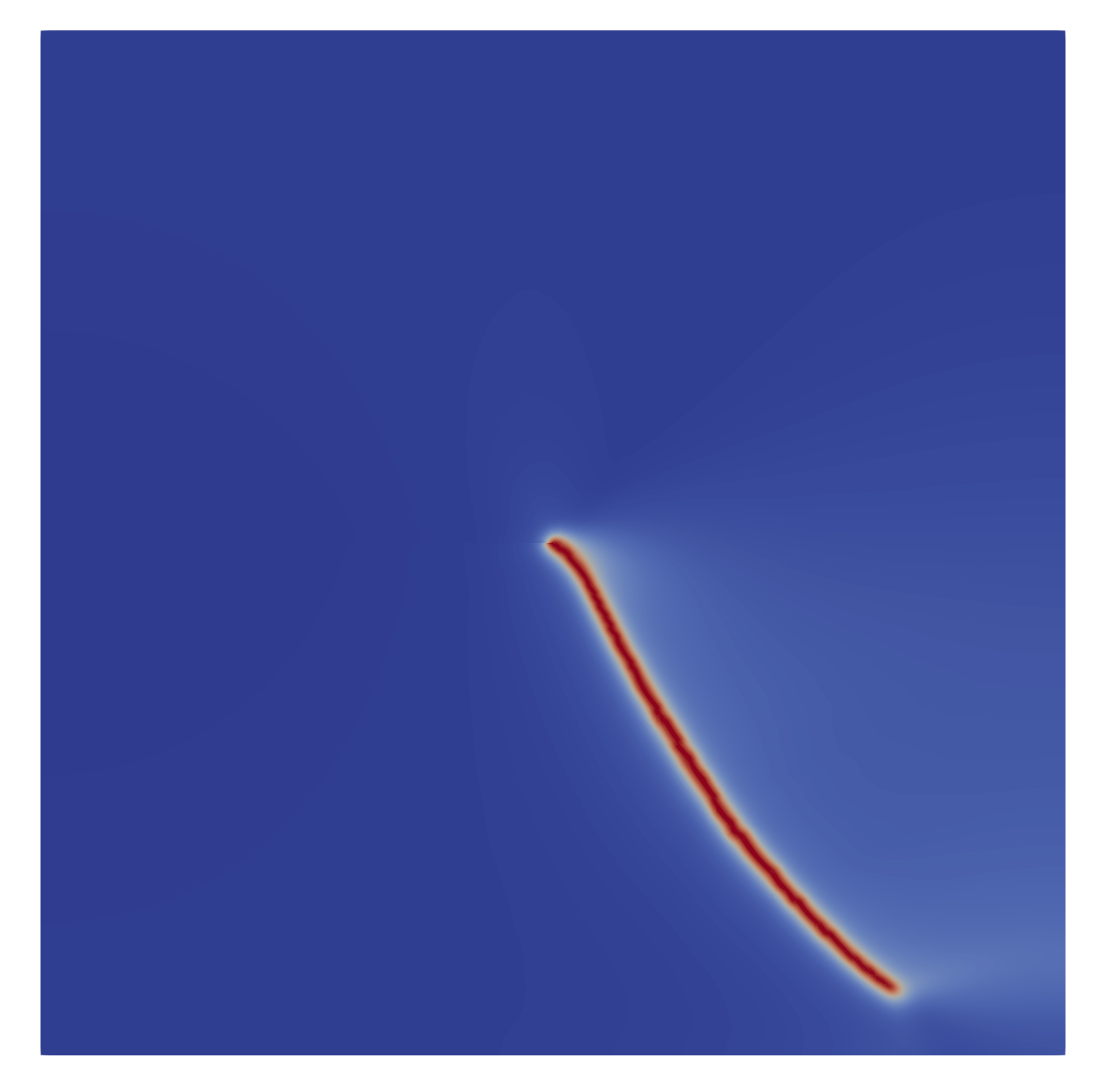}  %0.29
}

\centering

%\medskip
\vspace{-0.2cm}

\subfloat[\label{fig:CDG_beta3_Ex1}  $\beta_{s2} = 35.0\times10^{-5}$  ]{%
\includegraphics[width=0.30\textwidth]{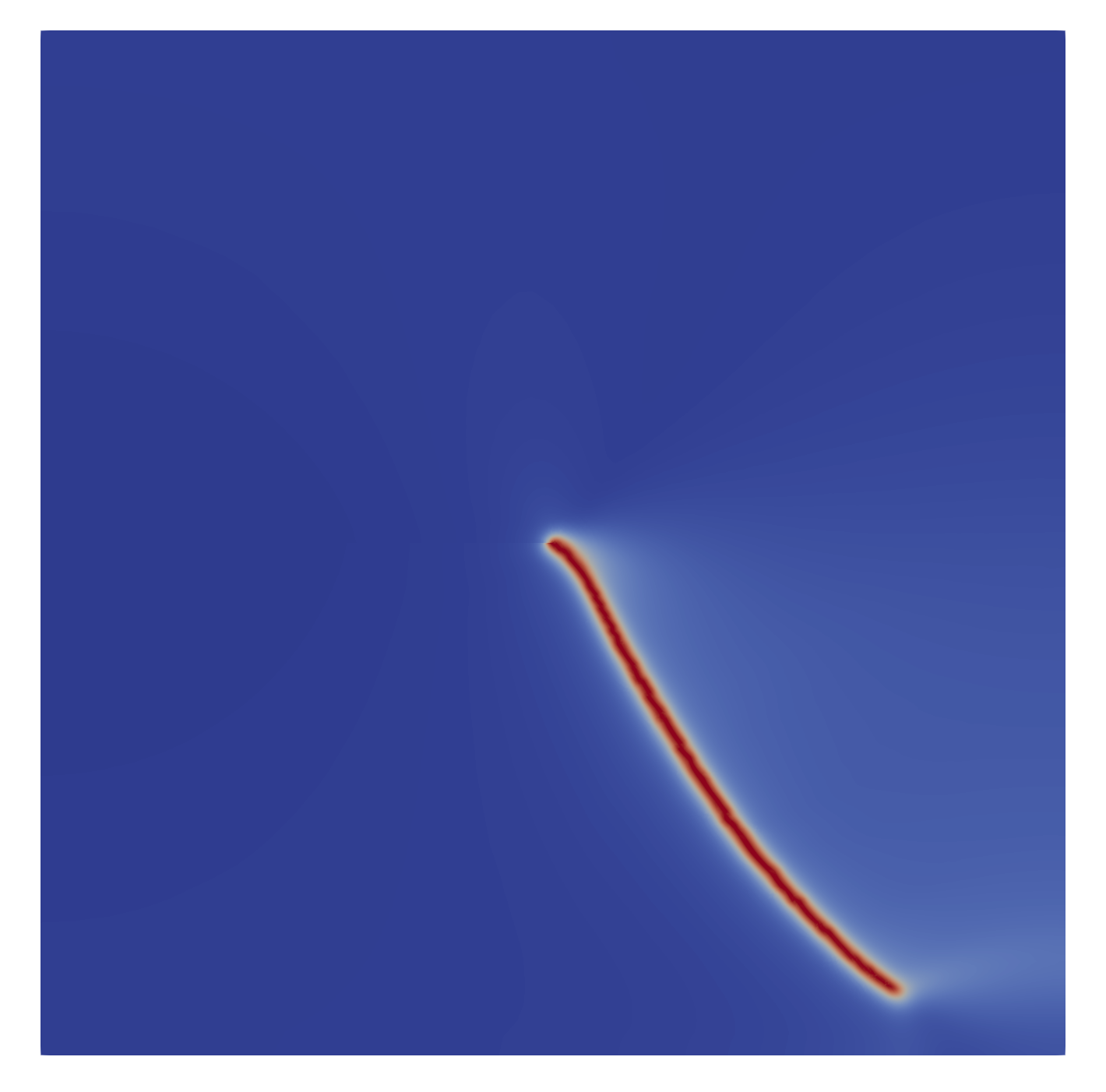}
}
\hspace{1.0cm}
\subfloat[\label{fig:CDG_beta4_Ex1}  $\beta_{s2} = 50.0\times10^{-5}$]{%
\includegraphics[width=0.30\textwidth]{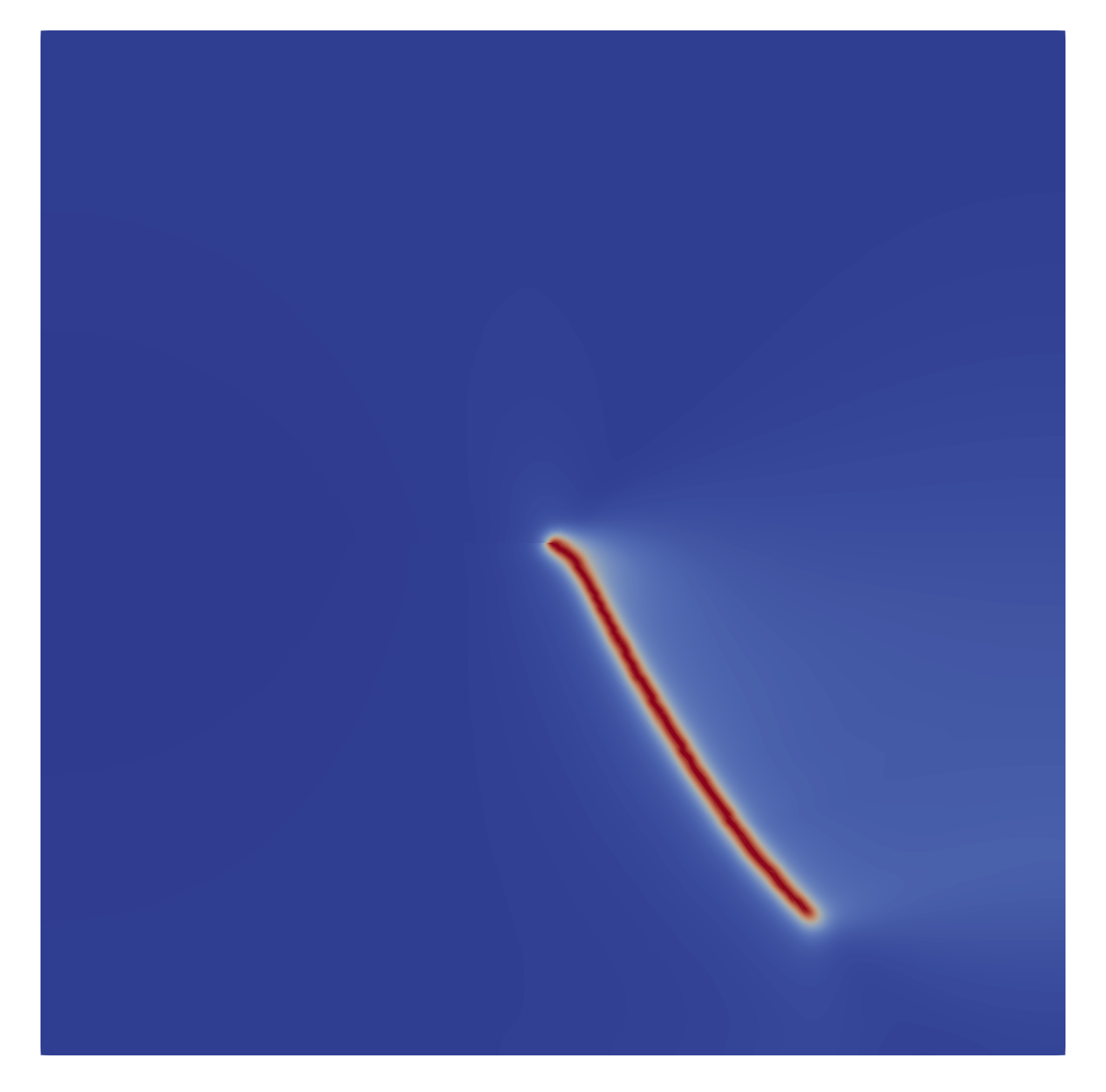}  %0.29
}

\captionsetup{width=.90\textwidth}
\caption{ A penalty-parameter study: phase-field fracture paths at $t=63.0\times10^{-6}$~s are depicted when different penalty parameters in Scheme \text{[Q4](Q9)} are used. } \label{fig:AllMethodsPFbetas_Ex1}
\end{figure}

A sensitivity analysis is required to examine the influence of the penalty parameter $\beta_{s2}$ on the phase-field solution when Scheme \text{[Q4](Q9)} is employed. The proper order of the penalty parameter cannot be estimated \textit{a priori} because it depends on the mesh, element type, and material parameters. Generally, the penalty parameter should be assigned the smallest value which sufficiently enforces boundary and regularity conditions, leading to stable solutions. Assigning excessively large values to the penalty parameter do not improve the accuracy, and may lead to numerical difficulties (due to ill-conditioning).

%LegendConfig = { Pr1M4F:   '[Q4](Q4Q4)',
%                 Pr1M9S:   '[Q4](Q9Q9)',
%                Pr1M9Pf1: '[Q4](Q9): $\\beta_{s2} = 5.0\times10^{-5}$',
%                 Pr1M9Pf2: '[Q4](Q9): $\\beta_{s2} = 20.0\times10^{-5}$',
%                 Pr1M9Pf3: '[Q4](Q9): $\\beta_{s2} = 35.0\times10^{-5}$',
%                 Pr1M9Pf4: '[Q4](Q9): $\\beta_{s2} = 50.0\times10^{-5}$'
%                }

The phase-field solutions computed using the C/DG scheme with different penalty parameter values, representing the predicted crack path at $t=63.0\times10^{-6}$~s, are shown in \cref{fig:AllMethodsPFbetas_Ex1}. For $\beta_{s2} = 5.0\times10^{-5}$, the crack path is fully developed as shown in \cref{fig:CDG_beta1_Ex1} but we have no evidence that this value is sufficient to accurately impose boundary and regularity conditions. The similar crack paths shown in \cref{fig:CDG_beta2_Ex1,fig:CDG_beta3_Ex1} demonstrate insensitivity to the penalty parameter in the range $\beta_{s2} \in [20.0\times10^{-5}, 35.0\times10^{-5}]$. Finally, for a larger value of the penalty parameter ($\beta_{s2} = 50.0\times10^{-5}$) the solution exhibits a time lag (\cref{fig:CDG_beta4_Ex1}), and is thus considered less accurate.  

% Figure: Force-Displacement and Plot over time (Pen)

\begin{figure}[t] %[H] % "[t!]" placement specifier just for this example
\centering
\hspace{3.0cm}
\includegraphics[width=0.40\textwidth]{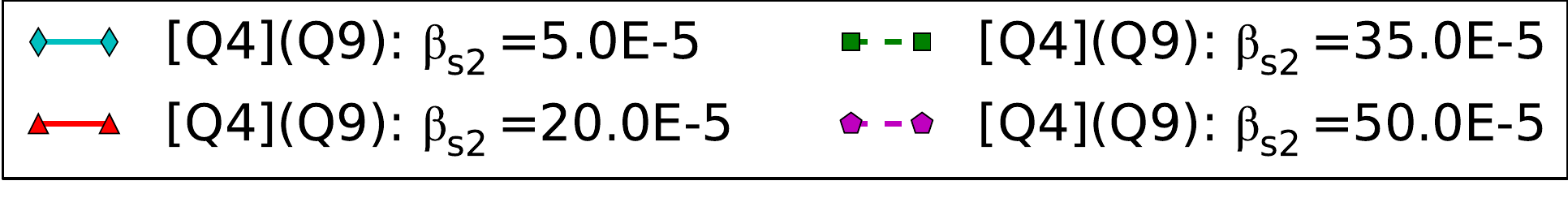}
\hspace{3.0cm}

\vspace{-0.35cm}

\centering
\subfloat[\label{fig:ForceDispPen_Ex1}  load-displacement curve]{%
\includegraphics[width=0.45\textwidth]{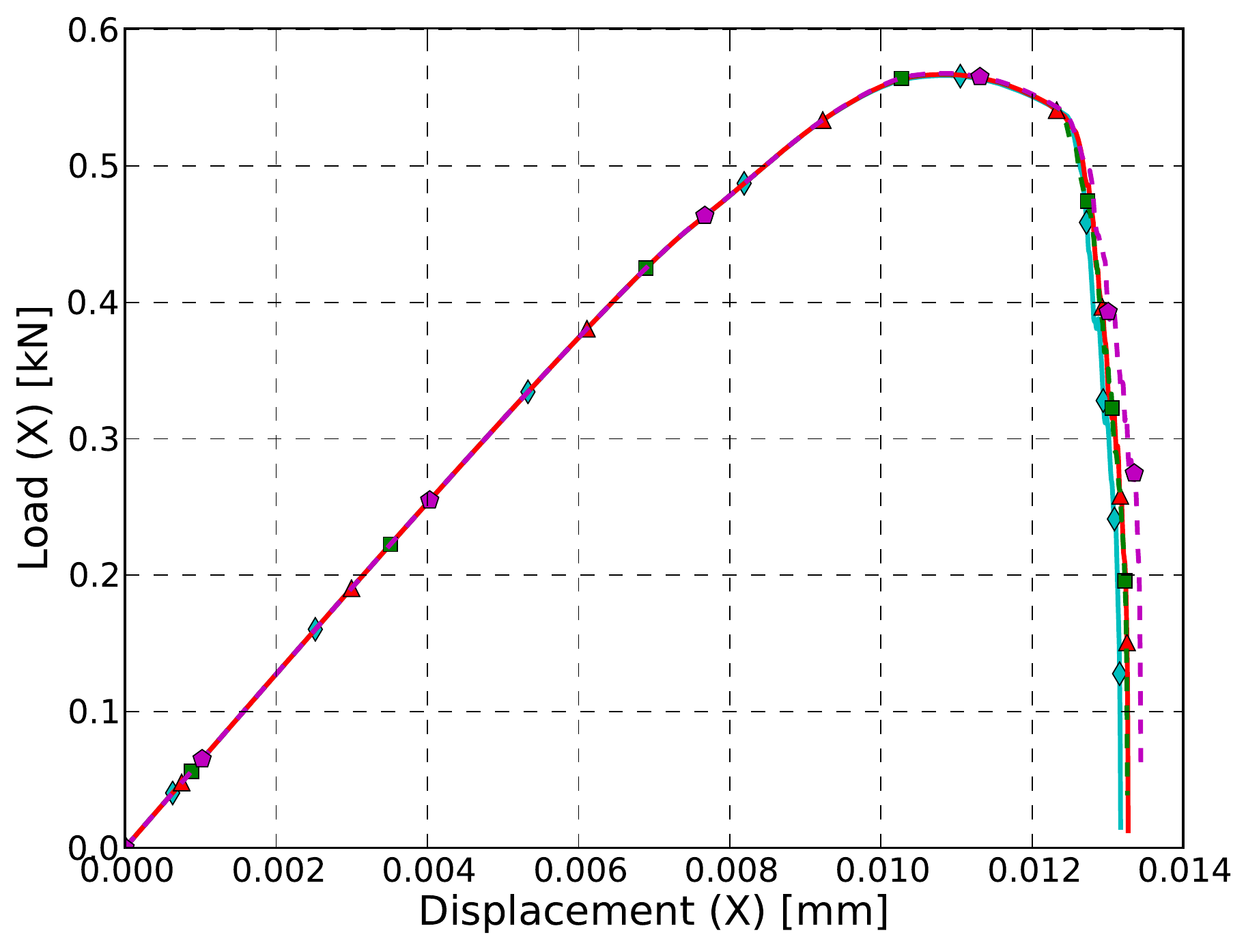}  %0.29
}
\hspace{1.0cm}
\subfloat[\label{fig:POTdPen_Ex1}  phase-field evolution over time  ]{%
\includegraphics[width=0.45\textwidth]{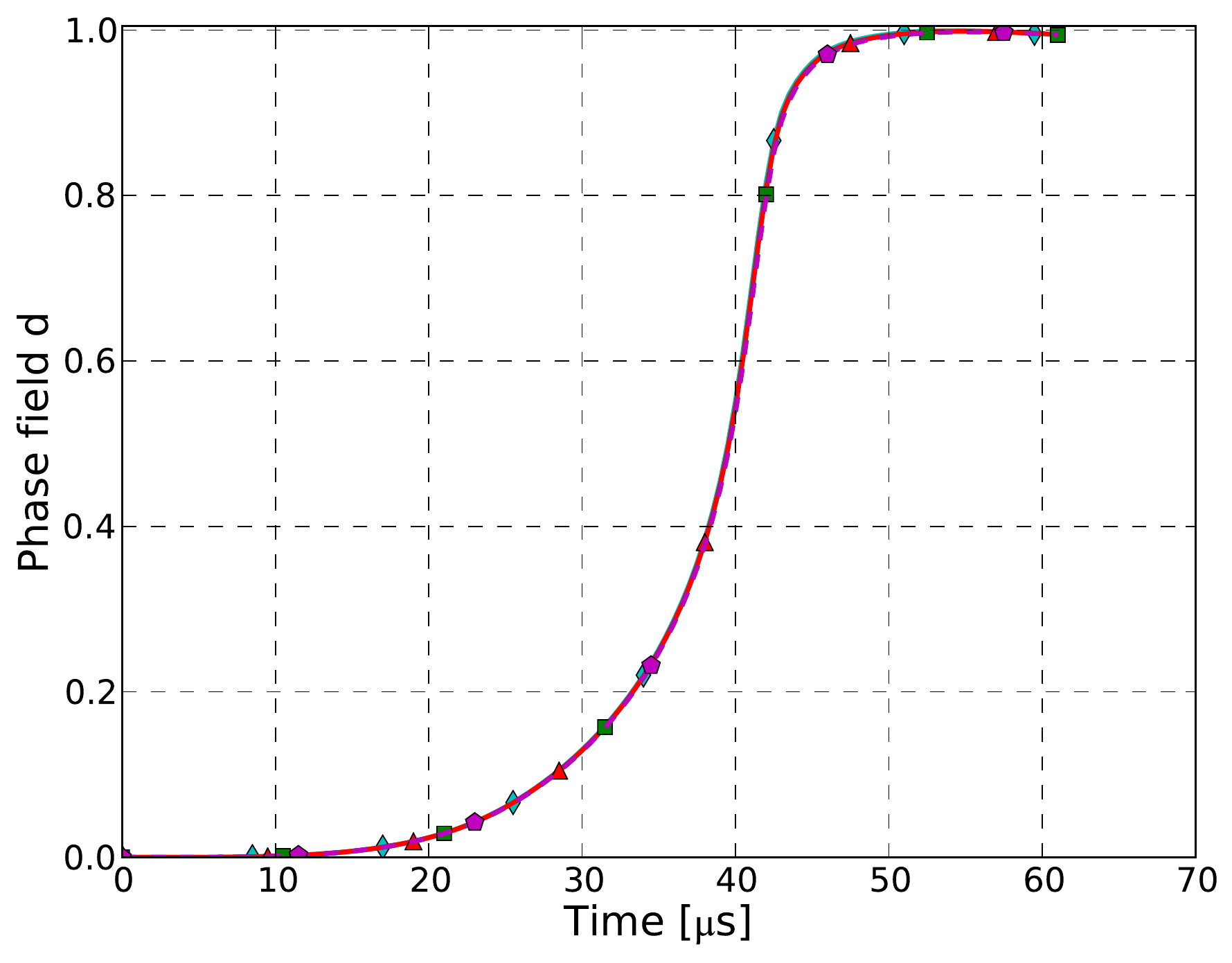}
}
\captionsetup{width=.90\textwidth}
\caption{(a) Force-displacement curve (b) phase-field solution as a function of time at crack tip (Point P) } \label{fig:ForceDispPOTdPen_Ex1}
\end{figure}

% Plot over line (Pen)

\begin{figure}[t] 
\centering
\includegraphics[width=0.45\textwidth]{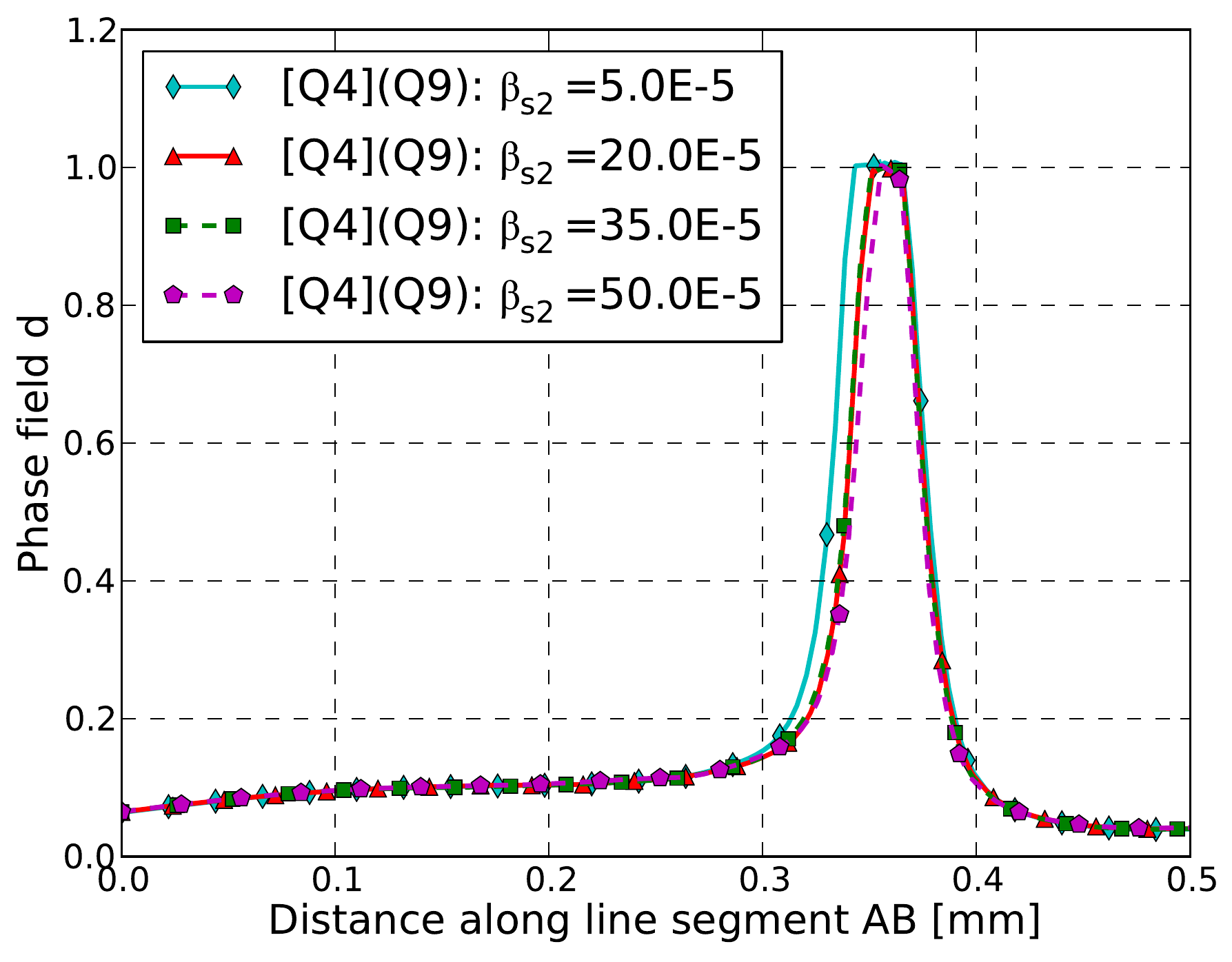}
\captionsetup{width=.90\textwidth}
\caption{ Phase-field solution along line segment AB}
\label{fig:POL_AB_Pen_Ex1}
\end{figure}

The influence of the penalty term on the mechanical response of the plate is reflected by the load--displacement curves shown for a range of penalty parameter values in \cref{fig:ForceDispPen_Ex1}. As seen, the response is linear prior to damage initiation, and is not sensitive to the penalty parameter until the abrupt loss of load carrying capacity. In this latter regime, the load-displacement curves corresponding to $\beta_{s2} = 20.0\times10^{-5}$ and $35.0\times10^{-5}$ coincide, confirming that there is a range of penalty parameter values in which insensitive solutions are obtained. \cref{fig:POTdPen_Ex1} shows the phase-field evolution over time at the tip of the notch (point $P$ in \cref{fig:Geometry_Ex1}). Clearly, crack initiation is not affected by the value of the penalty parameters. However, the topology of fully-developed cracks does change with $\beta_{s2}$. \cref{fig:POL_AB_Pen_Ex1} shows the phase-field solution at time $t=63.0\times10^{-6}$~s along line segment $AB$ (\cref{fig:Geometry_Ex1}). The phase-field solution is too diffusive with $\beta_{s2} = 5.0\times10^{-5}$, but the profiles of $\PF$ converge as the value of the penalty parameter is increased.

\begin{figure}[p] %[H] % "[t!]" placement specifier just for this example
\centering
\includegraphics[width=0.30\textwidth]{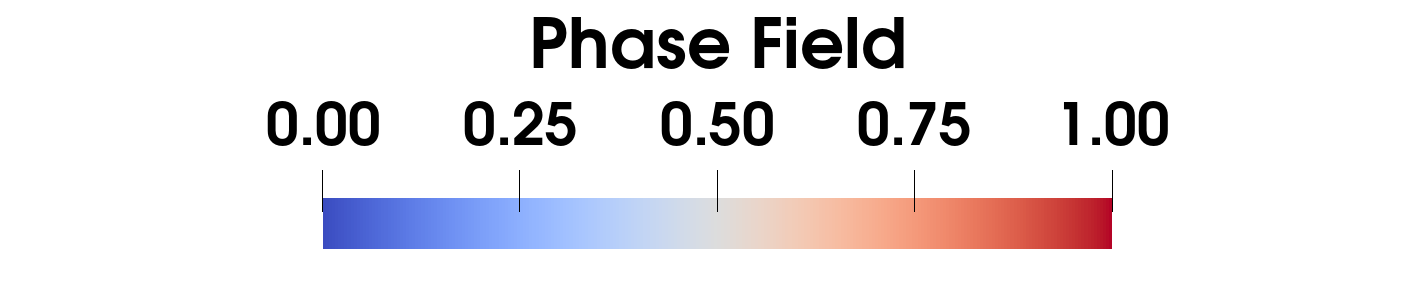}
\hspace{1.0cm}
\includegraphics[width=0.30\textwidth]{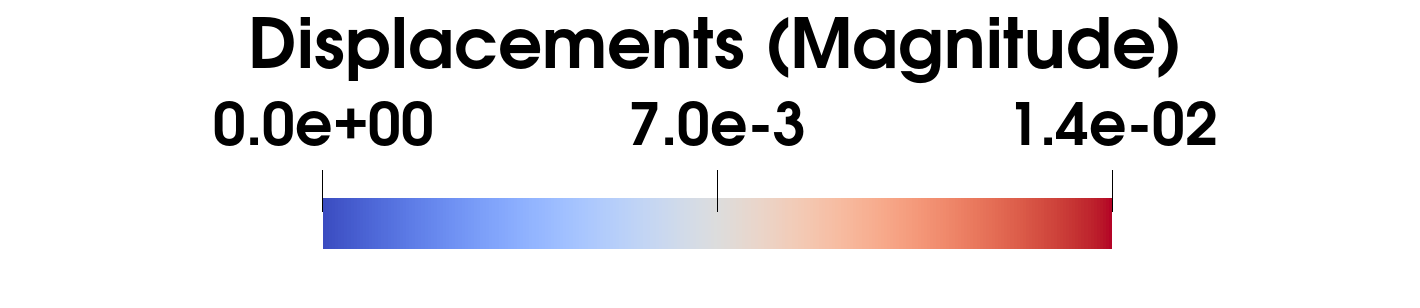}

\vspace{-0.35cm}

\centering
\subfloat[\label{fig:CDG_PF_Ex1} Scheme \text{[Q4](Q9)} - $\beta_{s2} = 20.0\times10^{-5}$]{%
\includegraphics[width=0.30\textwidth]{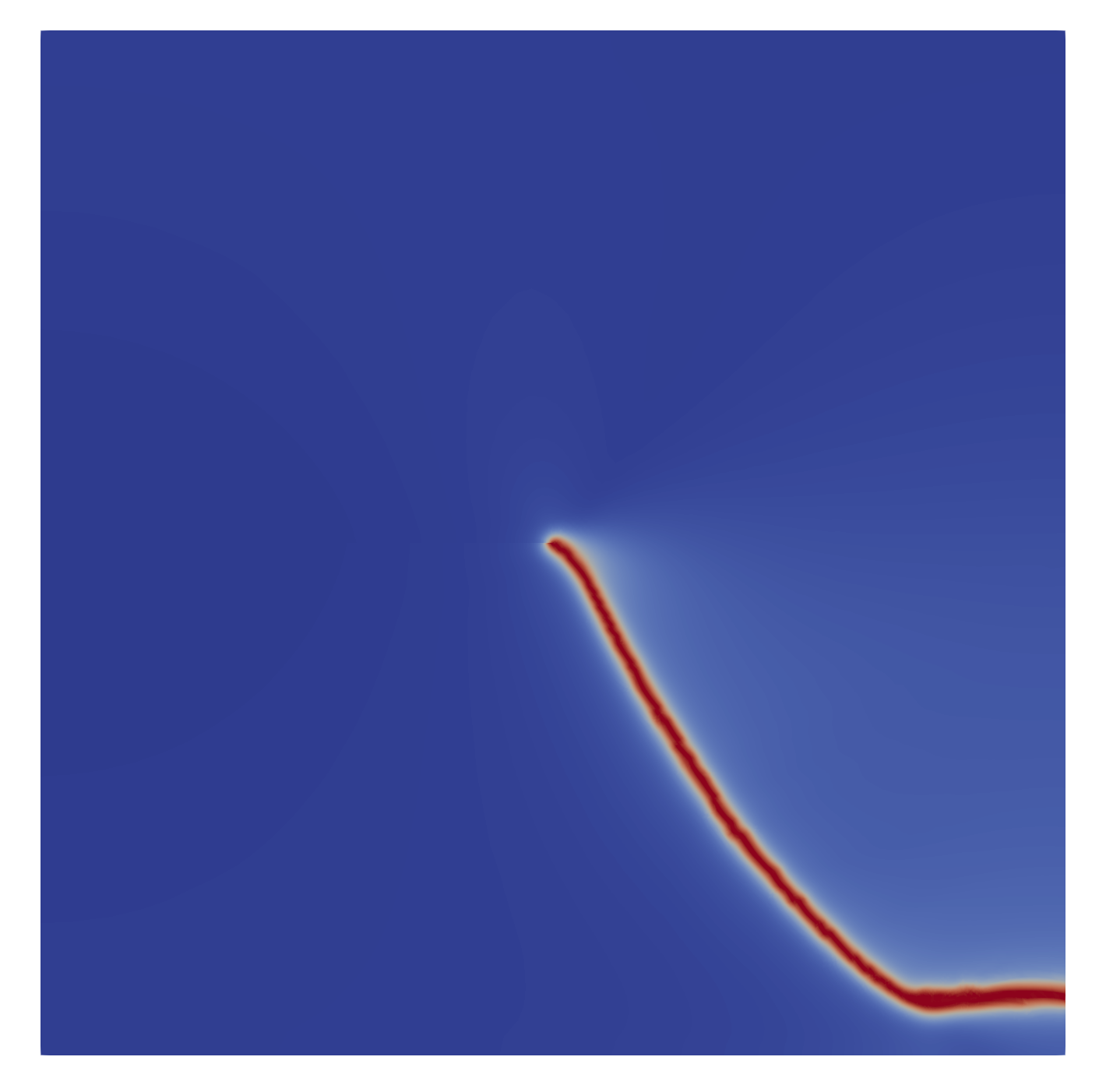}  %0.29
}
\hspace{1.0cm}
\subfloat[\label{fig:CDG_Disp_Ex1} Scheme \text{[Q4](Q9)} - $\beta_{s2} = 20.0\times10^{-5}$]{%
\includegraphics[width=0.30\textwidth]{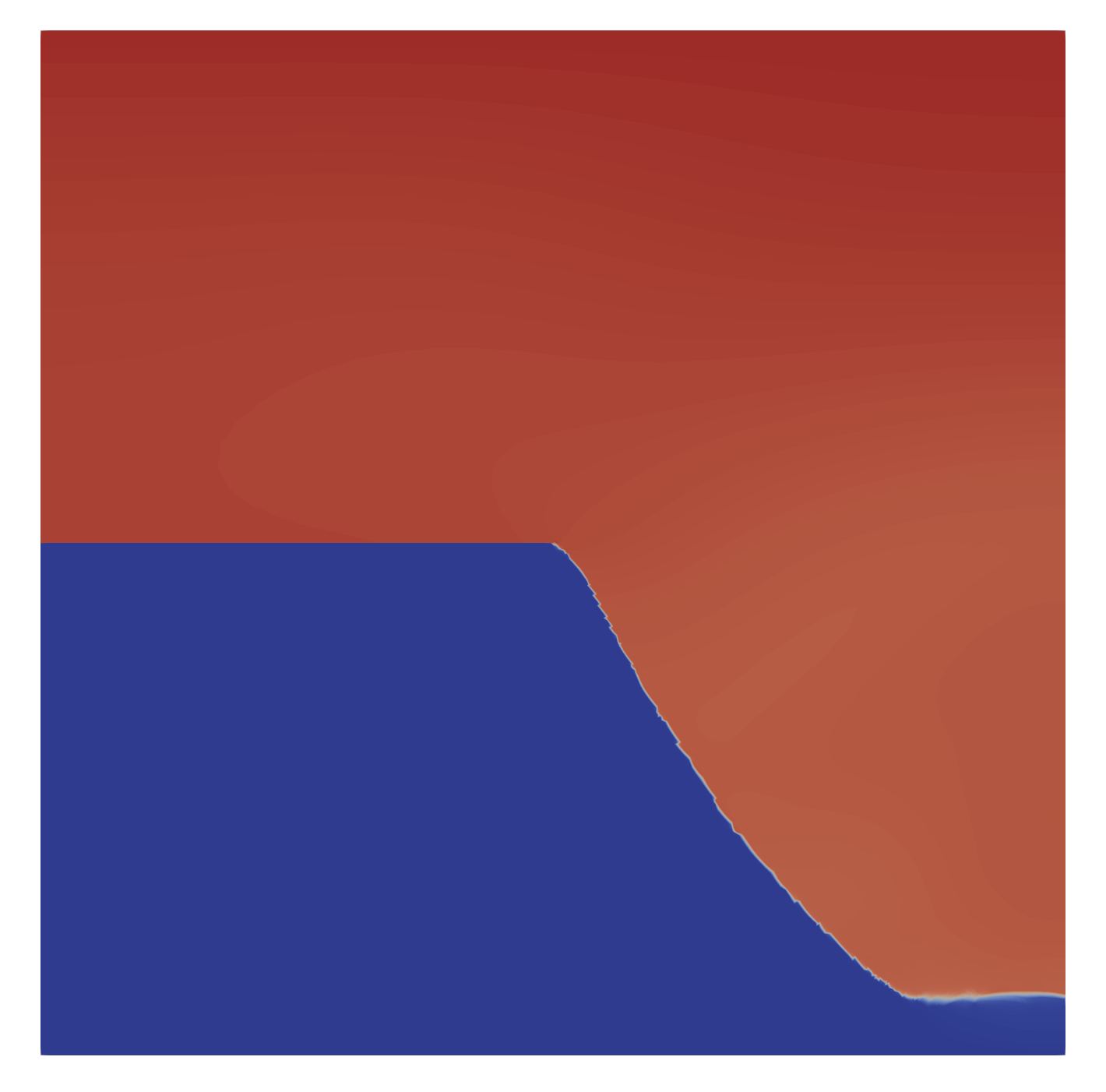}
}

%\medskip
\vspace{-0.2cm}

\centering
\subfloat[\label{fig:MixedQ9_PF_Ex1} Scheme \text{[Q4](Q9Q9)} ]{%
\includegraphics[width=0.30\textwidth]{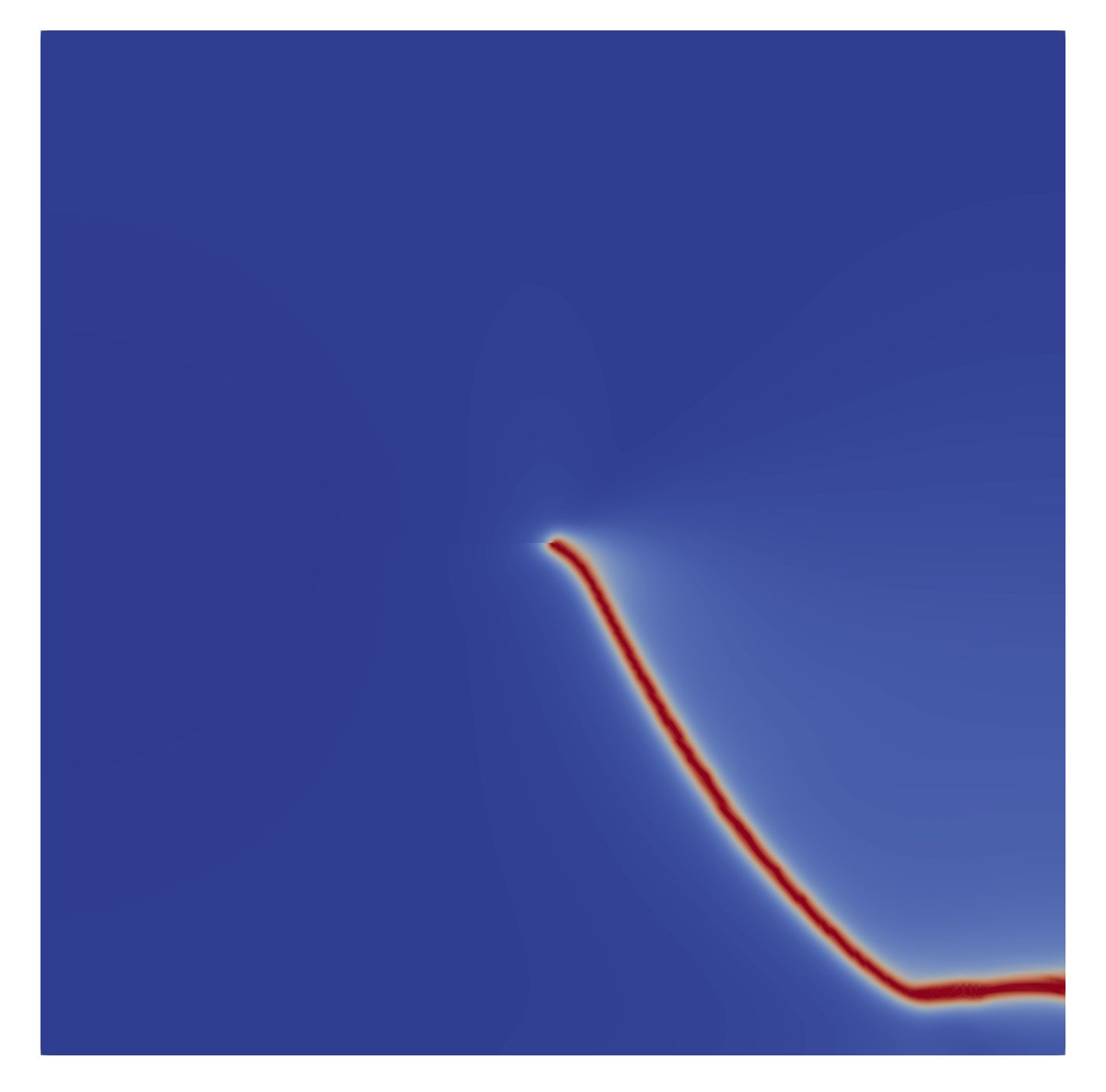}  %0.29
}
\hspace{1.0cm}
\subfloat[\label{fig:MixedQ9_Disp_Ex1} Scheme \text{[Q4](Q9Q9)} ]{%
\includegraphics[width=0.30\textwidth]{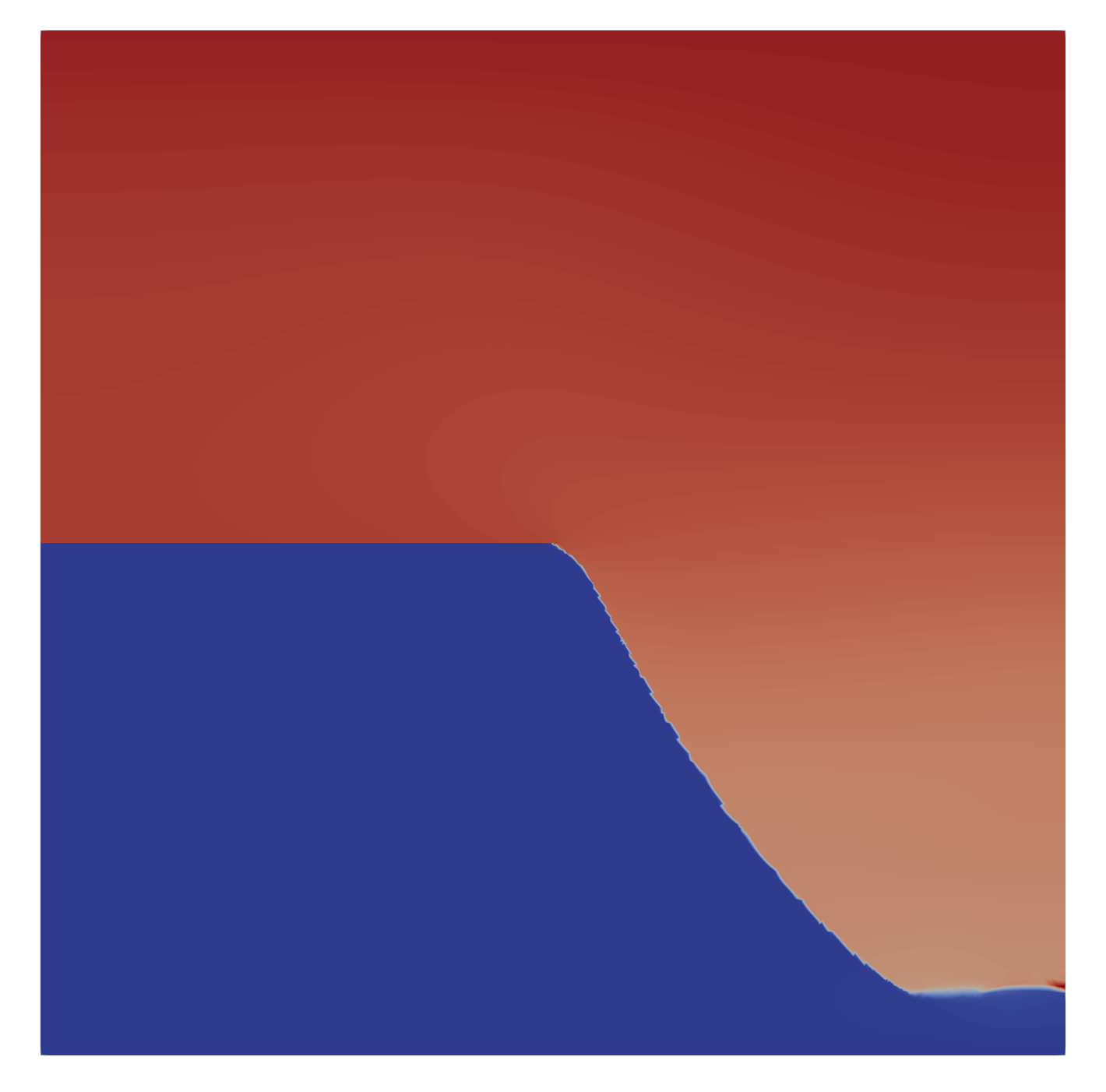}
}

%\medskip
\vspace{-0.2cm}

\centering
\subfloat[\label{fig:MixedQ4_PF_Ex1} Scheme \text{[Q4](Q4Q4)} ]{%
\includegraphics[width=0.30\textwidth]{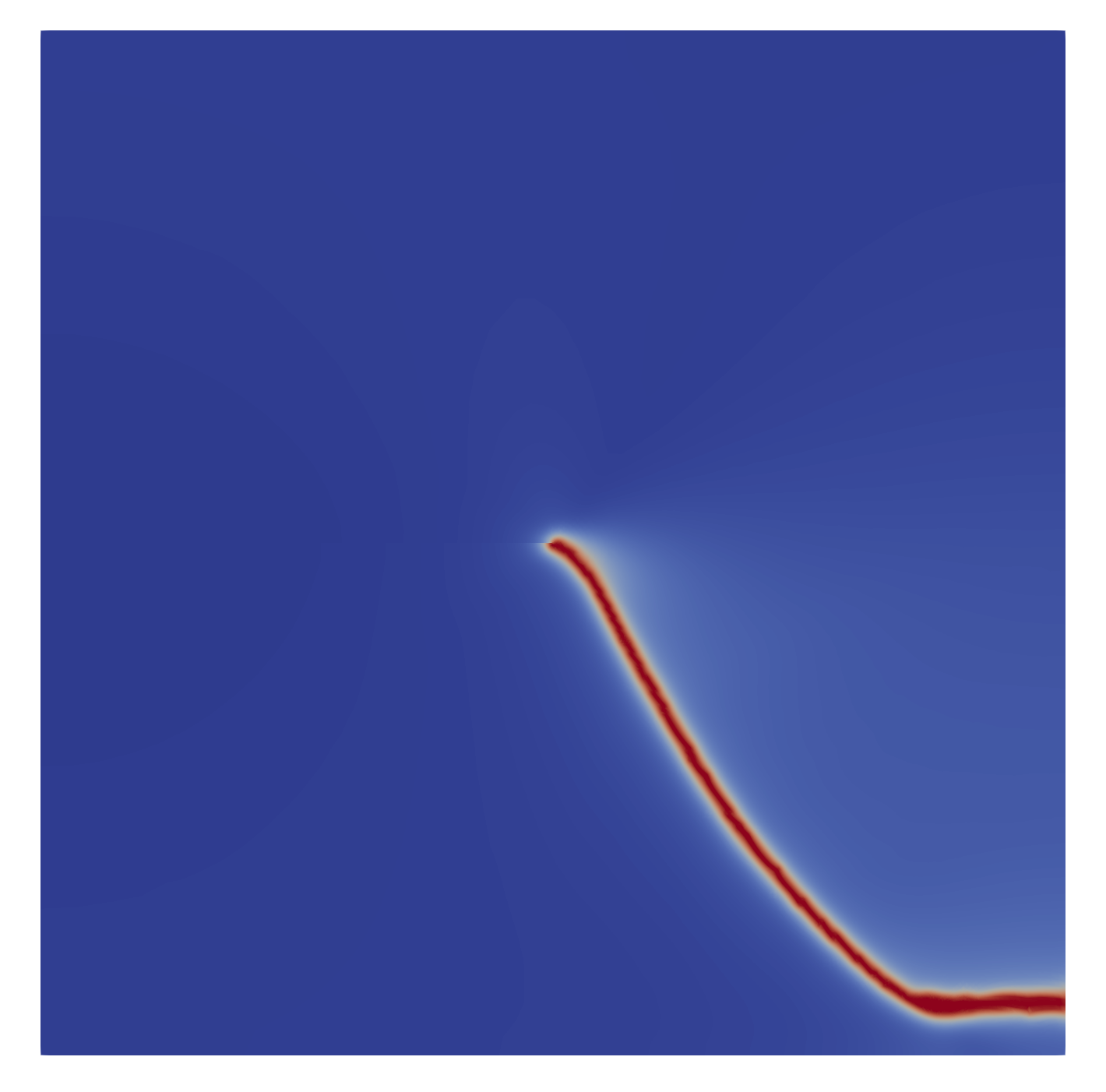}  %0.29
}
\hspace{1.0cm}
\subfloat[\label{fig:MixedQ4_Disp_Ex1} Scheme \text{[Q4](Q4Q4)} ] {%
\includegraphics[width=0.30\textwidth]{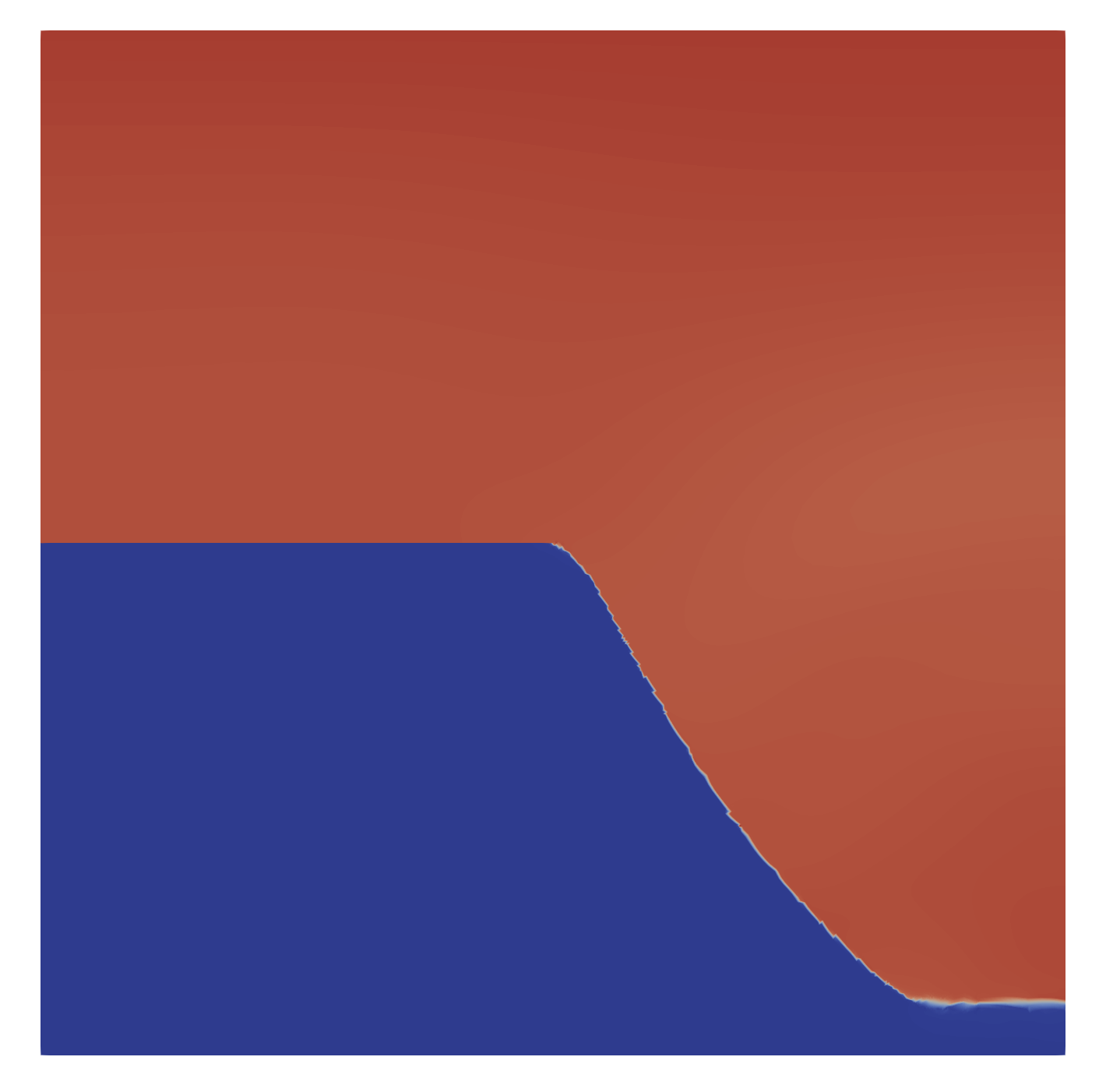}
}
\captionsetup{width=.90\textwidth}
\caption{ Comparison between the proposed Schemes (\text{[Q4](Q9)}, \text{[Q4](Q9Q9)} and \text{[Q4](Q4Q4)}) at the complete crack formation: (left) phase-field fracture paths (right) displacement contours.   } \label{fig:AllMethodsPFDisp_Ex1}
\end{figure}

This sensitivity analysis highlights the importance of selecting a penalty parameter value within the optimal range, where high accuracy is achieved without sacrificing computational efficiency or stability. Effective use of the proposed C/DG Scheme \text{[Q4](Q9)} requires this parameter calibration step, which may incur additional computational effort. In the present example problem, we set $\beta_{s2} = 20.0\times10^{-5}$ for comparison with other schemes in the next section.

\subsubsection{Comparison between numerical schemes}
\label{sss:ComparisonBetweenDifSchemes_Ex1}

%# Final time before split into two parts
%TimeRangeDict = { Pr1M4F:   [0.0,6.14e-05],
%                  Pr1M9S:   [0.0,6.461e-05],
%                  Pr1M9Pf1: [0.0,6.30e-05],
%                  Pr1M9Pf2: [0.0,6.34e-05],
%                  Pr1M9Pf3: [0.0,6.33e-05],
%                  Pr1M9Pf4: [0.0,6.40e-05] 
%                }    # 6.76\times10^{-5}

A comparison between the numerical results obtained using the three schemes under consideration is presented in \cref{fig:AllMethodsPFDisp_Ex1}. More specifically, phase-field and displacement contours, respectively, are shown at the moment where the crack reaches the fully-developed stage, corresponding to $t=63.3\times10^{-6}$~s for Scheme \text{[Q4](Q9)} (\cref{fig:CDG_PF_Ex1,fig:CDG_Disp_Ex1}), $t=64.1\times10^{-6}$~s for Scheme \text{[Q4](Q9Q9)} (\cref{fig:MixedQ9_PF_Ex1,fig:MixedQ9_Disp_Ex1}), and $t=61.4\times10^{-6}$~s for Scheme \text{[Q4](Q4Q4)} (\cref{fig:MixedQ4_PF_Ex1,fig:MixedQ4_Disp_Ex1}). Clearly, there is very good agreement between all three numerical schemes with regards to the final crack topology predicted. Moreover, the displacement discontinuity across the crack and the maximum displacement magnitude after crack formation are very similar in all cases. 

% Figure: Force-Displacement and Plot over time (Met)

\begin{figure}[t] %[H] % "[t!]" placement specifier just for this example
\centering
\hspace{3.0cm}
\includegraphics[width=0.40\textwidth]{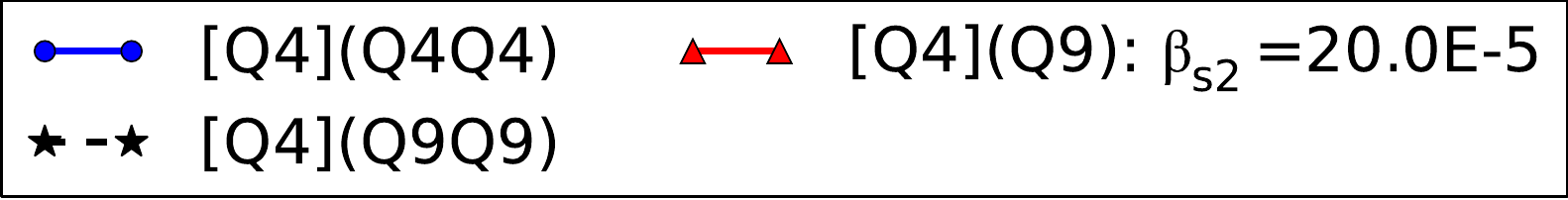}
\hspace{3.0cm}

\vspace{-0.35cm}

\centering
\subfloat[\label{fig:ForceDispMet_Ex1}  load-displacement curve]{%
\includegraphics[width=0.45\textwidth]{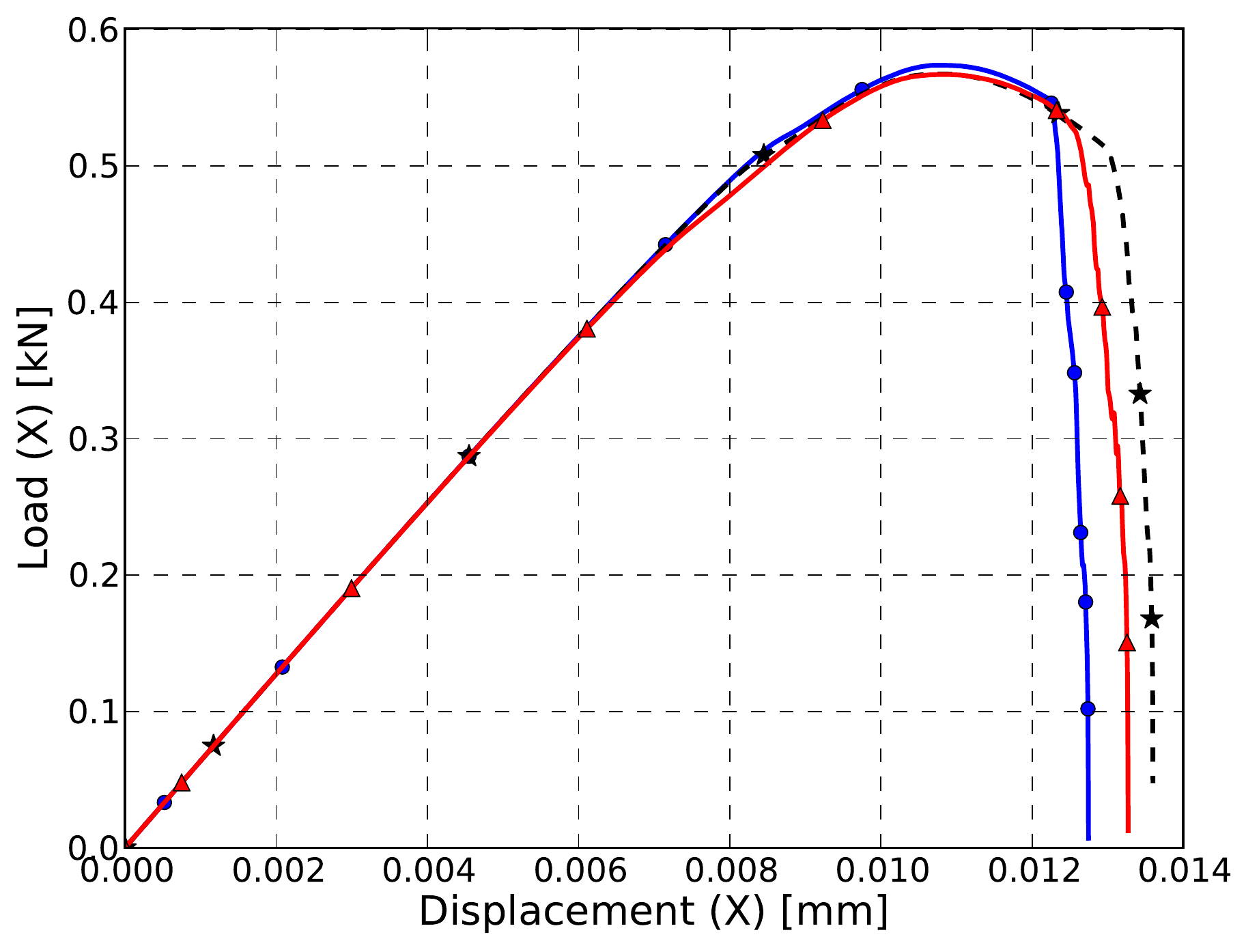}  %0.29
}
\hspace{1.0cm}
\subfloat[\label{fig:POTdMet_Ex1}  phase-field plot over time  ]{%
\includegraphics[width=0.45\textwidth]{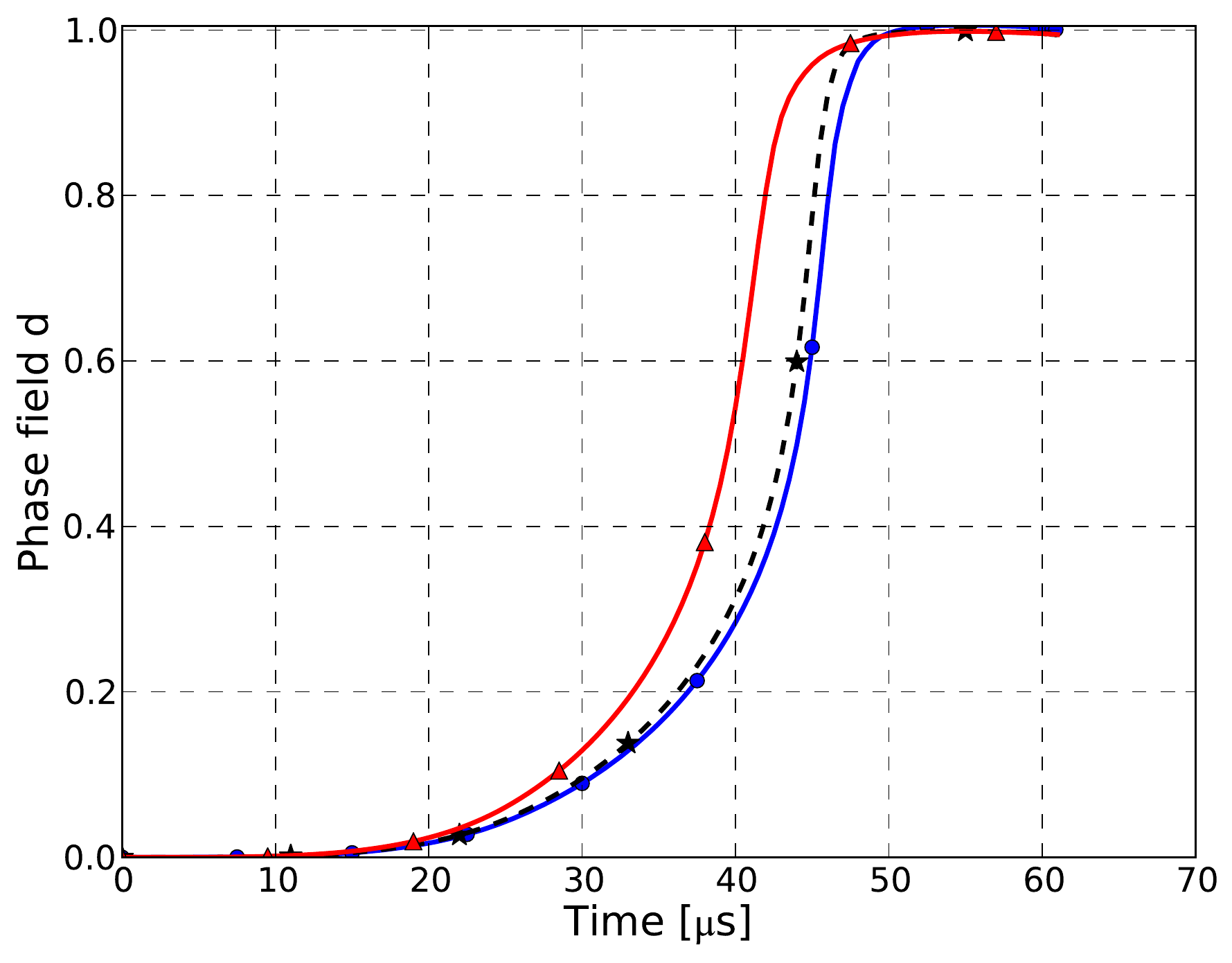}
}
\captionsetup{width=.90\textwidth}
\caption{ (a) Force-displacement curve (b) phase-field solution as a function of time at crack tip (point $P$) } \label{fig:ForceDispPOTdMet_Ex1}
\end{figure}

% Plot over line (Met)

%\begin{figure}[H] %[H] % "[t!]" placement specifier just for this example
%\centering
%\hspace{3.0cm}
%\includegraphics[width=0.40\textwidth]{Section4/Example1/POLineStudies/c_POL_ABMet.pdf}
%\hspace{3.0cm}
%
%\vspace{-0.35cm}
%
%\centering
%\includegraphics[width=0.40\textwidth]{Section4/Example1/POLineStudies/lgd_c_C3P1300vsTimeMet.pdf}
%
%\captionsetup{width=.90\textwidth}
%\caption{ } \label{fig:ForceDispPOTdMet_Ex1}
%\end{figure}

\begin{figure}[t] 
\centering
\includegraphics[width=0.45\textwidth]{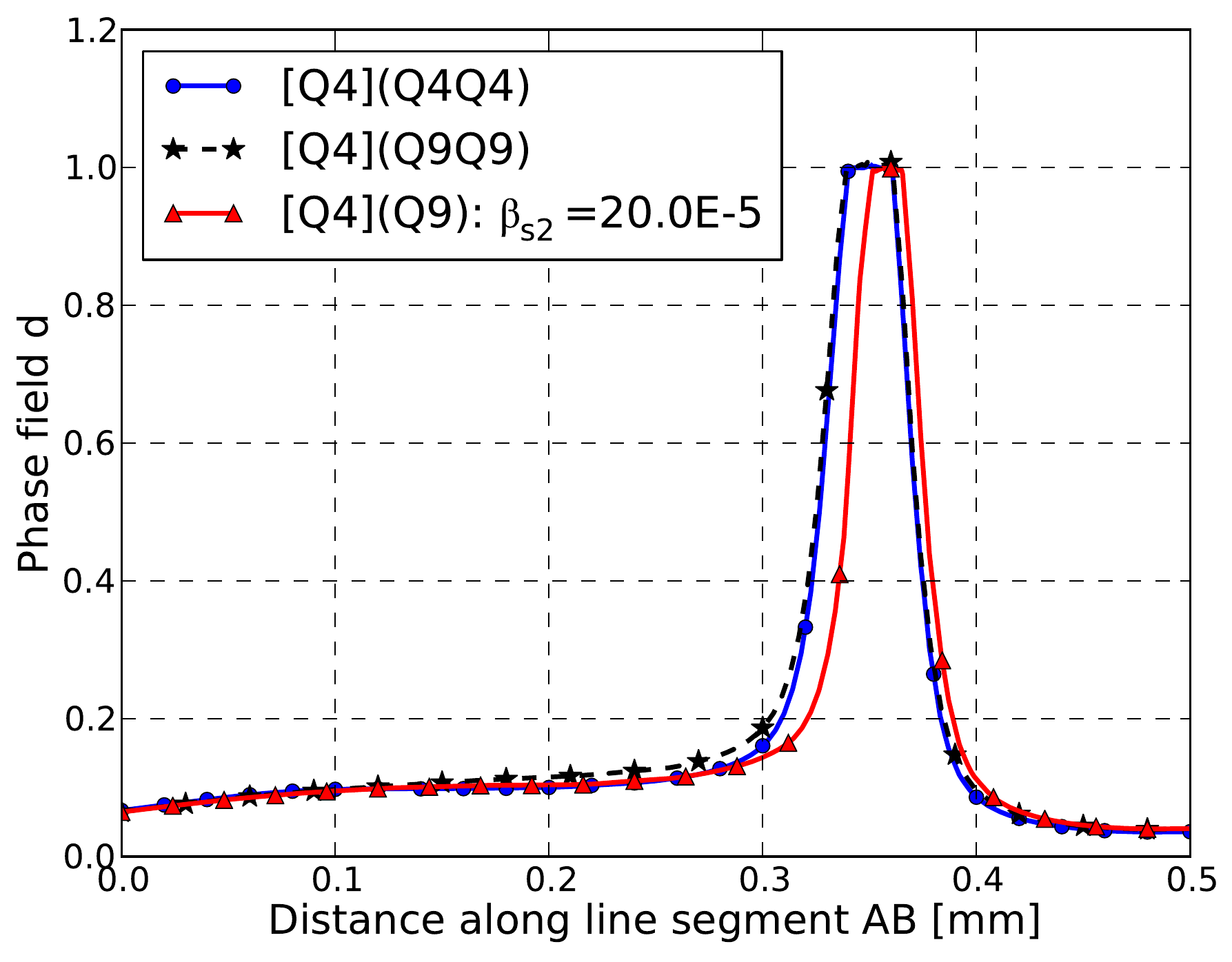}
\captionsetup{width=.90\textwidth}
\caption{ Phase-field solution along line segment $AB$ }
\label{fig:POL_AB_Met_Ex1}
\end{figure}

Load--displacement curves are shown in \cref{fig:ForceDispMet_Ex1} to assess the differences between the numerical schemes in (the more quantitative) terms of mechanical response. As expected, the initial linear behavior is identical for all numerical schemes, since damage is minimal and the same (conventional finite element) method is used to solve the momentum equation. However, the load capacity (or critical load) is slightly higher with Scheme \text{[Q4](Q4Q4)} compared with the two other schemes. Differences are also observed prior to critical loads. Scheme \text{[Q4](Q9)} deviates first from the other two schemes, and this deviation is attributed to the earlier damage initiation and accompanying material softening predicted by this scheme. Finally, it is noted that after the critical load is reached, the load--displacement curve computed using the C/DG Scheme \text{[Q4](Q9)} lies between the two curves obtained using the mixed FEM schemes. 

With the two second-order schemes, \text{[Q4](Q9)} and \text{[Q4](Q9Q9)}, complete loss of load-carrying capacity takes place later compared to Scheme \text{[Q4](Q4Q4)}. This is due to their different convergence behavior when used in conjunction with an adaptive time-stepping staggered approach. Given that the spatial discretization controls the convergence behavior within a time step, and noting that the adaptive staggered approach is sensitive to the solution at each time step, this sensitivity increases the dissimilarities among the solutions produced by the different numerical schemes.

\cref{fig:POTdMet_Ex1} shows the phase-field evolution over time at the notch tip (point $P$ in \cref{fig:Geometry_Ex1}). The C/DG Scheme \text{[Q4](Q9)} predicts the fastest crack initiation, leading to a different softening behavior even before reaching the critical load, as shown in \cref{fig:ForceDispMet_Ex1}. The phase-field solution along line segment $AB$ (\cref{fig:Geometry_Ex1}), at the moment where the crack reaches the
fully-developed stage, is shown in \cref{fig:POL_AB_Met_Ex1}. This figure demonstrates very good agreement, in terms of the predicted crack path, between Schemes \text{[Q4](Q4Q4)} and \text{[Q4](Q9Q9)}, which differ slightly from Scheme \text{[Q4](Q9)}.

Finally, we note that the present example problem is a variant of those presented in \cite{miehe2010phase,borden2012phase}. In these two works, the problem was treated using a second-order phase-field theory coupled to small-strain elasticity theory with a strain-based spectral decomposition of the elastic energy under quasi-static conditions. Although the geometry and boundary conditions are identical, our different modeling approach is reflected in our numerical results. More specifically, we observe a slightly different exit point of the crack path at the bottom right corner, which can be attributed mainly to our adoption of finite deformation kinematics and the volumetric/deviatoric energy decomposition described in \cref{ssec:ElasticConstitutiveLaw}, especially in such a shear-dominated problem.

%\subsubsection{Scalability analysis}

%%\input{Tables/Section4/Example1/Table_Parallel.tex}

%[Total time, Efficiency or SpeedUp]

% Second Example

\subsection{Double cantilever beam experiment}

\subsubsection{Problem description and discretization}

\begin{figure}[t] 
\centering
\subfloat[\label{fig:Geometry_Ex2}]{%
\includegraphics[width=0.55\textwidth]{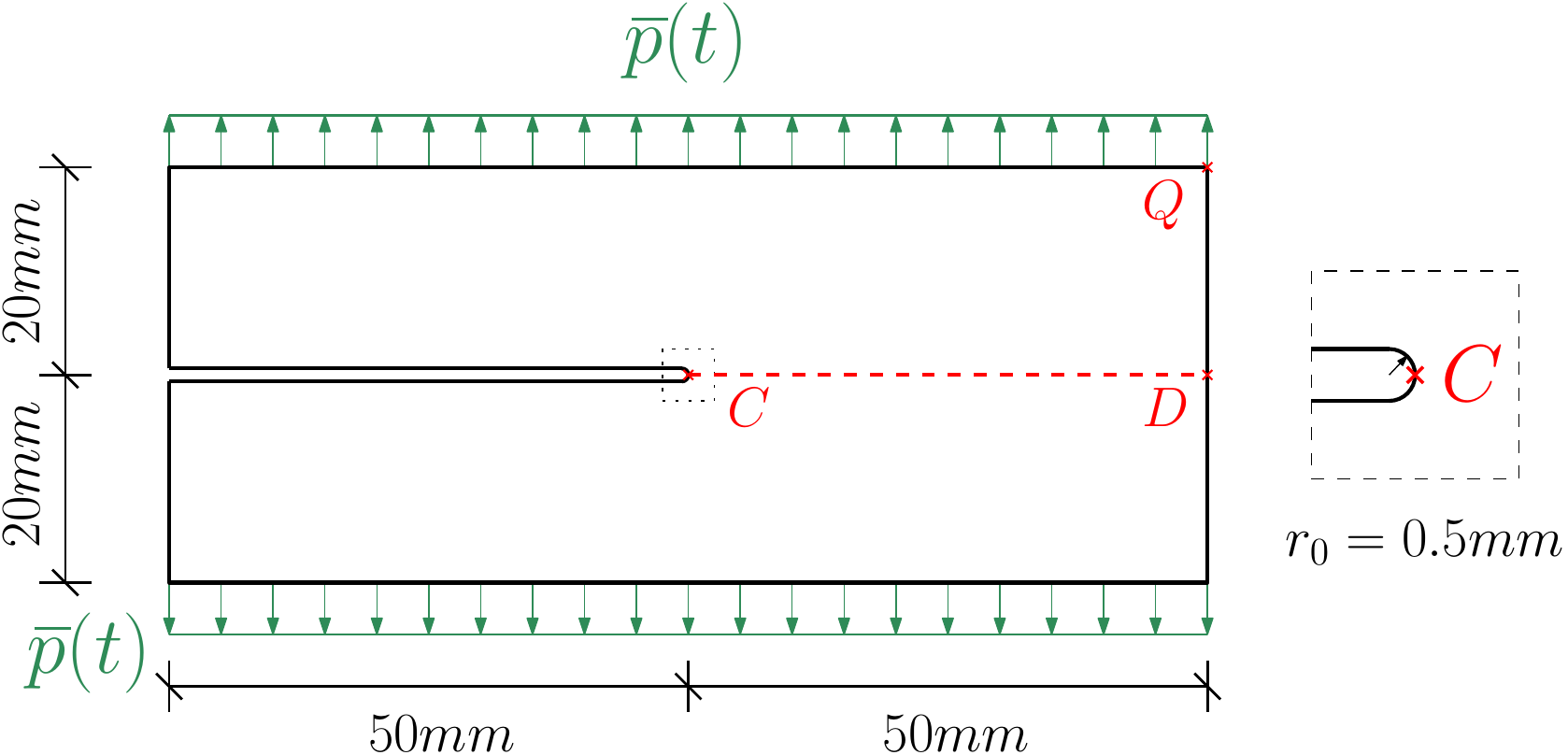}  %0.29
}
\hspace{1.0cm}
\subfloat[\label{fig:Load_Ex2}]{%
\includegraphics[width=0.35\textwidth]{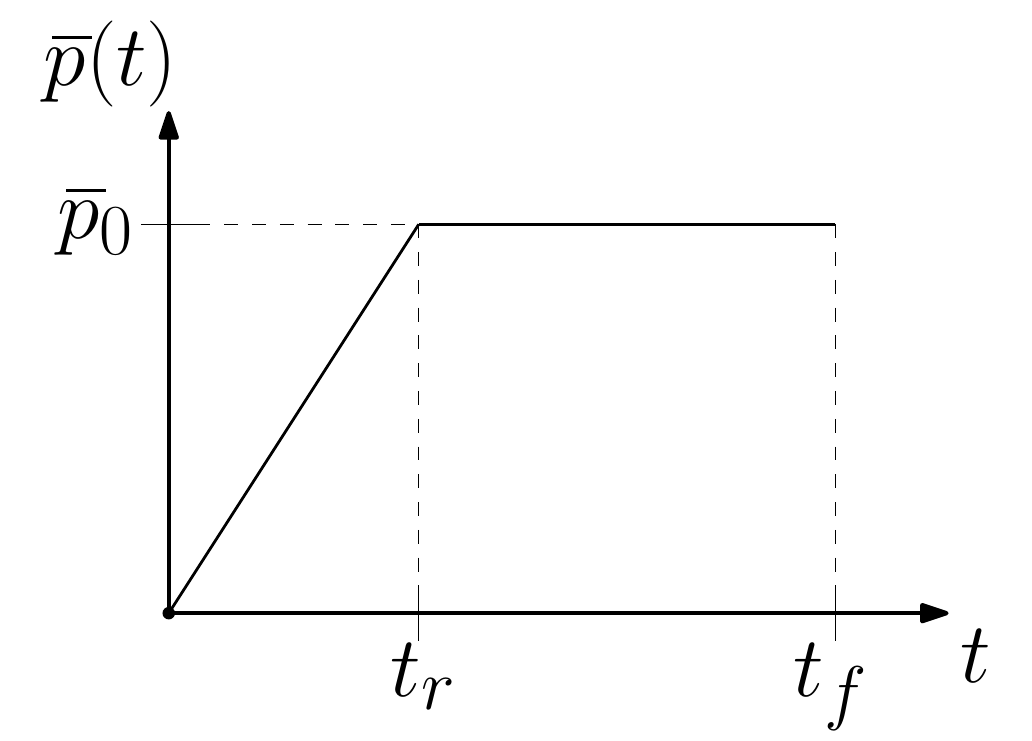}
}
\captionsetup{width=.90\textwidth}
\caption{ (a) Geometry, boundary conditions and notch detail of a rectangular plate loaded in tension, (b) The time-dependent traction load, imposed on the top and bottom surfaces.}
\label{fig:GeometryLoad_Ex2}
\end{figure}

In this numerical example, a uniformly distributed tensile load $\thickbar{p}$ is dynamically applied to a notched rectangular plate with dimensions $40 \times 100$~mm. The geometry and boundary conditions are presented in \cref{fig:Geometry_Ex2}. As shown in this figure, the plate has a pre-existing notch with a rounded tip of radius $r_0=0.5$~mm. This leads to stress concentration at the tip of the notch (point $C$), causing crack initiation to take place at that point. The traction imposed on the top and bottom surfaces, varies with time in the manner shown in \cref{fig:Load_Ex2}, with ramp time $t_r=1.0\times10^{-6}$~s, final time $t_f=95.0\times10^{-6}$~s, and maximum pressure $\thickbar{p}_0=1$~MPa. The remaining boundaries are traction-free (with no constraints on displacements). Under this dynamic loading, the crack propagates with increasing speed and branches into two or more cracks as observed in experiments \cite{ramulu1985mechanics}. Since capturing such crack branching accurately is a non-trivial challenge (e.g., see \cite{borden2012phase,song2008comparative}), this is considered a standard benchmark problem in dynamic fracture modeling. 

% This problem is commonly studied in the literature because its loading conditions lead to dynamic crack branching as shown in \cite{borden2012phase,song2008comparative}. 

%Units: mm, s, N, MPa, Mg, K, mJ
\begin{table}[b]
\centering
\caption{Material properties and simulation parameters}
\label{tbl:simulation_parameters_Ex2}
\resizebox{0.55\textwidth}{!}{%
\begin{tabular}{|lccc|}
\hline
Property & Notation  &  Value  &  Unit  \Tstrut\Bstrut\\
\hline
Young's Modulus & $E$ & 32.0E+3 & $MPa$    \Tstrut\\
Poisson's Ratio & $\nu$ & 0.2 &  -  \Tstrut\\
Mass Density & $\rho$ & 2.45E-9 & $Mg$ ${(mm)}^{-3}$  \Tstrut\\
Critical Energy Release Rate & $G_c$ & 3.0E-3 & $mJ$ ${{(mm)}^{-2}}$ \Tstrut\\
Length Scale & $\lzr$ & 0.125 & $mm$ \Tstrut\\
Residual Stiffness Parameter & $\eta_0$ & 1.0E-6 & -  \Tstrut\Bstrut\\
%Newmark beta & $\beta$ & 0.3025 & - \Tstrut\\
%Newmark gamma & $\gamma$ & 0.6 & - \Tstrut\Bstrut\\
\hline
\end{tabular}
}
\end{table}

\begin{figure}[t] 
\centering
\subfloat[\label{fig:MeshQ4_Ex2} (Q4) mesh]{%
\includegraphics[width=0.45\textwidth]{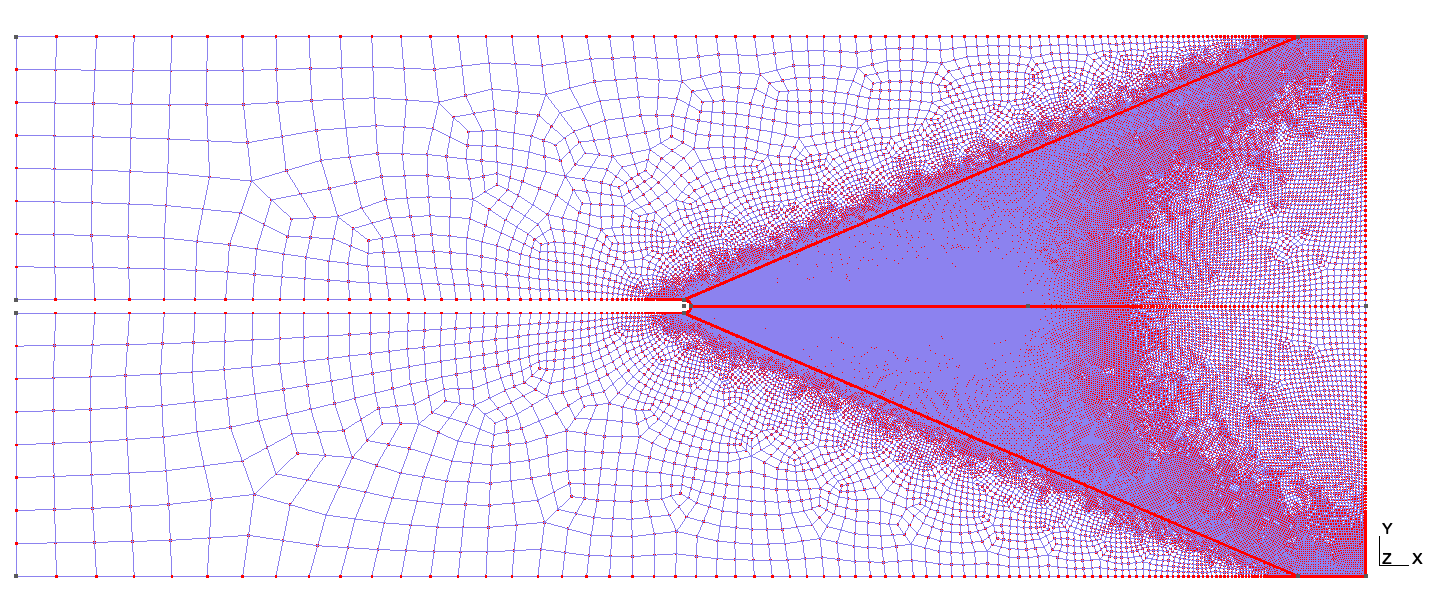}  %0.29
}
\hspace{1.0cm}
\subfloat[\label{fig:MeshQ9_Ex2} (Q9) mesh]{%
\includegraphics[width=0.45\textwidth]{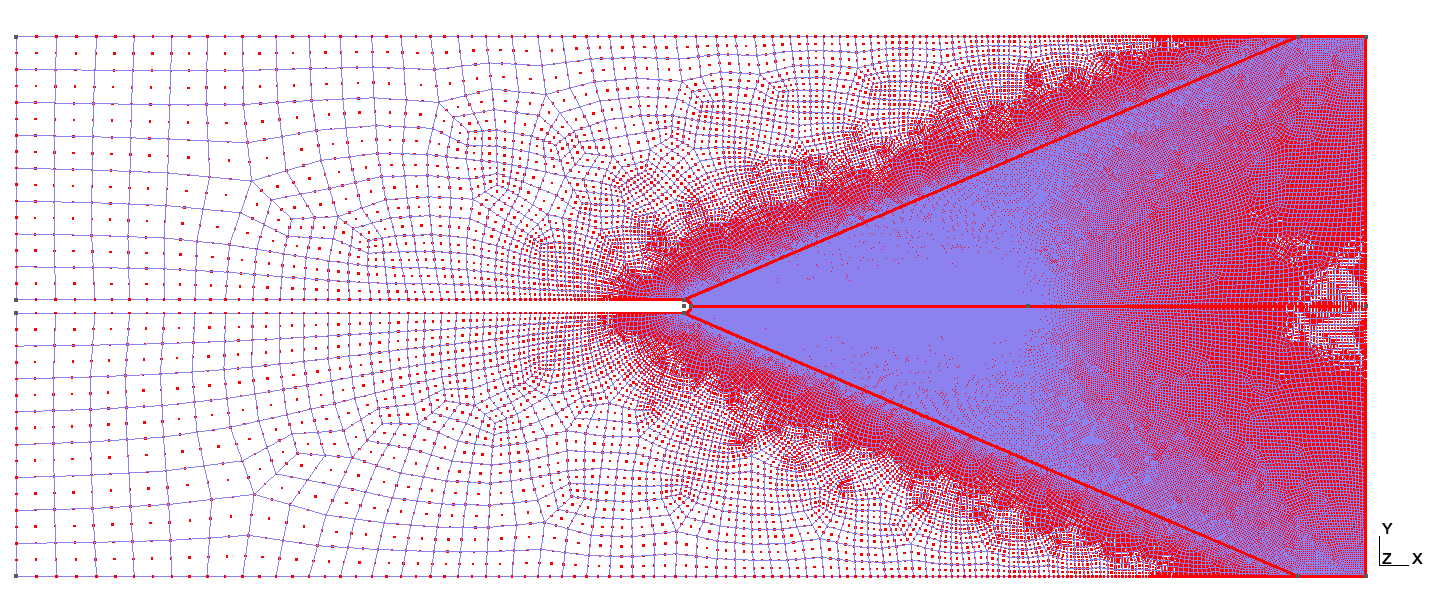}
}
\captionsetup{width=.90\textwidth}
\caption{Unstructured meshes used in this problem, refined over a triangular area behind the notch-tip where the crack branching is expected to develop: (a) four-node quadrilateral elements, and (b) nine-node quadrilateral elements.  }
\label{fig:MeshesQ4andQ9_Ex2}
\end{figure}

In this problem, two unstructured meshes are constructed to discretize the domain: (i) the mesh shown in \cref{fig:MeshQ4_Ex2}, consisting of 99,456 Q4 elements and 101,364 nodes and (ii) the mesh shown in \cref{fig:MeshQ9_Ex2}, comprising 99,456 Q9 elements and 401,636 nodes. Both meshes are considerably refined in a triangular area behind the notch-tip where crack branching is expected to develop. As in the previous problem, solution of the linear momentum balance equation involves two unknown displacements per node, while the phase-field equation involves either one unknown per node in the case of the C/DG method or two unknowns per node in the case of mixed FEM. The total number of DoFs are 597,971 for the C/DG Scheme [Q4](Q9), 996,412 for Scheme [Q4](Q9Q9), and 399,060 for Scheme [Q4](Q4Q4). As noted in the previous section, this illustrates an advantage of the C/DG scheme over the bi-quadratic mixed FEM scheme [Q4](Q9Q9) with regards to computational efficiency, despite being closely comparable to it in terms of accuracy.

% reported in \cref{tbl:DoFsPerSchemes_Ex2}. 

%\begin{table}[H]
%\centering
%\caption{The total number of Degrees of Freedom for each scheme}
%\label{tbl:DoFsPerSchemes_Ex2}
%%\resizebox{0.3\textwidth}{!}{%
%\begin{tabular}{|l|r|}
%\hline
%Numerical Schemes  &  DoFs  \Tstrut\Bstrut\\
%\hline
%Scheme  [Q4](Q9) & 597,971 \Tstrut\\
%Scheme [Q4](Q9Q9) & 996,412   \Tstrut\\
%Scheme [Q4](Q4Q4)  & 399,060   \Tstrut\Bstrut\\  
%\hline
%\end{tabular}
%%}
%\end{table}

The material parameters used are given in \cref{tbl:simulation_parameters_Ex2}. The mass density and and the elastic and fracture parameters are chosen in accordance with other works in the literature focused on reproducing dynamic crack branching experiments. We also use the same Newmark parameters as in the previous numerical example (see \cref{tbl:simulation_parameters_Ex1}). 

\subsubsection{Sensitivity analysis: Influence of the penalty parameter}

Since the penalty parameter is generally dependent on the mesh (length-scale $h_e$), element type, and material parameters (e.g., see \cite{mourad2007bubble,versino_globallocal_2014,versino_globallocal_2015}), and as a result, cannot be estimated \textit{a priori}, we conduct a sensitivity analysis to study the influence of this penalty parameter on the phase-field solution obtained for this problem using Scheme \text{[Q4](Q9)}.

\begin{figure}[t] %[H] % "[t!]" placement specifier just for this example
\centering
\hspace{3.0cm}
\includegraphics[width=0.35\textwidth]{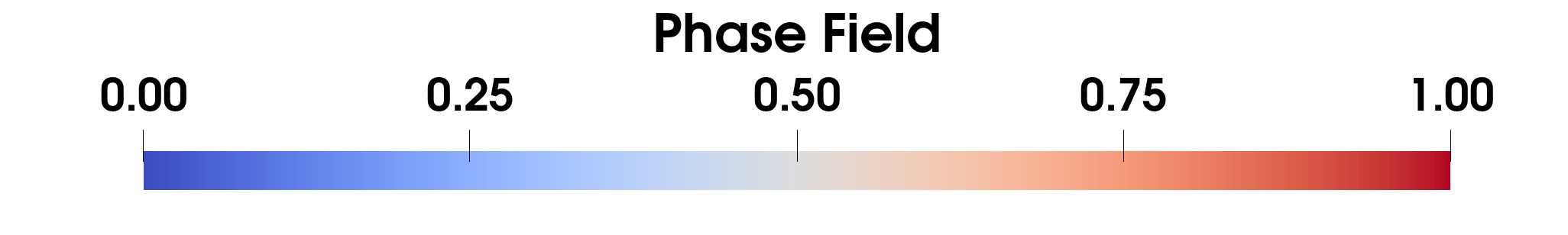}
\hspace{3.0cm}

\vspace{-0.35cm}

\centering
\subfloat[\label{fig:CDG_beta0_Ex2}  $\beta_{s2} = 5.0\times10^{-4}$]{%
\includegraphics[width=0.45\textwidth]{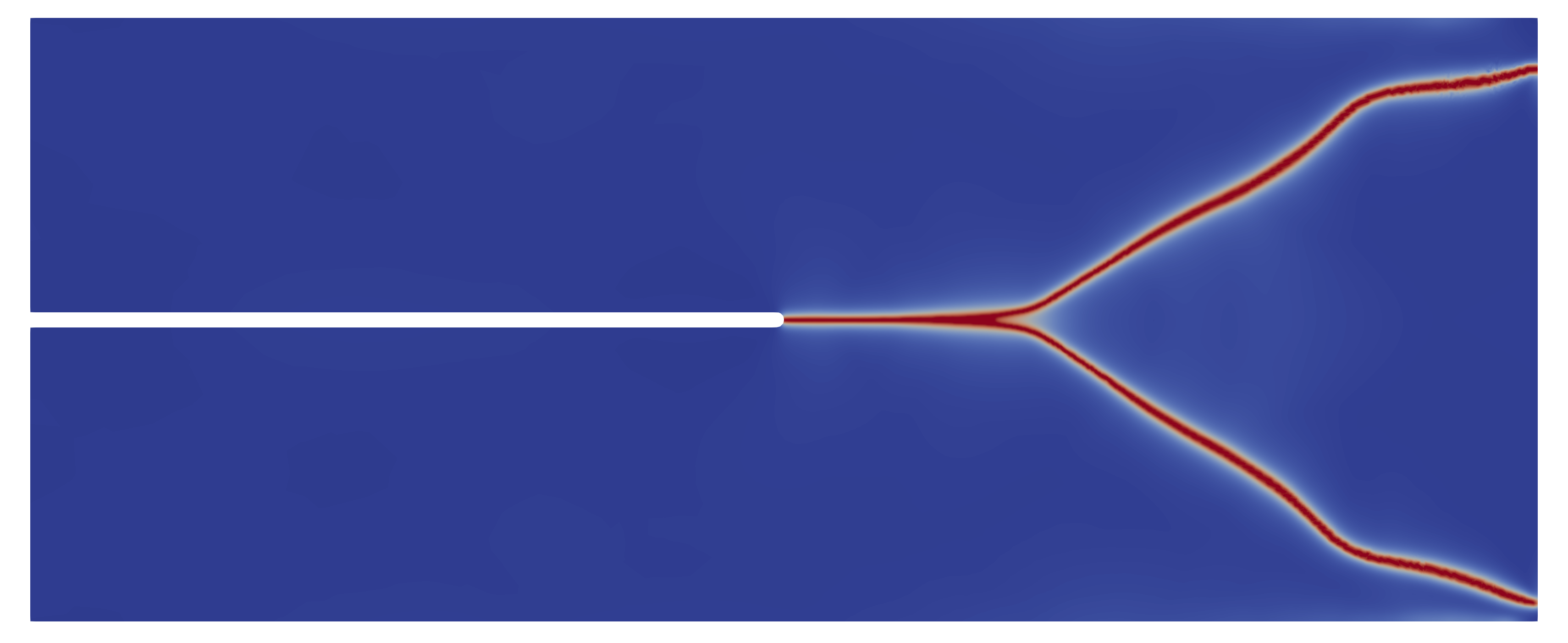}  %0.29
}
\hspace{1.0cm}
\subfloat[\label{fig:CDG_beta1_Ex2}  $\beta_{s2} = 2.0\times10^{-2}$  ]{%
\includegraphics[width=0.45\textwidth]{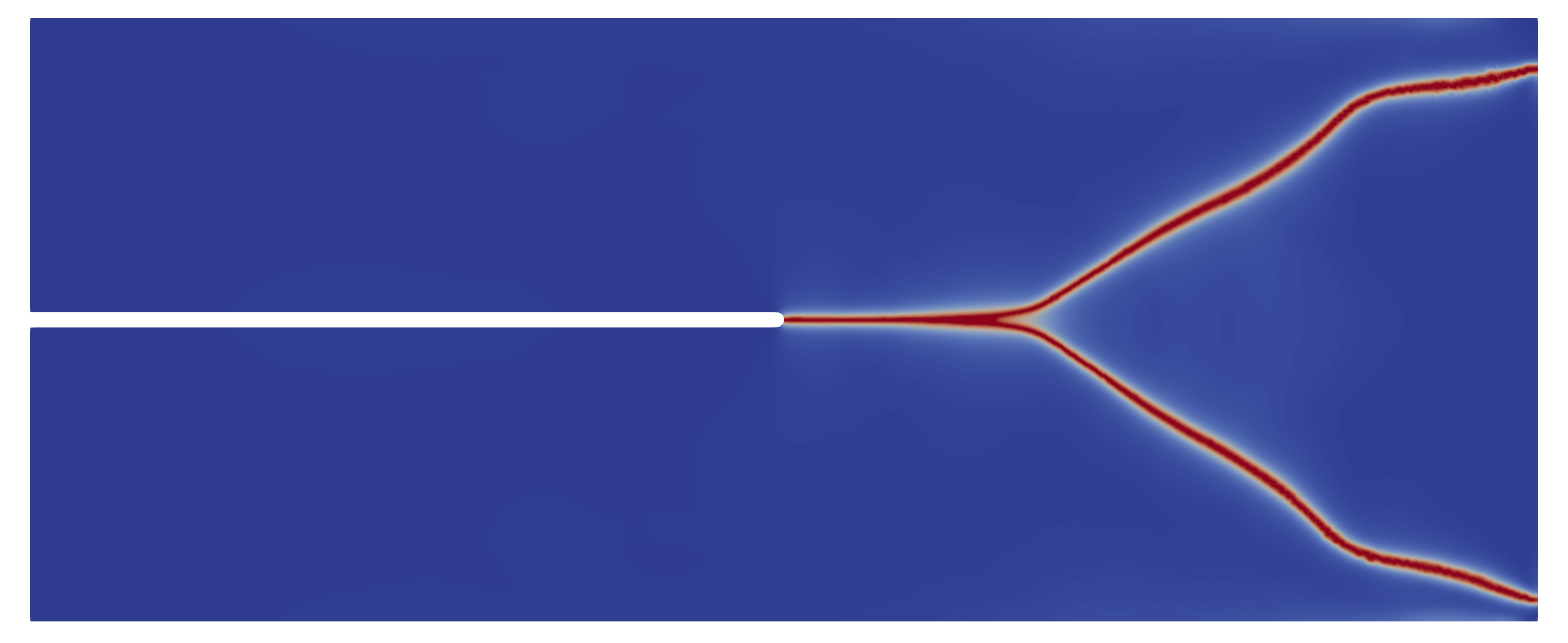}
}

%\medskip
\vspace{-0.2cm}

\centering
\subfloat[\label{fig:CDG_beta2_Ex2}  $\beta_{s2} = 3.5\times10^{-2}$]{%
\includegraphics[width=0.45\textwidth]{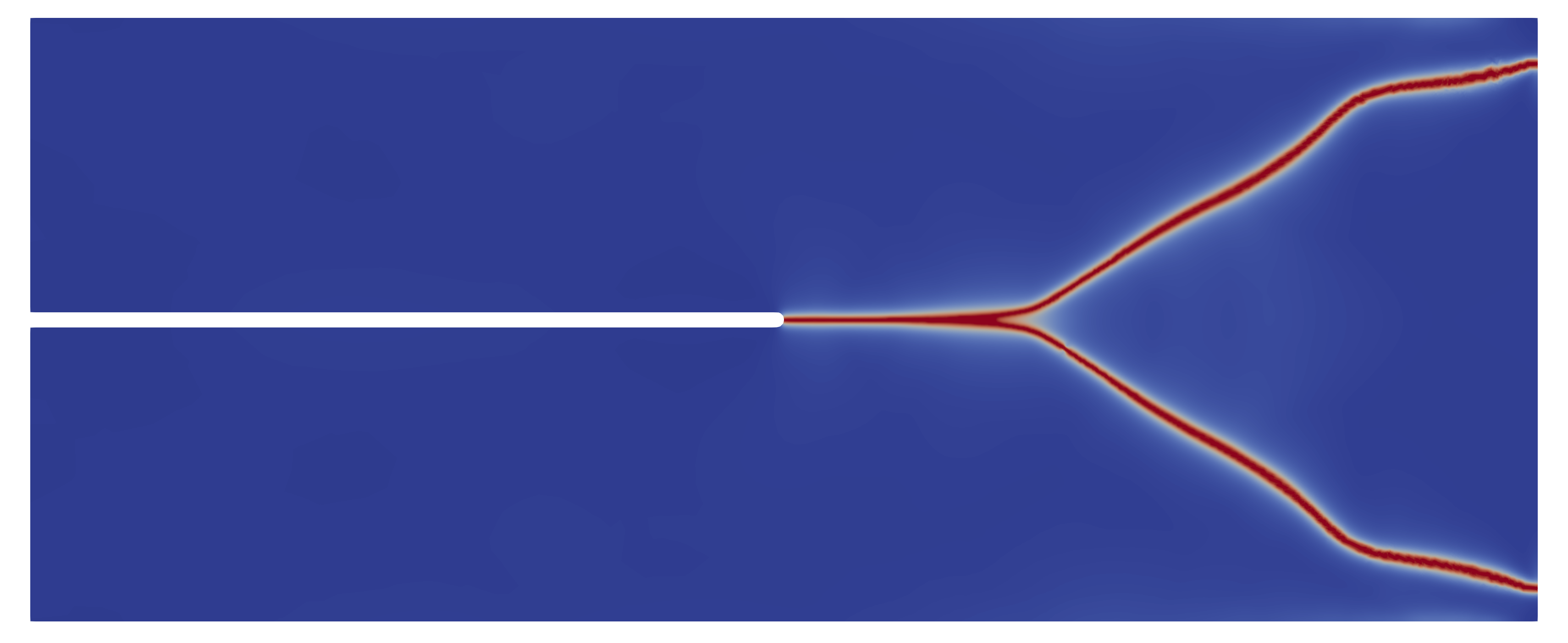}  %0.29
}
\hspace{1.0cm}
\subfloat[\label{fig:CDG_beta3_Ex2}  $\beta_{s2} = 5.0\times10^{-2}$  ]{%
\includegraphics[width=0.45\textwidth]{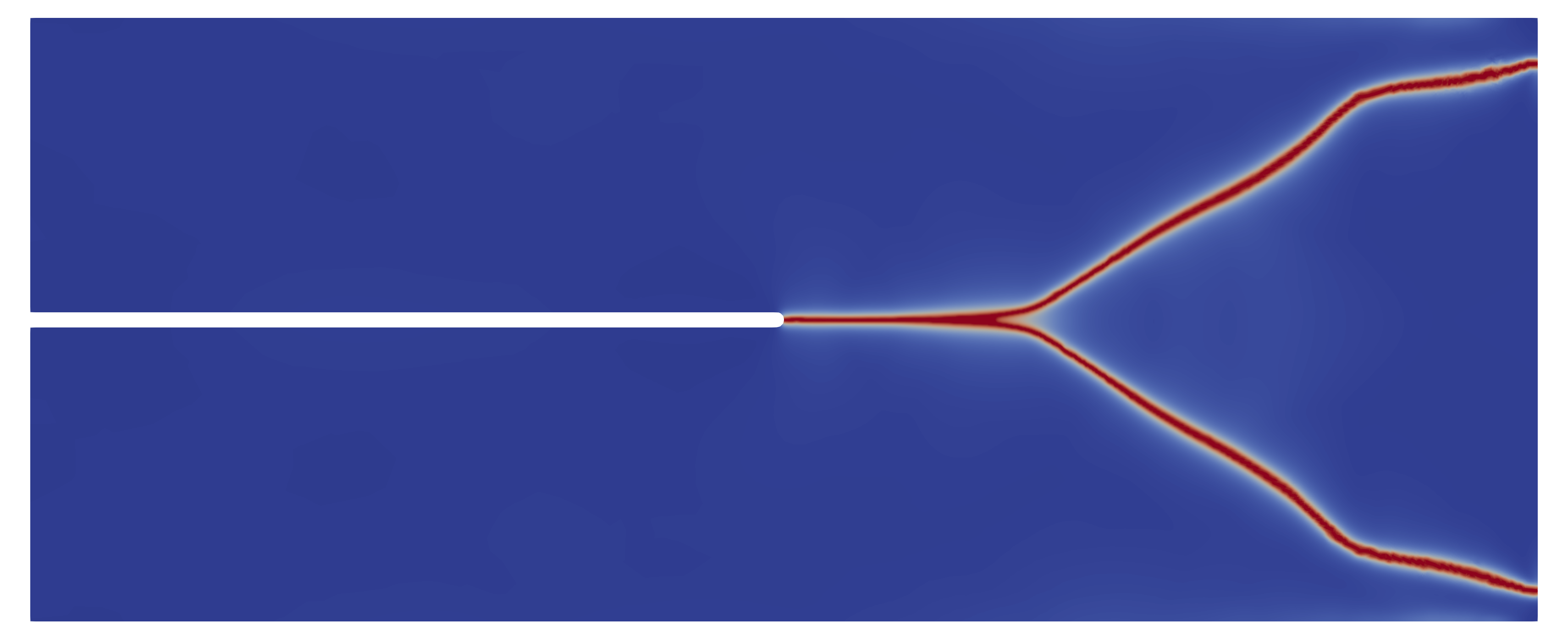}
}

\captionsetup{width=.90\textwidth}
\caption{ A penalty-parameter study: phase-field fracture paths at $t=95.0\times10^{-6}$~s are depicted when different penalty parameters in Scheme \text{[Q4](Q9)} are used. } \label{fig:AllMethodsPFbetas_Ex2}
\end{figure}

The phase-field solutions, computed using the C/DG scheme with different penalty parameter values, and representing the predicted crack path at  $t=t_f=95.0\times10^{-6}$~s, are shown in \cref{fig:AllMethodsPFbetas_Ex2}. With a very small penalty parameter value, $\beta_{s2} = 5.0\times10^{-4}$, crack branching is captured as can be seen in \cref{fig:CDG_beta0_Ex2}, but the two crack branches are not symmetric with respect to the horizontal center-line of the plate. This lack of symmetry is also observed with $\beta_{s2} = 2.0\times10^{-2}$ (\cref{fig:CDG_beta1_Ex2}), but is almost imperceptible when regularity requirements are enforced more rigorously by increasing the penalty parameter to $\beta_{s2} = 3.5\times10^{-2}$ and then $5.0\times10^{-2}$ (\cref{fig:CDG_beta2_Ex2,fig:CDG_beta3_Ex2}, respectively). Note that increasing the penalty parameter further gives rise to convergence difficulties due to ill-conditioning. We also note that the phase-field solution here does not exhibit the time lag observed in the first example problem (see \cref{sss:SensitivityAnalysisEx1}).

% Figure: Displacement and Phase-field Plot over time (Pen)

\begin{figure}[t]  %[H] % "[t!]" placement specifier just for this example
\centering
\hspace{3.0cm}
\includegraphics[width=0.40\textwidth]{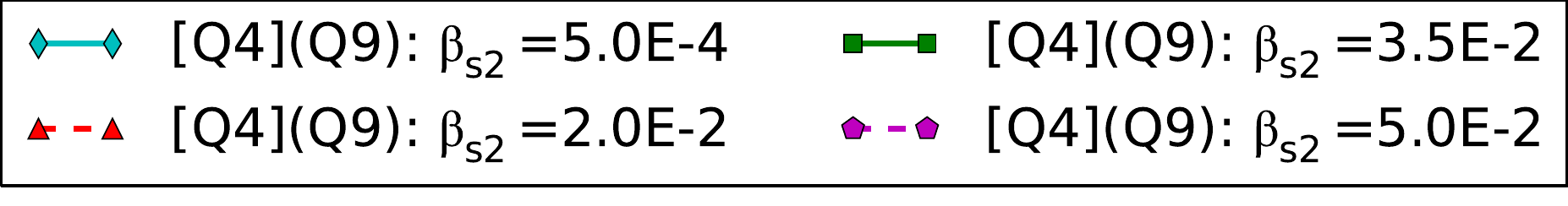}
\hspace{3.0cm}

\vspace{-0.35cm}

\centering
\subfloat[\label{fig:POTdispPen_Ex2}  displacement magnitude plot over time]{%
\includegraphics[width=0.45\textwidth]{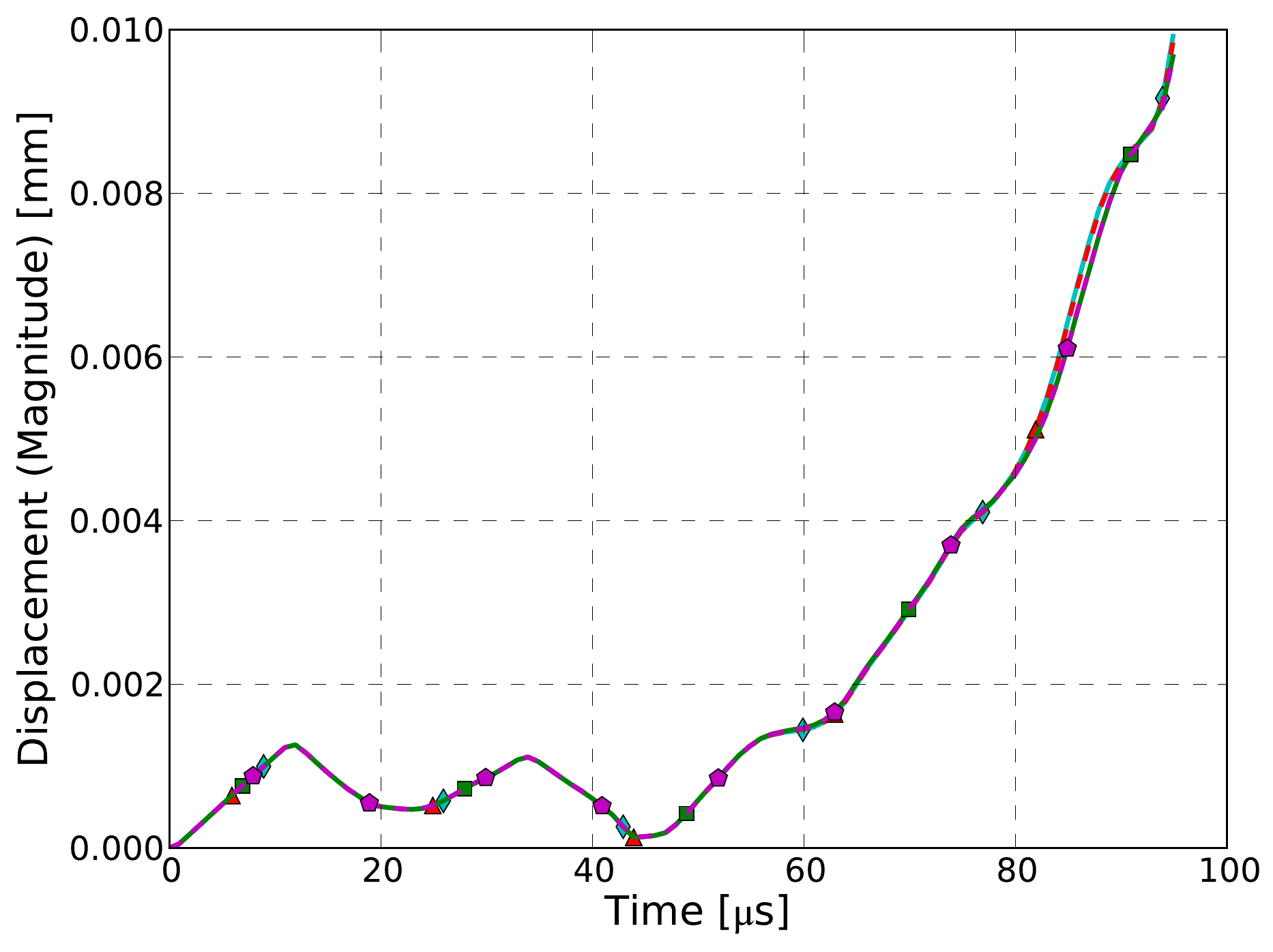}  %0.29
}
\hspace{1.0cm}
\subfloat[\label{fig:POTdPen_Ex2}  phase-field plot over time  ]{%
\includegraphics[width=0.45\textwidth]{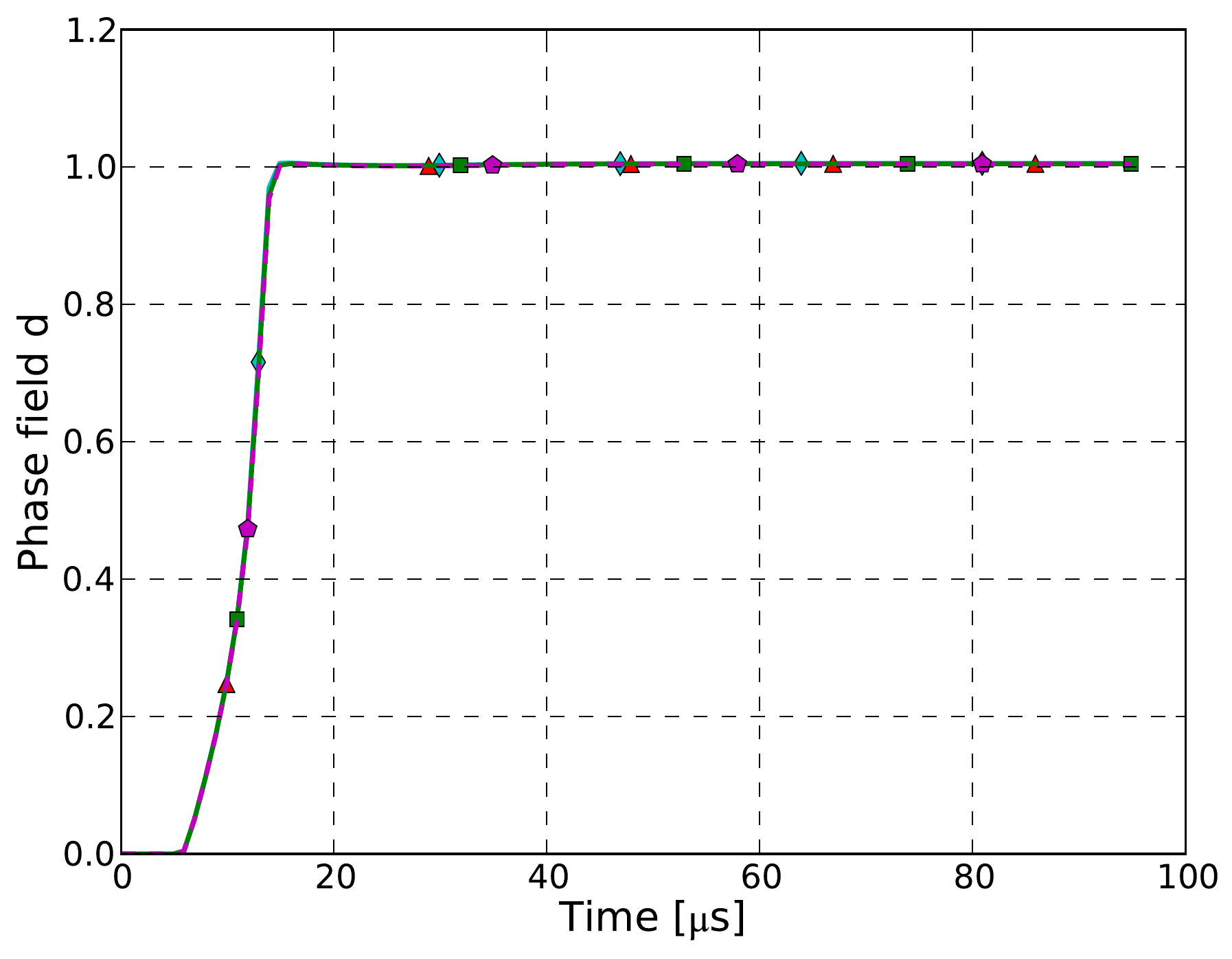}
}
\captionsetup{width=.90\textwidth}
\caption{ plot over time: (a) displacement magnitude at the upper right corner (point $Q$) (b) phase-field solution at crack tip (point $C$) } \label{fig:ForceDispPOTdPen_Ex2}
\end{figure}

% Plot over line (Pen)

\begin{figure}[p] 
\centering
\includegraphics[width=0.45\textwidth]{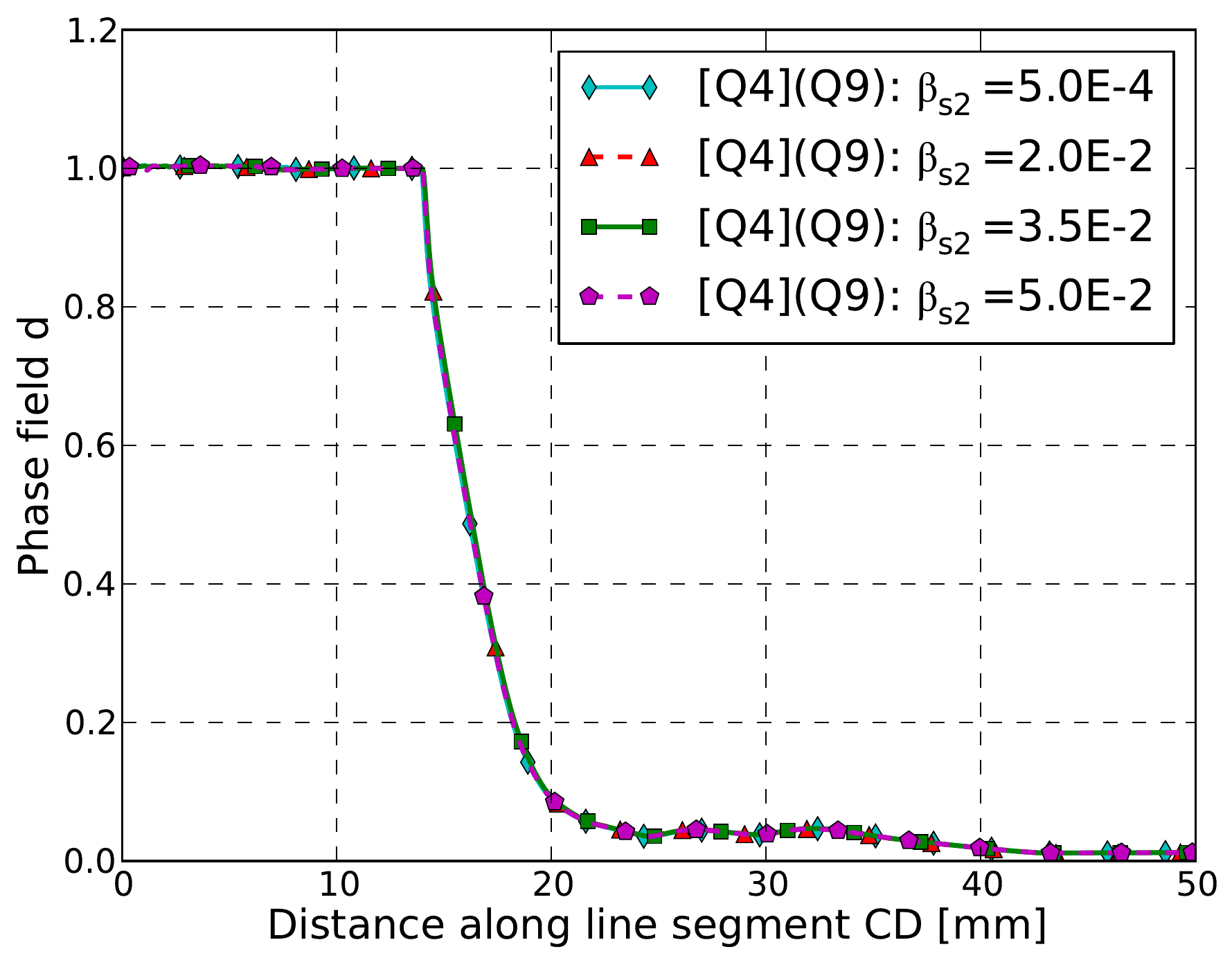}
\captionsetup{width=.90\textwidth}
\caption{ Phase-field solution along line segment CD }
\label{fig:POL_CD_Pen_Ex2}
\end{figure}

\begin{figure}[p] %[H] % "[t!]" placement specifier just for this example
\centering
\includegraphics[width=0.45\textwidth]{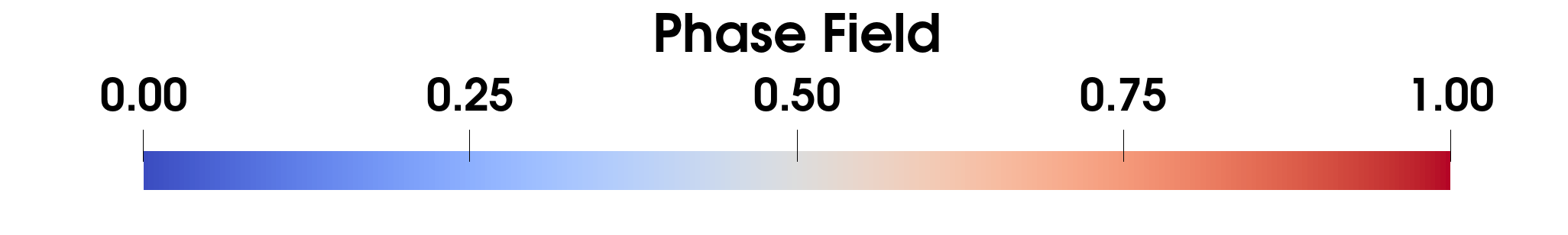}
\hspace{0.5cm}
\includegraphics[width=0.45\textwidth]{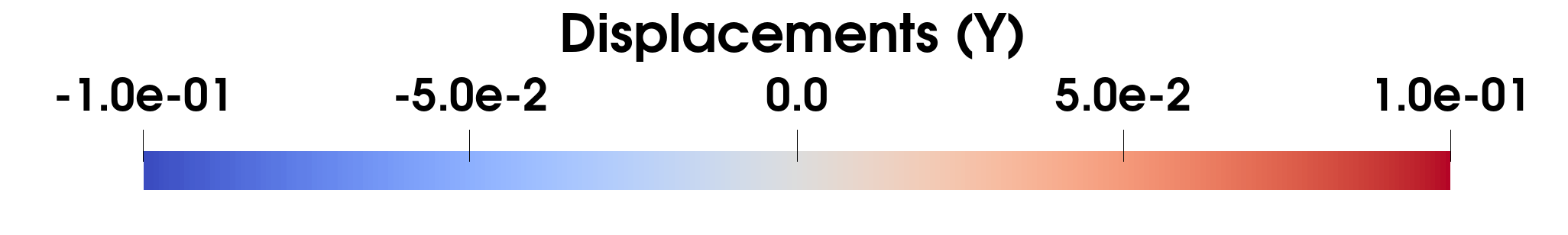}

\vspace{-0.35cm}

\centering
\subfloat[\label{fig:CDG_PF_Ex2} Scheme \text{[Q4](Q9)} - $\beta_{s2} = 5.0\times10^{-2}$]{%
\includegraphics[width=0.45\textwidth]{Section4/Example2/FinalSolution/Pr16M9Pf3/c95_crop.png}  %0.29
}
\hspace{1.0cm}
\subfloat[\label{fig:CDG_Disp_Ex2} Scheme \text{[Q4](Q9)} - $\beta_{s2} = 5.0\times10^{-2}$]{%
\includegraphics[width=0.45\textwidth]{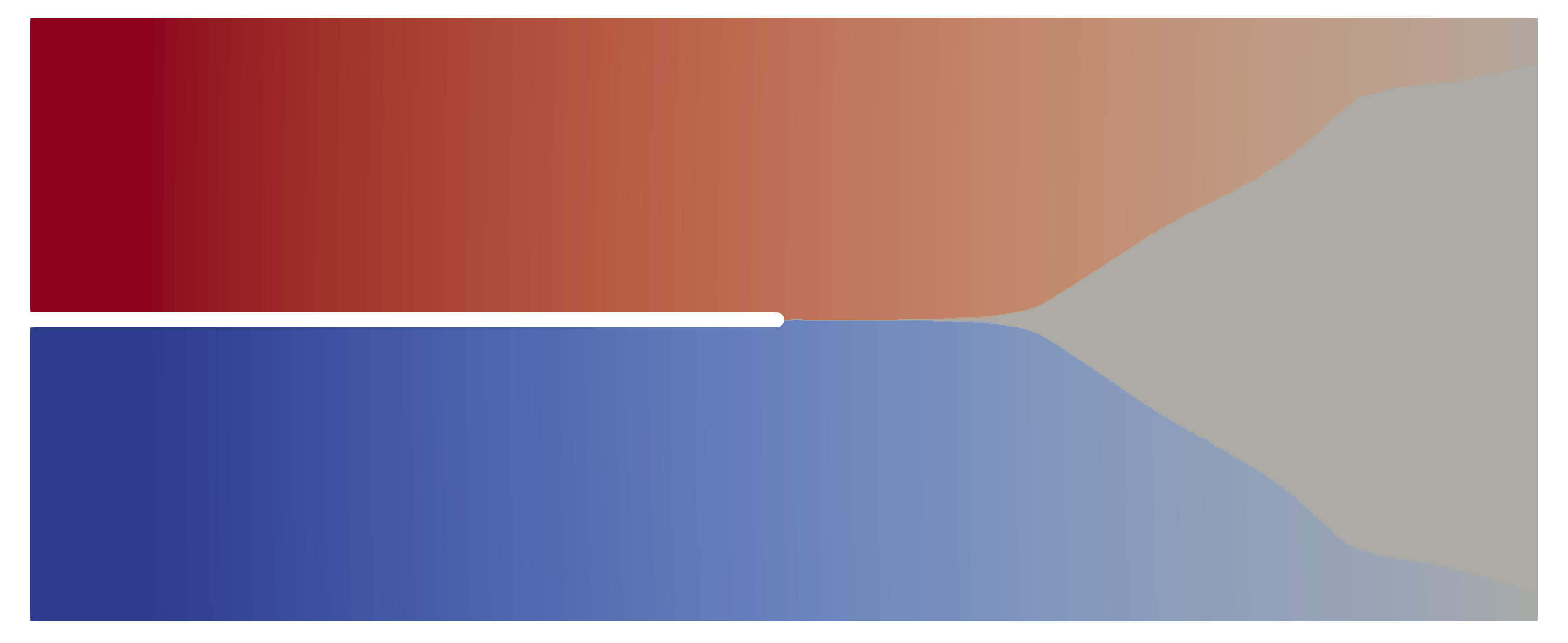}
}

%\medskip
\vspace{-0.2cm}

\centering
\subfloat[\label{fig:MixedQ9_PF_Ex2} Scheme \text{[Q4](Q9Q9)} ]{%
\includegraphics[width=0.45\textwidth]{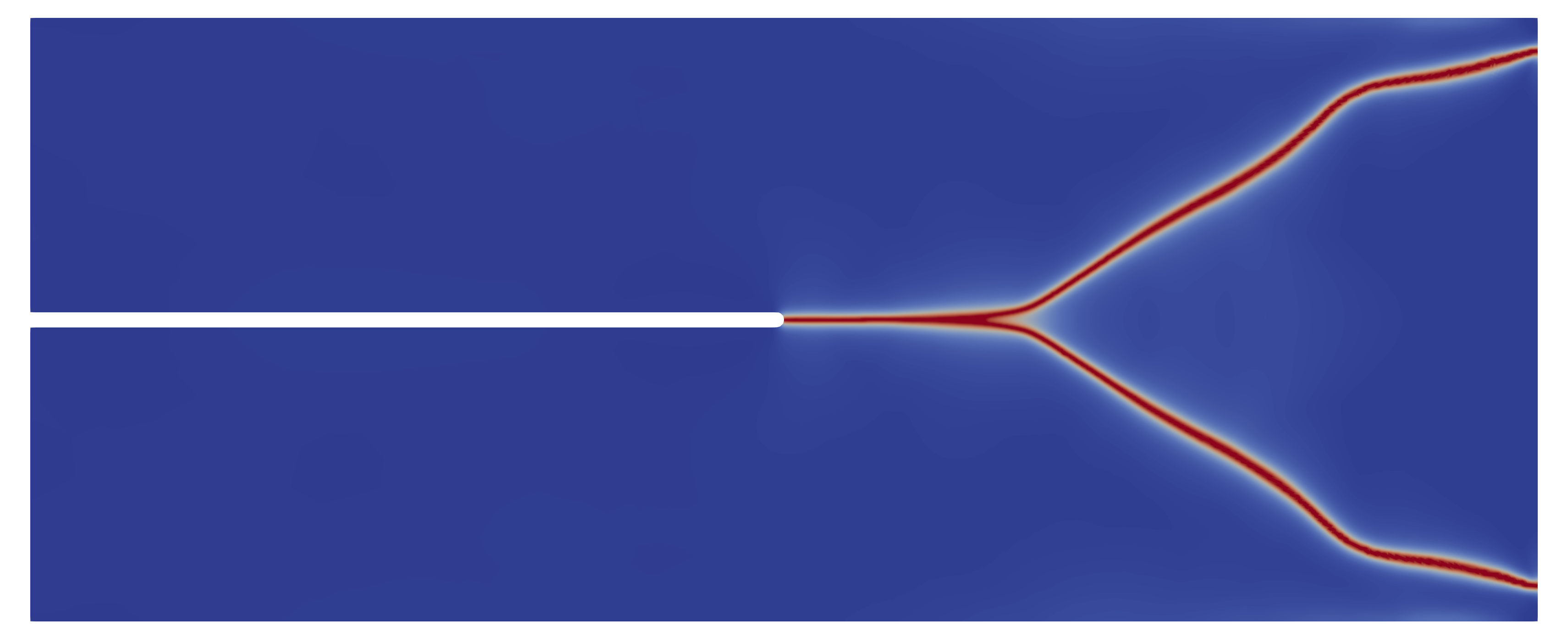}  %0.29
}
\hspace{1.0cm}
\subfloat[\label{fig:MixedQ9_Disp_Ex2} Scheme \text{[Q4](Q9Q9)} ]{%
\includegraphics[width=0.45\textwidth]{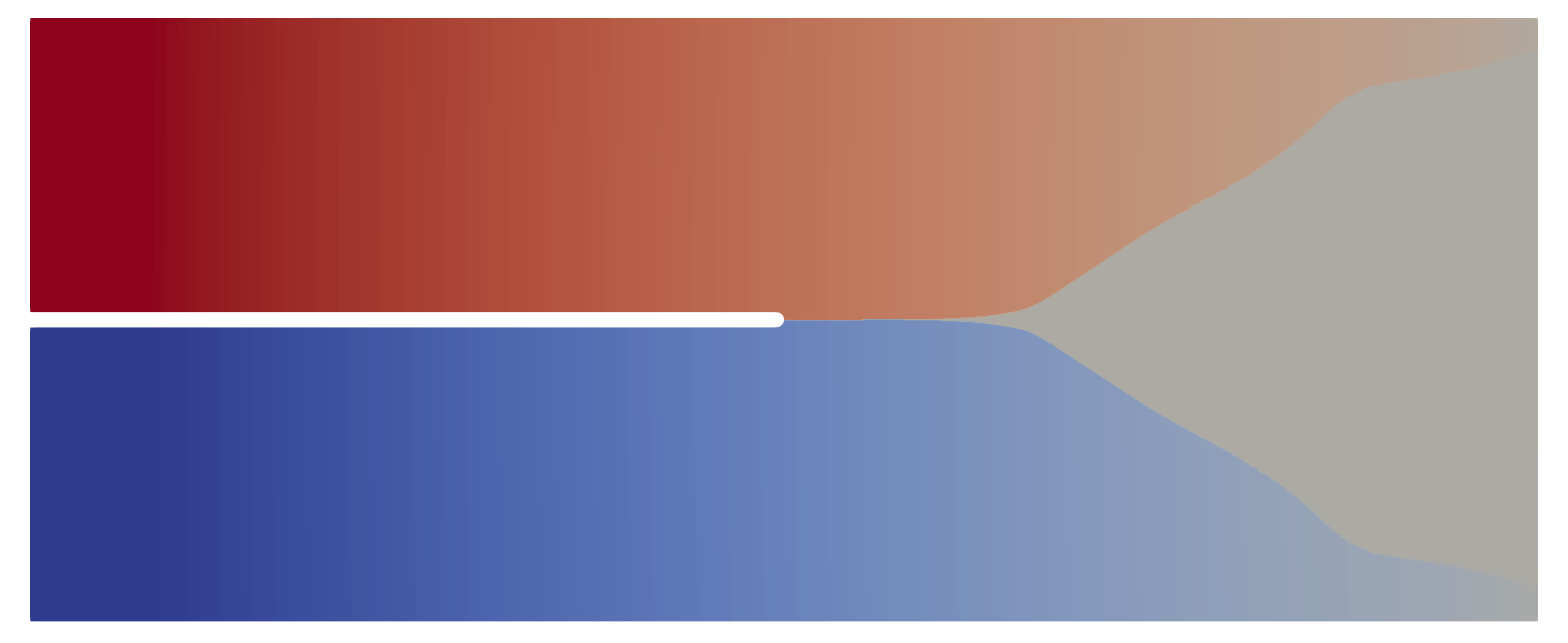}
}

%\medskip
\vspace{-0.2cm}

\centering
\subfloat[\label{fig:MixedQ4_PF_Ex2} Scheme \text{[Q4](Q4Q4)} ]{%
\includegraphics[width=0.45\textwidth]{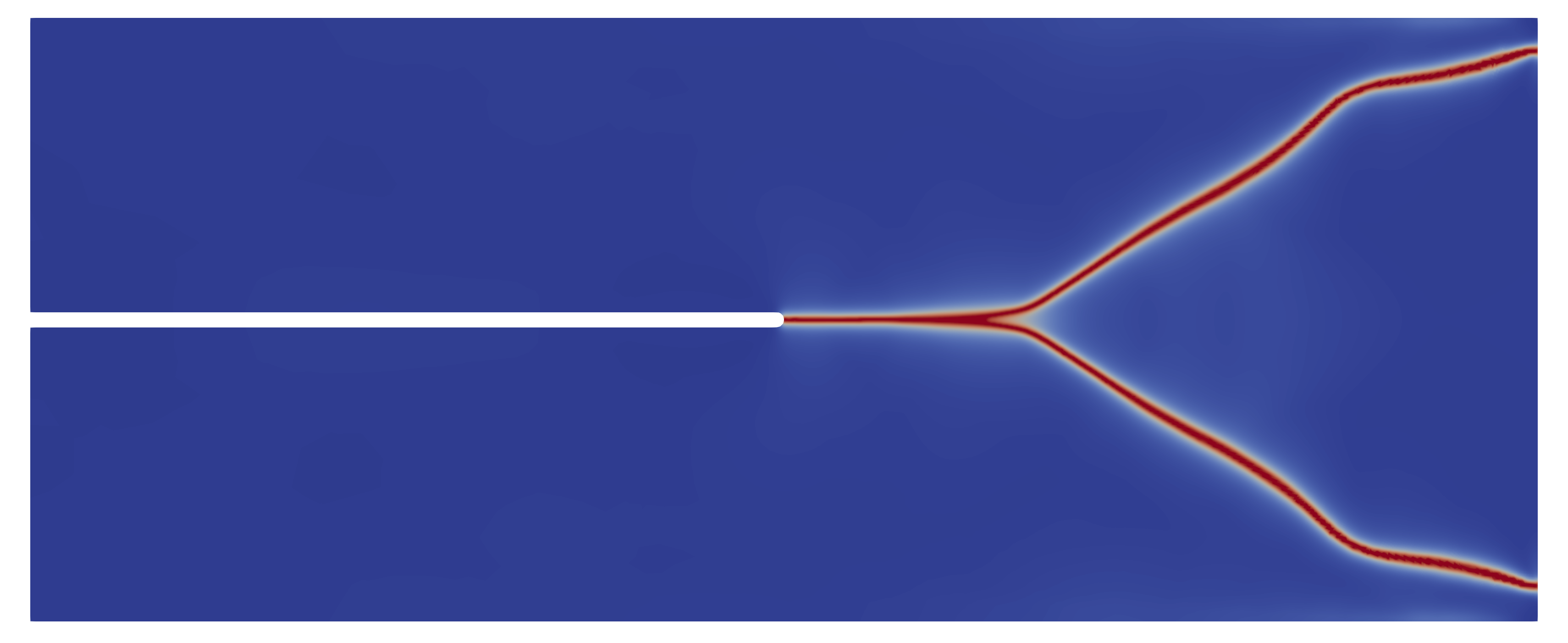}  %0.29
}
\hspace{1.0cm}
\subfloat[\label{fig:MixedQ4_Disp_Ex2} Scheme \text{[Q4](Q4Q4)} ] {%
\includegraphics[width=0.45\textwidth]{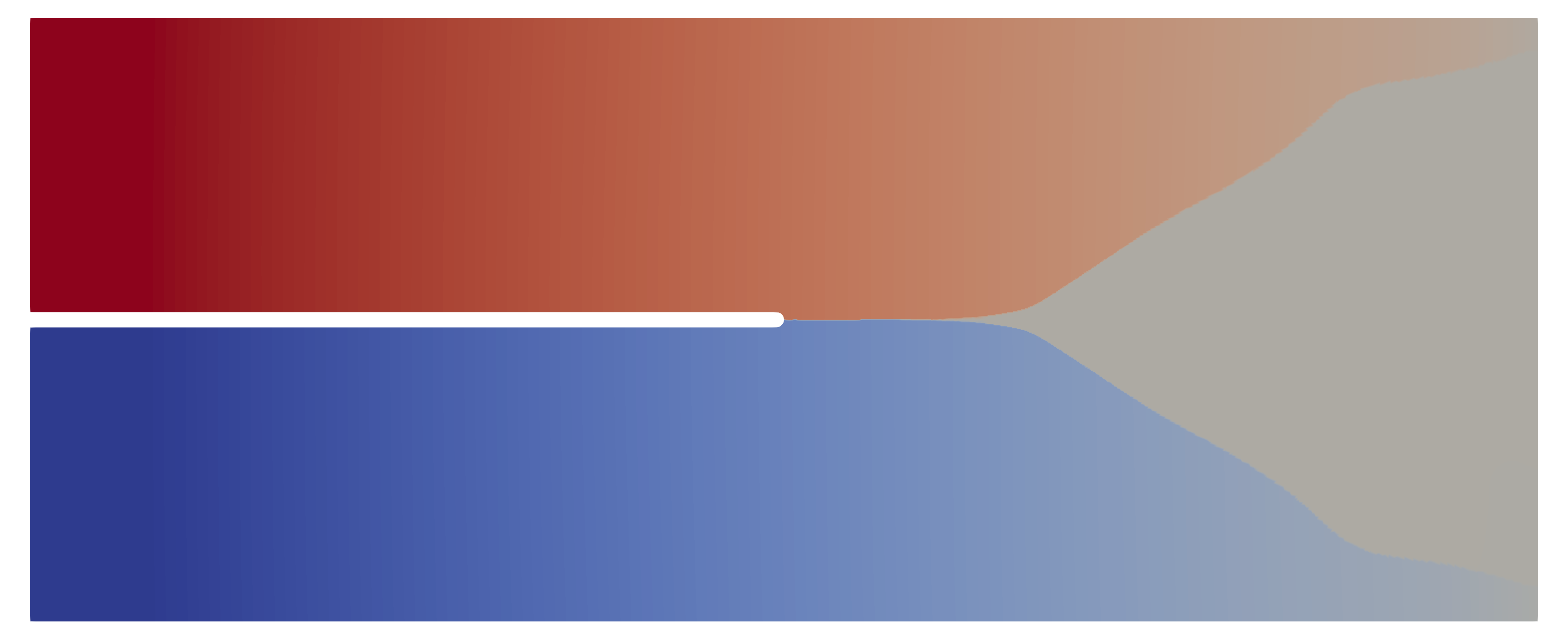}
}
\captionsetup{width=.90\textwidth}
\caption{ Comparison between the proposed Schemes (\text{[Q4](Q9)}, \text{[Q4](Q9Q9)} and \text{[Q4](Q4Q4)}) upon complete crack formation: (left) phase-field fracture paths (right) vertical displacement contours.   } \label{fig:AllMethodsPFDisp_Ex2}
\end{figure}

% Figure: Displacement and Phase-field Plot over time (Met)

\begin{figure}[t]  %[H] % "[t!]" placement specifier just for this example
\centering
\hspace{3.0cm}
\includegraphics[width=0.40\textwidth]{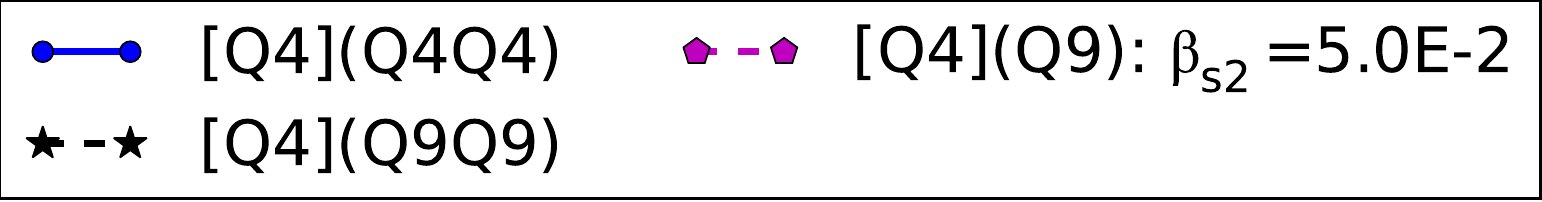}
\hspace{3.0cm}

\vspace{-0.35cm}

\centering
\subfloat[\label{fig:POTdispMet_Ex2}  displacement magnitude plot over time]{%
\includegraphics[width=0.45\textwidth]{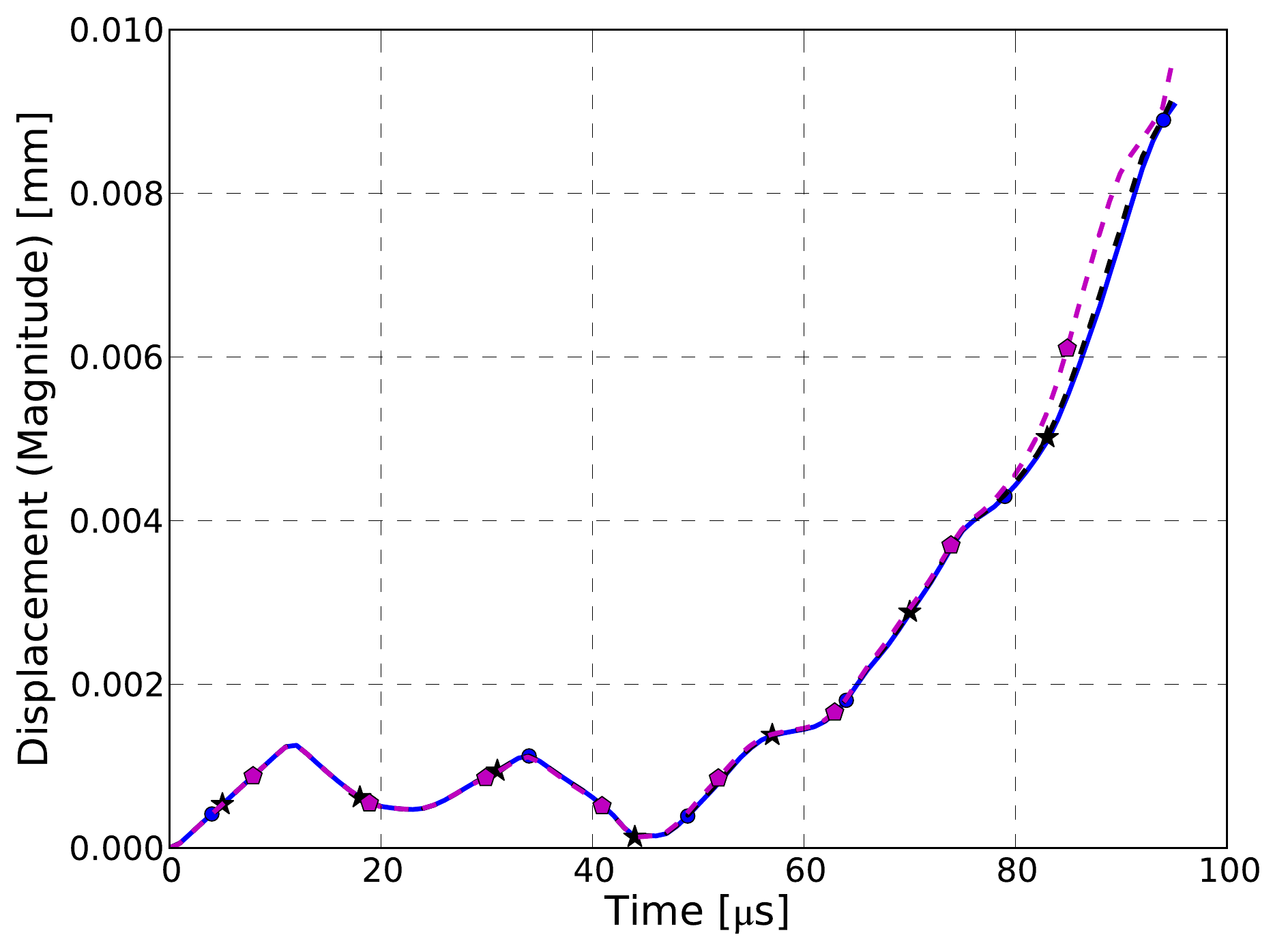}  %0.29
}
\hspace{1.0cm}
\subfloat[\label{fig:POTdMet_Ex2}  phase-field plot over time  ]{%
\includegraphics[width=0.45\textwidth]{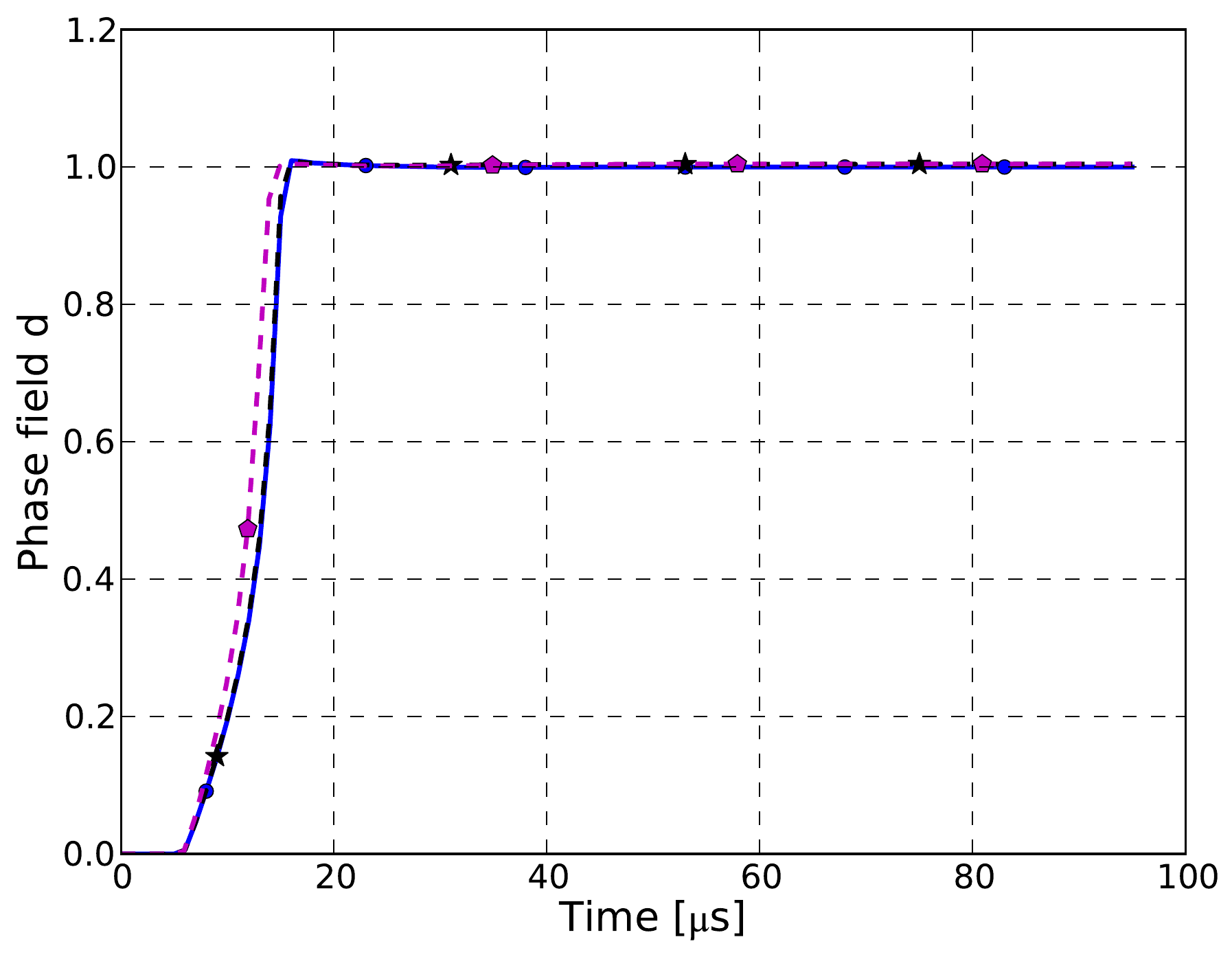}
}
\captionsetup{width=.90\textwidth}
\caption{ plot over time: (a) displacement magnitude at the upper right corner (point $Q$) (b) phase-field solution at notch tip (point $C$)  } \label{fig:ForceDispPOTdMet_Ex2}
\end{figure}

% Plot over line (Met)

\begin{figure}[t] 
\centering
\includegraphics[width=0.45\textwidth]{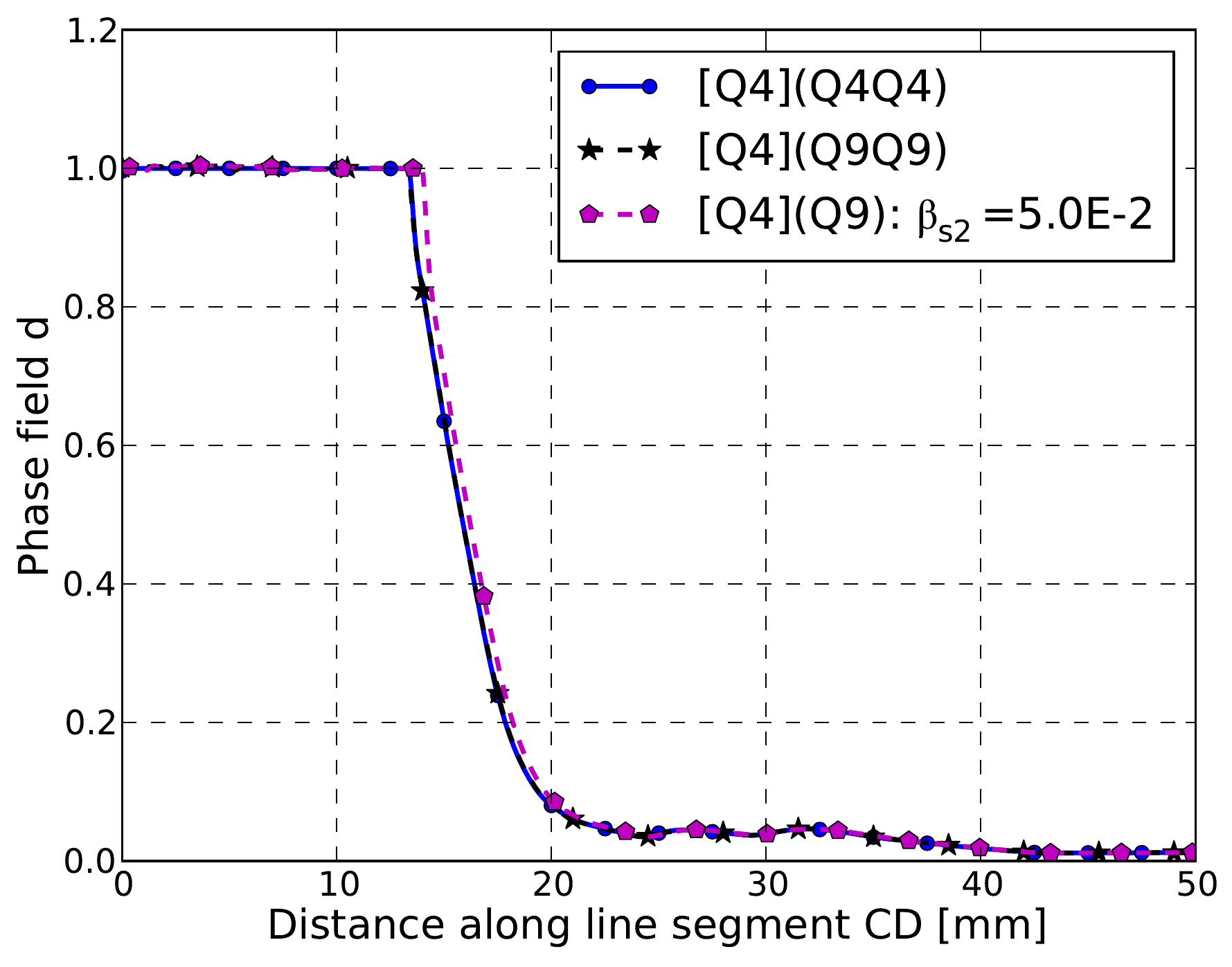}
\captionsetup{width=.90\textwidth}
\caption{ Phase-field solution along line segment CD }
\label{fig:POL_CD_Met_Ex2}
\end{figure}

Focusing on the influence of the penalty parameter on the displacement field, we examine the time evolution of the displacement magnitude, shown in \cref{fig:POTdispPen_Ex2}, at the upper right corner of the plate (point $Q$ in \cref{fig:Geometry_Ex2}). Note that the solutions are indistinguishable until $t=80.0\times10^{-6}$~s with negligible differences between them after that point. Then, we turn our attention to the influence of the penalty parameter on the predicted crack topology, and examine the time evolution of the phase field at the tip of the pre-existing notch (point $C$), presented in \cref{fig:POTdispPen_Ex2}. The excellent agreement between the different curves in this plot leads to the conclusion that crack initiation is insensitive to the penalty parameter value. A similar observation holds for crack branching: As can be seen from \cref{fig:POL_CD_Pen_Ex2}, crack branching takes place at the same point (located approximately $14.0$~mm from point $C$ along line segment $CD$; cf.\ \cref{fig:Geometry_Ex2}), at $t \approx 37.9\times10^{-6} $~s, regardless of the penalty parameter value. 

Based on the sensitivity analysis presented in this section, we choose the penalty parameter value $\beta_{s2} = 5.0\times10^{-2}$ for comparison with the mixed FEM schemes in the next section. The analysis results demonstrate that excessively low values of the penalty parameter lead to a loss of symmetry in the solution of the present symmetric crack branching problem, whereas excessively high values result in numerical difficulties, highlighting once more the importance of selecting the proper penalty parameter value, despite the added computational effort required to do so.

%, because symmetric and insensitive results are found.

\subsubsection{Comparison between numerical schemes}

A comparison between the numerical results obtained using the three schemes under consideration is presented in \cref{fig:AllMethodsPFDisp_Ex2}. More specifically, phase-field and vertical displacement contours at $t = t_f = 95.0\times10^{-6}$~s are shown, respectively, in \cref{fig:CDG_PF_Ex2,fig:CDG_Disp_Ex2} for Scheme \text{[Q4](Q9)}, in \cref{fig:MixedQ9_PF_Ex2,fig:MixedQ9_Disp_Ex2} for Scheme \text{[Q4](Q9Q9)}, and in \cref{fig:MixedQ4_PF_Ex2,fig:MixedQ4_Disp_Ex2} for Scheme \text{[Q4](Q4Q4)}. Clearly, there is excellent agreement between all three numerical schemes in terms of the final crack topology predicted. It is also noteworthy that all three schemes agree very well with regards to the discontinuity in vertical displacement across the crack, and more generally, the distribution of vertical displacements after crack formation, including in the region of the plate that is separated due to crack branching (shaded in gray color in \cref{fig:AllMethodsPFDisp_Ex2}) where these displacements are zero.

To assess the differences in mechanical response as predicted by the different numerical schemes, the time evolution of the displacement magnitude at point $Q$ is plotted in \cref{fig:POTdispMet_Ex2}. Note that the oscillations that take place over the first $50\times10^{-6}$~s reflect the dynamic response of the plate. These oscillations (stress waves) are responsible for crack branching. We note that the same behavior is captured by all schemes almost throughout the simulation; more specifically, the curves coincide until $t=80.0\times10^{-6}$~s, after which the C/DG Scheme \text{[Q4](Q9)} starts deviating slightly from the other two schemes. To study this deviation, we examine \cref{fig:POTdMet_Ex2} which shows the time evolution of the phase field at the notch tip (point $C$), and conclude that Scheme \text{[Q4](Q9)} predicts slightly earlier crack initiation in comparison to the other schemes. This suggests a potential sensitivity of the C/DG scheme to the mesh topology around point $C$. In addition, it leads the C/DG scheme to predict crack branching slightly further to the right along line segment $CD$, as can be seen from \cref{fig:POL_CD_Met_Ex2} which shows a plot of the phase-field solution along line segment $CD$ at $t=95.0\times10^{-6}$~s (when crack branches reach the fully-developed stage). 

Finally, we note that, in contrast to the previous example problem, no solution differences that can be attributed to the staggered approach are observed here. This is because a constant time step size is used throughout all the simulations described in the present section, thus reducing the interplay between the time integration procedure and the convergence behavior of each spatial discretization scheme.

%Hence, the crack paths in Scheme \text{[Q4](Q9)} are in a good agreement with the results of Schemes \text{[Q4](Q4Q4)} and \text{[Q4](Q9Q9)}. 

%The different spatial discretizations control the convergence behavior at each time step. However, an adaptive staggered approach is sensitive to the solution at each time step. This sensitivity increases the solution dissimilarities among the numerical schemes.

%Although the crack initiation is insensitive to penalty parameters, the crack branches and paths are 

%keeping the time step constant we reduce the impact of the staggered approach
%constant time step is kept, and consequently, the time integration is not affected by the spatial convergence behavior. 

%Given that an adaptive time-stepping method is used to advance the system in time, the aforementioned differences are mainly caused by the staggered approach and its sensitivity to the spatial discretization.

%\subsubsection{Scalability analysis}

% Third Study

%\subsection{Convergence analysis: Method of Manufactured Solutions}

\section{Conclusions}
\label{sec:Conclusions}

We have derived a fourth-order phase-field theory of fracture in two independent ways (from Hamilton's principle and a micromechanics-based approach). We have shown that the latter approach is more general since (i) it is equivalent to the former in the case of a system with no dissipative mechanisms, (ii) dissipative terms can be naturally included, and (iii) it provides a physical interpretation of the boundary conditions in high-order damage theories.

We have coupled the fourth-order phase-field equation with the momentum balance equation to model dynamic fracture of hyper-elastic materials using a volumentric/deviatoric decomposition of the elastic energy. To solve the momentum equation, the standard finite element method is adopted with bilinear shapes functions.

To solve the fourth-order phase-field equation, special numerical methods are required to enforce $C^1$-continuity of the unknown phase field. In this work, we have used a C/DG method where the key idea is that Lagrange shape functions provide $C^0$-continuity to the solution, while gradient discontinuities at inter-element boundaries are prevented and the required higher-order regularity conditions are enforced with the aid of additional variational and interior penalty terms in the weak form of the problem. The C/DG method is compared with mixed FEM, where the order of the PF equation is reduced by treating the Laplacian of the PF as an additional unknown. The conclusions of the present work are as follows:
\begin{itemize}
\item The C/DG method establishes a middle ground between first-order and second-order mixed FEMs. The number of degrees of freedom are decreased compared to the second-order mixed FEM while maintaining second-order accuracy of the interpolation functions.
\item The C/DG method requires additional effort to tune the penalty parameters and produce accurate results.
\item For low values of penalty parameters, regularity conditions are not enforced in a sufficiently accurate manner, leading to instability and/or inaccurate solutions that may be sensitive to these parameter values. This is well understood from the convergence analysis of C/DG methods. However, there is a range of penalty parameters in which numerical solutions with insensitivity to this parameter are achieved. On the other hand, large penalty values lead to ill-conditioning of the resulting equation system and are thus detrimental to numerical stability and robustness.
\item The results of the C/DG method are in very good agreement with the results of mixed FEMs in terms of crack initiation and branching. Small deviations are observed when the total crack length increases with C/DG results displaying a slight sensitivity to the quality of the mesh used.
\end{itemize}

In the present paper, the focus is placed on brittle materials, allowing a hyperelastic material model to be used in representing the stress response of the undamaged material. From a constitutive modeling perspective, this keeps complexity at a tractable level, and allows us to sidestep non-trivial issues related to the coupling between plasticity and damage. Such issues will be examined in future work that seeks to generalize the approach presented herein to the case of dynamic ductile fracture.

\section*{Acknowledgments}
The work presented in this article was funded by the Laboratory Directed Research and Development (LDRD) program of Los Alamos National Laboratory under project number 20220129ER. The authors gratefully acknowledge this support. Los Alamos National Laboratory is operated by Triad National Security, LLC, for the National Nuclear Security Administration of U.S. Department of Energy (Contract No.\ 89233218CNA000001).

\clearpage
\appendix
\renewcommand{\theHsection}{A\arabic{section}}
\appendix

\section{Details on derivation from Hamilton's principle }
\label{apdx:HamiltonDerivationDetails}

In this appendix, we provide details on the way we derive governing equations and boundary conditions from Hamilton's principle. Applying integration by parts and Green's second identity, \cref{eq:IntegrandHamiltonPrinciple} is rewritten as
\begin{equation} \label{eq:IntegrandAfterIntbyPartsandGreen}
\begin{split}
\delta \mathcal{L} + \delta \mathcal{W}_{ext} & =  \CompInt{\RefDom}{\rho_0\dt{\DispRef} \cdot \delta \dt{\DispRef}}{V} + \CompInt{\RefDom}{ \MicroDens  \dt{\PF} \, \delta \dt{\PF} }{V}  \\
& + \CompInt{\RefDom}{ \Div \left( \partder{\DegElastEnerDen}{\DeformGrad} \right) \cdot \delta \DispRef }{V} - \CompInt{\ExtBoun_N^\DispRef}{ \left( \partder{\DegElastEnerDen}{\DeformGrad} \cdot \NormVecRef \right) \cdot  \delta \DispRef}{A} - \CompInt{\RefDom}{ \partder{\DegElastEnerDen}{\PF} \virtPF }{V}   \\
& - \CompInt{\RefDom}{ G_c \partder{\gamma_{\lzr}}{\PF} \virtPF}{V} + \CompInt{\RefDom}{ G_c \Div \left(\partder{\gamma_{\lzr}}{\PFgrad} \right)   \virtPF }{V}-  \CompInt{\partial\RefDom}{ G_c   \partder{\gamma_{\lzr}}{\PFgrad} \cdot \NormVecRef   \virtPF }{A} \\
&- \CompInt{\RefDom}{ G_c \Lapl \left( \partder{\gamma_{\lzr}}{\PFlap} \right)  \virtPF }{V} - \CompInt{\partial \RefDom}{ G_c  \partder{\gamma_{\lzr}}{\PFlap}  \delta \PFgrad \cdot  \NormVecRef }{A} + \CompInt{\partial \RefDom}{ G_c  \Grad \left( \partder{\gamma_{\lzr}}{\PFlap} \right) \cdot \NormVecRef    \virtPF }{A}  \\
& + \CompInt{\RefDom}{ \BodyForceRef \cdot \delta \DispRef }{V} + \CompInt{\ExtBoun_N^\DispRef}{ \SurfTractionRef \cdot \delta \DispRef }{A}, \\
\end{split}
\end{equation}
where $\Div$ denotes the divergence operator with respect to the reference configuration, and $\NormVecRef$ is the outward normal vector on the boundary $\partial \RefDom$. Substituting \cref{eq:IntegrandAfterIntbyPartsandGreen} into \cref{eq:HamiltonPrinciple} and taking into account that $t_1,t_2$ are arbitrary times, the extended Hamilton principle in the system takes the form 
\begin{equation}
\begin{split}
& \SpaceTimeInt{\RefDom}{- \rho_0\ddt{\DispRef} + \Div \left( \partder{\DegElastEnerDen}{\DeformGrad} \right) + \BodyForceRef}{ \cdot \, \delta \DispRef}{V}  + \SpaceTimeInt{\ExtBoun_N^\DispRef} { \SurfTractionRef -  \partder{\DegElastEnerDen}{\DeformGrad} \cdot \NormVecRef   }{ \cdot \, \delta \DispRef}{A} \\
& + \SpaceTimeInt{\RefDom}{ - \MicroDens  \ddt{\PF} - \partder{\DegElastEnerDen}{\PF} - G_c \partder{\gamma_{\lzr}}{\PF} + G_c \Div \left(\partder{\gamma_{\lzr}}{\PFgrad} \right) - G_c \Lapl \left( \partder{\gamma_{\lzr}}{\PFlap} \right) }{ \virtPF}{V} \\
& \SpaceTimeInt{\partial\RefDom} { G_c   \Grad \left( \partder{\gamma_{\lzr}}{\PFlap} \right) -  G_c \partder{\gamma_{\lzr}}{\PFgrad}  }{\cdot \NormVecRef   \virtPF }{A} - \SpaceTimeInt{\partial\RefDom} { G_c  \partder{\gamma_{\lzr}}{\PFlap}   }{ \delta  \PFgrad \cdot \,  \NormVecRef }{A} = 0, \\
\end{split}
\end{equation}

Given that $\delta \DispRef$, $\virtPF$, and $\delta  \PFgrad \cdot \NormVecRef $ are virtual fields of the system, the governing equations and boundary conditions are derived as shown in \cref{eq:GoverningEqsHamilton,eq:BoundaryConditionsHamilton} respectively.

\section{Details on derivation from a micromechanics-based approach }
\label{apdx:MicromechanicalDerivationDetails}

In this part, we provide details on the way we derive the energy balance equation for materials with high-order stresses. Using an arbitrary open set $\ArbDom \subset \RefDom$ (with a smooth boundary $\partial \ArbDom$) as a control volume, the conservation of energy balance can be stated as follows
\begin{equation} \label{eq:EnergyBalanceIntergralForm}
\begin{split}
\dt{\mathcal{K}}_\ArbDom + \int_\ArbDom \rho_0 \dt{\hat{e}}  \ud V  & =  \int_{\ArbDom} \BodyForceRef \cdot \dt{\DispRef} \ud V + \int_{\partial \ArbDom}  ( \tens{P} \cdot \NormVecRef ) \cdot  \dt{\DispRef} \ud A \\
& + \int_{\partial \ArbDom} ( \thickbar{\vect{\Sigma}} \cdot \NormVecRef) \, \dt{\PF}    \ud A  
- \int_{\partial \ArbDom} ( \nabla \Phi \cdot \NormVecRef) \, \dt{\PF}    \ud A
+ \int_{\partial \ArbDom} \Phi \,  \nabla \dt{\PF} \cdot  \NormVecRef  \ud A ,    
\end{split}
\end{equation}
where $\mathcal{K}_\ArbDom$ denotes the restriction of the kinetic energy to the control volume $\ArbDom$, the first two terms on the right side represent the rate of work (power) of macro-forces, and the remaining three terms represent the power expended across $\partial \ArbDom$ due to micro-forces. Specifically, the low-order micro-stresses are power-conjugate to $\dt{\PF}$, while the higher-order ones to $\nabla \dt{\PF}$.

%micro-mechanical variables are ($\dt{\PF}$, $\tens{\Sigma}$) and ($\nabla \dt{\PF}$, $\Phi$). 

%Specifically, the pairs of power-conjugate micro-mechanical variables are ($\dt{\PF}$, $\tens{\Sigma}$) and ($\nabla \dt{\PF}$, $\Phi$). 

To obtain the local form of the energy balance equation, we follow the next steps. First, we take the time derivative of the kinetic energy expression \eqref{eq:KineticEnery}
\begin{equation} \label{eq:KineticEneryRestrictionDerivative}
\dt{\mathcal{K}}_\ArbDom = \int_\ArbDom \rho_0 \dt{\DispRef} \cdot \ddt{\DispRef}  \ud V + \int_\ArbDom \MicroDens \dt{\PF} \ddt{\PF}  \ud V .
\end{equation}
Second, we multiply \cref{eq:MacrobalanceLocalLagr} by $\dt{\DispRef}$ and integrate by parts over the domain $\ArbDom$ as follows
\begin{equation} \label{eq:MacrobalanceLocalLagrIntByParts}
\int_\ArbDom \rho_0 \dt{\DispRef} \cdot \ddt{\DispRef}  \ud V  = - \int_\ArbDom \tens{P} : \dt{\DeformGrad}  \ud V + \int_{\partial\ArbDom} ( \tens{P} \cdot \NormVecRef ) \cdot  \dt{\DispRef}  \ud A + \int_{\ArbDom} \BodyForceRef \cdot \dt{\DispRef}  \ud V .
\end{equation}
Similarly, we multiply \cref{eq:MicrobalanceLocalLagrHighOrd} by $\dt{\PF}$ and integrate by parts twice
\begin{equation} \label{eq:MicrobalanceLocalLagrHighOrdIntByParts}
\begin{split}
\int_\ArbDom \MicroDens \dt{\PF} \ddt{\PF}  \ud V  & = - \int_\ArbDom \Phi \Delta \dt{\PF}  \ud V - \int_{\partial \ArbDom} ( \nabla \Phi \cdot \NormVecRef) \, \dt{\PF}    \ud A + \int_{\partial \ArbDom} \Phi \,  \nabla \dt{\PF} \cdot  \NormVecRef  \ud A \\
& + \int_{\partial \ArbDom} ( \thickbar{\vect{\Sigma}} \cdot \NormVecRef) \, \dt{\PF}    \ud A - \int_{\ArbDom} \thickbar{\vect{\Sigma}} \cdot \nabla \dt{\PF}   \ud V + \int_{\ArbDom}  M^* \dt{\PF} \ud V. 
\end{split}
\end{equation}
Finally, we subtract \cref{eq:KineticEneryRestrictionDerivative,eq:MacrobalanceLocalLagrIntByParts,eq:MicrobalanceLocalLagrHighOrdIntByParts} from the energy balance equation (\ref{eq:EnergyBalanceIntergralForm}) and after simplifications we arrive at
\begin{equation}
\int_\ArbDom \rho_0 \dt{\hat{e}}  \ud V=  \int_\ArbDom \tens{P} : \dt{\DeformGrad} \ud V + \int_\ArbDom \thickbar{\vect{\Sigma}} \cdot \nabla \dt{\PF} \ud V - \int_\ArbDom M^* \dt{\PF} \ud V + \int_\ArbDom \Phi \Delta \dt{\PF}  \ud V ,
\end{equation}
which yields the local form \eqref{eq:LocalEnergyBalance}.

\section{Equivalence of boundary conditions of phase-field equations}
\label{apdx:BoundaryconditionsPFequation}

The boundary conditions of phase-field equation are identical in both formulations (presented in \cref{ssec:HamiltonDerivationGovEqFourthOrderTheory} and \cref{ssec:MicroMechanicsBasedInterpretationFourthOrderT}) as one can see by comparison with the following derivation
\begin{equation} \label{eq:BoundaryConditionsThermodyDeriv}
\left. \begin{aligned}
\left( - \nabla \Phi +  \thickbar{\vect{\Sigma}} \right) \cdot \NormVecRef =0 &     \\
\Phi = 0   &   
\end{aligned} \right\} \Rightarrow
\left\{ \begin{aligned}
- \left[ \nabla \left( \partder{\SurfEnerDenFun}{\PFlap} \right) - \partder{\SurfEnerDenFun}{\PFgrad}  \right] \cdot \NormVecRef  =0 &     \\
\partder{\SurfEnerDenFun}{\PFlap} = 0   &   
\end{aligned} \right\}  \text{ on } \partial \RefDom .
\end{equation}
For the special case of the fourth-order crack density function (\cref{eq:CrackDensityFunction}), the boundary conditions read
\begin{equation} \label{eq:BoundaryCondPhaseFieldHighOrderThermo}
\left. \begin{aligned}
\nabla \left(  \alpha_2 \PFlap + \alpha_1  \PF  \right) \cdot \vect{N}  =0 \Rightarrow  \nabla \left(  \lzr^2 \PFlap- 2  \PF  \right) \cdot \vect{N} & =0   \\
  \PFlap & = 0  
\end{aligned} \right\}  \text{ on } \partial \RefDom  \,
\end{equation}
where the surface energy density function $\SurfEnerDenFun$ is given in \cref{eq:DefineSurfEnerDenFun}.

\section{Linearization of the variational formulation}
\label{apdx:Linearization}

In this part, we briefly present the linearization of the variational formulation which is expressed in \cref{eq:WeakFormLinearMomBalance}. 
The linearization is stated for the weak form with respect to the initial configuration.

Denoting $\bar{\delta}\vect{u}$ the increment of displacement field and assuming that the external load is conservative, the directional derivative of the residual $R_\DispRef$ in the direction of $\bar{\delta}\vect{u}$ can be calculated by considering only the term $\innerprod{\delta \tens{E}_h}{\tens{S}}{\RefDom}$. Specifically, the linearization yields
\begin{equation} \label{eq:LinearAppdTwoCompo}
\tens{J}_{\vect{u}\vect{u}} = D \, R_\DispRef(\vect{w}_h,\motion) \cdot \bar{\delta}\vect{u}= \innerprod{\nabla \vect{w}_h}{\tens{S}  \nabla \bar{\delta}\vect{u} }{\RefDom} + \innerprod{\delta \tens{E}_h}{\tenf{C} : \bar{\delta}\tens{E} }{\RefDom} ,
\end{equation}
where the increment of the Green-Lagrange strain tensor is given by
\begin{equation}
\bar{\delta} \tens{E} = \frac{1}{2} \left( \tens{F}\transpose \Grad \bar{\delta}\vect{u} + \Grad\transpose \bar{\delta}\vect{u} \tens{F} \right) ,
\end{equation}
and the \emph{degraded} elasticity tensor
\begin{equation}
\tenf{C} = 4 \frac{\partial^2  \DegElastEnerDen}{ \partial \tens{C} \partial \tens{C}} = 2 \partder{\tens{S}}{\tens{C}}  =  2 \partder{\tens{S}_e}{\tens{C}} + \left[ g(\PF) - 1 \right]  \left[ 2 \partder{\tens{S}^+}{\tens{C}}  + 2 \partder{\thickbar{\tens{S}}}{\tens{C}}   \right] .
\end{equation}
Notice that the linearization in \cref{eq:LinearAppdTwoCompo} has two terms: (i) initial stress term and (ii) a term associated with the \emph{degraded} elasticity tensor (i.e., incremental constitutive tensor). The elasticity tensor is decomposed into volumetric and isochoric components, $ \tenf{C}_\circ$ and  $\thickbar{\tenf{C}}$ respectively.
\begin{equation}
\tenf{C}_e = \tenf{C}_\circ + \thickbar{\tenf{C}}
\end{equation}
Based on the latter decomposition, the dilative term associated with $U^+(J)$, defined in \cref{eq:DilativeTermDef0}, is express as
\begin{equation}
\tenf{C}^+ = H(J-1) \tenf{C}_\circ .
\end{equation}
In compact form, the \emph{degraded} elasticity tensor is given as follows
\begin{equation}
\tenf{C} = \tenf{C}_e +  \left[ g(\PF) - 1 \right]  \left[ \tenf{C}^+ + \thickbar{\tenf{C}} \right] ,
\end{equation}
where
\begin{equation} \label{eq:DefinitionsOfConstTensor}
\left. \begin{aligned}
 \tenf{C}_\circ  &=2 \partder{\tens{S}_\circ }{\tens{C}} \\
\thickbar{\tenf{C}} &= 2 \partder{\thickbar{\tens{S}}}{\tens{C}} 
\end{aligned} \right\}.
\end{equation}
% \tenf{C}_e &=2 \partder{\tens{S}_e}{\tens{C}}  \\

Plugging \cref{eq:StressesInCompactFormandExplan} into \cref{eq:DefinitionsOfConstTensor} and after some algebraic manipulations, the volumetric component reads
\begin{equation}
 \tenf{C}_\circ = \left[ J U^\prime(J) + J^2 U^{\prime\prime} (J) \right]  \tens{C}^{-1} \otimes  \tens{C}^{-1} -2 J U^\prime (J)  \tenf{I}_{\tens{C}^{-1}}
\end{equation}
where the tensor $ \tenf{I}_{\tens{C}^{-1}}$ is introduced to shorten the notation of the derivative of the inverse right Cauchy-Green tensor as follows
\begin{equation}
\partder{\tens{C}^{-1}_{IJ}}{\tens{C}_{KL}}  = - \left( \tenf{I}_{\tens{C}^{-1}} \right)_{IJKL}= - \frac{1}{2} \left( \tens{C}^{-1}_{IK} \tens{C}^{-1}_{JL} + \tens{C}^{-1}_{IL} \tens{C}^{-1}_{JL} \right) .
\end{equation}
Similarly, the isochoric component is calculated as follows
\begin{equation}
\thickbar{\tenf{C}} = -  \frac{2}{3} \left[ \tens{C}^{-1} \otimes \thickbar{\tens{S}} + \thickbar{\tens{S}} \otimes  \tens{C}^{-1}    \right]  - \frac{2}{3} J^{-\frac{2}{3}} \left( 2 \partder{\thickbar{W}}{\thickbar{\tens{C}}} : \tens{C} \right) \left[ \frac{1}{3}  \tens{C}^{-1} \otimes  \tens{C}^{-1} - \tenf{I}_{\tens{C}^{-1}}   \right] + 4 \tenf{Q} : \frac{\partial^2  \thickbar{W}}{ \partial \thickbar{\tens{C}} \partial \thickbar{\tens{C}}} : \tenf{P}  ,
\end{equation}
where the fourth-order tensors $\tenf{P}$ and $\tenf{Q}$ have symmetries expressed as follows
\begin{equation}
\tenf{Q}_{IJKL} =\tenf{P}_{KLIJ} = J^{-\frac{2}{3}} \left( \tenf{I}_{IJKL} - \frac{1}{3} \tens{C}^{-1}_{IJ} \tens{C}_{KL} \right) ,
\end{equation}
and the fourth-order identity tensor is given
\begin{equation}
\partder{\tens{C}_{IJ}}{\tens{C}_{KL}} =\tenf{I}_{IJKL} = \frac{1}{2} \left(\delta_{IK} \delta_{JL} + \delta_{IL} \delta_{JK} \right)
\end{equation}
where $\delta_{IJ}$ denotes the second-order identity tensor.

For the special form of volumetric and deviatoric contributions given in \cref{eq:VolDevContributionsSpec}, the components of the elasticity tensor read
\begin{equation}
\tenf{C}_\circ = \kappa ( J^2 -1 ) \left(  \tens{C}^{-1} \otimes  \tens{C}^{-1} - \tenf{I}_{\tens{C}^{-1}}   \right) ,
\end{equation}
and
\begin{equation}
\thickbar{\tenf{C}} = -  \frac{2}{3} \left[ \tens{C}^{-1} \otimes \thickbar{\tens{S}} + \thickbar{\tens{S}} \otimes  \tens{C}^{-1}    \right]  - \frac{2}{3} \mu J^{-\frac{2}{3}} \trace{\tens{C}}  \left[ \frac{1}{3}  \tens{C}^{-1} \otimes  \tens{C}^{-1} - \tenf{I}_{\tens{C}^{-1}}   \right]  .
\end{equation}

\section{Consistency of C/DG Method}
\label{apdx:CDG_Consistency}

Using integration by parts multiple times in \cref{eq:CDGWeakFormCDGMethod} and taking $\beta_{s2}=\beta_s \alpha_2$, we arrive at the relation
\begin{equation}  \label{eq:CDGConsistencyCDG}
\begin{split}
\innerprod{c_h}{\MicroDens  \ddt{\PF}_h}{\RefDom}   &  + \innerprod{c_h}{\MicroDamp  \dt{\PF}_h}{\RefDom} + \innerprod{c_h }{\alpha_2 \Delta^2 \PF_h)}{\TildeRefDom} + \innerprod{ c_h }{ \alpha_1 \PFlap_h }{\RefDom} + \innerprod{c_h}{\alpha_0  \PF_h}{\RefDom}  \\
& + \innerprod{c_h}{ g^\prime(\PF_h) \mathcal{H}}{\RefDom} +\innerprod{\nabla c_h  \cdot \vect{N} }{\alpha_2 \PFlap_h}{\partial  \RefDom} - \innerprod{c_h }{\nabla \left(  \alpha_2 \PFlap + \alpha_1  \PF  \right) \cdot \vect{N}}{\ExtBoun_{R_2}^\PF} \\
 & + \innerprod{\frac{\beta_{s}}{\avrg{h_e}} \nabla c_h  \cdot \vect{N}-\Delta c_h}{ \alpha_2 \PFgrad_h \cdot \vect{N}}{\ExtBoun_{D_2}^\PF}  -\innerprod{c_h }{ \jump{\nabla \left(  \alpha_2 \PFlap + \alpha_1  \PF  \right)}}{\IntOnlyBoun} \\ 
 &+ \innerprod{\avrg{\nabla c_h}}{\alpha_2 \jump{\PFlap_h}}{\IntOnlyBoun} + \innerprod{\frac{\beta_{s}}{\avrg{h_e}} \jump{\nabla  c_h} - \avrg{\Delta c_h}}{ \alpha_2 \jump{\PFgrad_h}}{\IntOnlyBoun}  = 0 , \quad \forall c_h \in V_{\PF,1}^h 
\end{split}
\end{equation}
The above condition leads to the weak imposition of the governing equation and continuity requirements as follows
\begin{equation}
\left . \begin{aligned}
 \MicroDens  \ddt{\PF}_h +  \MicroDamp  \dt{\PF}_h + \alpha_2 \Delta^2 \PF_h + \alpha_1 \PFlap_h + \alpha_0  \PF_h + g^\prime(\PF_h) \mathcal{H} &= 0 \text{ in } \RefDom \\
\PFlap_h &= 0 \text{ on } \partial  \RefDom \\
\nabla \left(  \alpha_2 \PFlap + \alpha_1  \PF  \right) \cdot \vect{N}  &= 0  \text{ on } \ExtBoun_{R_2}^\PF \\
\PFgrad_h \cdot \vect{N} &= 0 \text{ on } \ExtBoun_{D_2}^\PF  \\
\jump{\nabla \left(  \alpha_2 \PFlap + \alpha_1  \PF  \right)} & = 0  \text{ on } \IntOnlyBoun \\
\jump{\PFlap_h} & = 0   \text{ on } \IntOnlyBoun  \\
\jump{\PFgrad_h}  & = 0 \text{ on } \IntOnlyBoun  \\
\end{aligned} \right\}
\end{equation}
It is noteworthy that the strong form of the phase-field formulation is imposed along with the corresponding boundary and regularity conditions. This procedure highlights the importance of jump and average terms for deriving the strong form of the fourth-order phase-field equation. 

\clearpage
%%%%%%\section*{References}

%\bibliography{bReferences}

\end{document}